\documentclass{amsart}
\usepackage[left=3.25cm, right=3.25cm, top=3cm, bottom=3cm]{geometry}

\usepackage{lipsum}
\usepackage{amsfonts}
\usepackage{graphicx}
\usepackage{float}
\usepackage{algorithm}
\usepackage{algorithmic}
\usepackage{multirow,makecell}
\usepackage{bm,amsfonts,booktabs,amssymb}
\usepackage{amsmath}
\makeatletter
\g@addto@macro\normalsize{%
  \setlength\abovedisplayskip{7pt}
  \setlength\belowdisplayskip{7pt}
  \setlength\abovedisplayshortskip{8pt}
  \setlength\belowdisplayshortskip{8pt}
}
\makeatother
\usepackage{mathtools}
\usepackage{enumitem}
\usepackage{subcaption}
\usepackage[font={small,it}]{caption}
\usepackage{comment}
\usepackage{nicematrix}
\usepackage[hidelinks]{hyperref}
\usepackage{cleveref}
\usepackage{autonum}
\usepackage{longtable}
\usepackage{tikz}
\usetikzlibrary{positioning}
\usetikzlibrary{arrows}

\ifpdf
  \DeclareGraphicsExtensions{.eps,.pdf,.png,.jpg}
\else
  \DeclareGraphicsExtensions{.eps}
\fi

\newcommand{\retainlabel}[1]{\label{#1}\sbox0{\ref{#1}}}

\newcommand*{\vertbar}{\rule[-1ex]{0.5pt}{2.5ex}}

\newcolumntype{P}[1]{>{\centering\arraybackslash}p{#1}}
\newcolumntype{M}[1]{>{\centering\arraybackslash}m{#1}}

\newcommand{\revadd}[1]{#1}
\newcommand{\revdel}[1]{\textcolor{red}{}}

\newtheorem{assump}{Assumption}
\newtheorem{dfn}{Definition}[section]
\newtheorem{thm}{Theorem}[section]
\newtheorem{cor}{Corollary}[section]
\newtheorem{lem}{Lemma}[section]
\newtheorem{prop}{Proposition}[section]
\newtheorem{rmk}{Remark}[section]

\newcommand{\proofoffirst}[1]{\noindent \underline{\textbf{Proof of #1.}}}
\newcommand{\proofof}[1]{\smallskip \noindent \underline{\textbf{Proof of #1}}.}

\newcommand{\argmin}{\mathop{\mathrm{argmin}}}

\newcommand{\Tr}{{\mathrm{Tr}}}

\newcommand{\EXP}{{\mathrm{Exp}}}

\newcommand{\grad}{{\mathrm{grad}}}
\newcommand{\vecz}{{\mathrm{vec}}}
\newcommand{\Hess}{{\mathrm{Hess}}}
\newcommand{\blockdiag}{{\mathrm{block\text{-}diag}}}
\newcommand{\Skew}{{\mathrm{Skew}}}
\newcommand{\Sym}{{\mathrm{Sym}}}
\newcommand{\rank}{{\mathrm{rank}}}

\newif\iftodos

\newcommand{\citea}[2][]{\cite{#2}}
\newcommand{\citeb}[2][]{\cite[#1]{#2}}

\makeatletter
\newcommand*{\addFileDependency}[1]{
\typeout{(#1)}
\@addtofilelist{#1}
\IfFileExists{#1}{}{\typeout{No file #1.}}
}\makeatother

\usepackage{amsopn}
\DeclareMathOperator{\diag}{diag}

\everypar{\looseness=-1}

\ifpdf
\hypersetup{
  pdftitle={Non-degenerate Rigid Alignment in a Patch Framework},
  pdfauthor={D. Kohli, G. Mishne and A. Cloninger}
}
\fi

\begin{document}
\renewcommand{\thefootnote}{\fnsymbol{footnote}}
\footnotetext[1]{Department of Mathematics, UC San Diego (dhkohli@ucsd.edu, acloninger@ucsd.edu)}
\footnotetext[2]{Halicio\u{g}lu Data Science Institute, UC San Diego (gmishne@ucsd.edu).}

\title[Non-degenerate Rigid Alignment in a Patch Framework]{Non-degenerate Rigid Alignment in a Patch Framework}
\author[D. Kohli, G. Mishne, A. Cloninger]{Dhruv Kohli${}^{\ast}$, Gal Mishne${}^\dagger$, Alexander Cloninger${}^{\ast, \dagger}$}

\renewcommand{\thefootnote}{\arabic{footnote}}

\begin{abstract}%
Given a set of overlapping local views (patches) of a dataset, we consider the problem of finding a rigid alignment of the views that minimizes a $2$-norm based alignment error. In general, the views are noisy and a perfect alignment may not exist. In this work, we characterize the non-degeneracy of an alignment in the noisy setting based on the kernel and positivity of a certain matrix. This leads to a polynomial time algorithm for testing the non-degeneracy of a given alignment. Subsequently, we focus on Riemannian gradient descent for minimizing the alignment error, providing a sufficient condition on an alignment for the algorithm to converge (locally) linearly to it. \revadd{Additionally, we provide an exact recovery and noise stability analysis of the algorithm}. In the case of noiseless views, a perfect alignment exists, resulting in a realization of the points that respects the geometry of the views. Under a mild condition on the views, we show that a non-degenerate perfect alignment \revadd{characterizes the infinitesimally rigidity of a realization, and thus the local rigidity of a generic realization}. By specializing the non-degeneracy conditions to the noiseless case, we derive necessary and sufficient conditions on the overlapping structure of the views for \revadd{a perfect alignment to be non-degenerate and equivalently, for the resulting realization to be infinitesimally rigid}. Similar results are also derived regarding the uniqueness of a perfect alignment and global rigidity.
\end{abstract}


\maketitle
\section{Introduction}
\label{sec:intro}
There exist a plethora of problems \citea{fang2012using,cucuringu2012sensor,sharp2004multiview,williams2000simultaneous} which involve the task of aligning local views (also known as registration of point clouds) so as to obtain a global view of the data that respects the geometry of the local views. Although it is not uncommon for the correspondence between the points and the views to not be known apriori \citea{szeliski2007image}, nevertheless, we assume that the overlapping structure of the local views is available. A suitable example of our setup is rendered by the bottom-up manifold learning techniques, LTSA \citea{zhang2004principal}, LDLE \citea{kohli2021ldle}, improved MDS \citea{shang2004improved} etc., which first construct low distortion local views of high dimensional data into lower dimension followed by the alignment of the views to obtain a low dimensional global embedding of the data. Here we focus on the alignment step.

The correspondence between the local views and the points is captured by a bipartite graph, which when combined with the local coordinates of the points due to the views containing them, form a \textit{patch framework} $\Theta$ \citea{chaudhury2015global,gortler2010affine}. The alignment of the views amounts to finding a suitable transformation for each view so that the local coordinates of a point in its transformed views are close to each other. While LTSA seeks affine transformations, LDLE seeks rigid ones for alignment. 

In the case where the views undergo rigid transformation, the problem of aligning views can be posed as the minimization of a $2$-norm based alignment error, given by a quadratic $F$ over the product of orthogonal groups $\mathbb{O}(d)^m$ (see Section~\ref{sec:setup}).  With our definition of the alignment error, we identify an ``alignment" of the local views by an element of $\mathbb{O}(d)^m$. Note that rotating/reflecting each view by the same amount does not affect the alignment error. For a given alignment $\mathbf{S} \in \mathbb{O}(d)^m$, this translates to the alignment error being the same due to $\mathbf{S}\mathbf{Q}$ for all $\mathbf{Q} \in \mathbb{O}(d)$ i.e. $F(\mathbf{S}) = F(\mathbf{S}\mathbf{Q})$. In this sense, every alignment is degenerate and every optimal alignment is non-unique. With a slight abuse of convention, \textit{we define a non-degenerate alignment to be a local minimum of $F$ which is non-degenerate up to an orthogonal transformation}. \textit{Similarly, we say that an optimal alignment (a global minimum of $F$) is unique if every other optimal alignment can be obtained by an orthogonal transformation of it}. This brings us to our first contribution.

\textit{In Section~\ref{sec:non_deg}, we derive a characterization of non-degenerate alignment of (possibly noisy) local views that can be tested in polynomial time. \revadd{Additionally, we derive a bound on the size of the neighborhood around a non-degenerate alignment in which the Hessian of $\widetilde{F}$, the function induced by $F$ on a quotient manifold, is positive definite.}
}

Given an alignment $\mathbf{S} \in \mathbf{O}(d)^m$ of the local views, a consensus representation $\Theta(\mathbf{S})$ of the points can be obtained by averaging the local coordinates of the points due to the (rigidly transformed) views containing them. In the noiseless setting where the local views are clean measurements of the data (obtained by applying an unknown rigid transformation to a subset of data points), a perfect alignment of views is possible. Equivalently, when the views are noiseless, a value of zero for $F$ is attainable, and an $\mathbf{S}$ that achieves it is called a ``perfect alignment". Clearly, a perfect alignment is an optimal one, while the converse may not hold. To be consistent with previous works \citea{gortler2010affine,hendrickson1992conditions}, the consensus representation of the points $\Theta(\mathbf{S})$ due to a perfect alignment $\mathbf{S}$ of the views is called a realization of the framework.

An understanding of \revadd{affine, global, local and infinitesimal rigidity} \citea{toth2017handbook,gortler2010affine} of a realization $\Theta(\mathbf{S})$ has importance in several areas such as molecular dynamics \citea{lee2008geometric,clementi2000topological,cucuringu2012eigenvector} or sensor network localization \citea{zhu2010universal,zhang2010rigid}. Under a mild assumption on the structure of the local views, \citea{chaudhury2015global,zha2009spectral,gortler2010affine} characterized the affine rigidity of a realization by the rank of a certain matrix derived from the framework $\Theta$.\revdel{It was shown in \citea{saxe1979embeddability} that deriving a similar characterization of global rigidity is NP-Hard. Nevertheless a characterization of the local and infinitesimal rigidity is useful from an algorithmic standpoint as we show in this work.} Furthermore, necessary and sufficient conditions on the overlapping structure of the local views for affine rigidity were derived in \citea{zha2009spectral}. Similar results in the context of infinitesimal, local and global rigidity form our second set of contributions.

\revadd{\textit{In Section~\ref{sec:noiseless_non_deg_results}, under a mild assumption on the structure of the local views, we show that a non-degenerate perfect alignment $\mathbf{S}$ characterizes the infinitesimal rigidity of the resulting realization $\Theta(\mathbf{S})$, and also its local rigidity if the realization is generic. By specializing the characterization of a non-degenerate alignment to noiseless setting, we derive necessary and sufficient conditions on the overlapping structure of the noiseless local views for a perfect alignment to be non-degenerate, and thus, for the resulting realization to be infinitesimally rigid. Similar results are also derived for the uniqueness of a perfect alignment, and equivalently for global rigidity of a realization.}}

Several algorithms exist to obtain an approximate solution of the alignment problem at hand. A few important ones are: semidefinite programming (SDP) \citea{chaudhury2015global,bandeira2014multireference}, spectral relaxation (SPEC) \citea{chaudhury2015global,bandeira2014multireference}, Procrustes analysis (PROC) \citea{kohli2021ldle,proc_algo_1}, generalized power method (GPM) \citea{ling2021generalized} and Riemannian gradient descent (RGD) \citea{krishnan2005global}. It is common to obtain stability and convergence guarantees of an algorithm under rigidity constraints over the framework. For example, \citea{chaudhury2015global} derived stability guarantees on the SPEC and SDP solutions under affine rigidity constraints. The authors of \citea{ling2021generalized,ling2021near} derived stability guarantees and showed global linear convergence of GPM with SPEC initialization under the setting where each view is affinely non-degenerate and contains all the points. It is easy to deduce that this framework also exhibits \revadd{affine} rigidity. Moreover, a recent study \cite{zhu2023rotation} employs a quotient manifold approach to establish convergence guarantees for RGD for a related problem of rotation group synchronization. Our last contribution is along similar lines.

\textit{\revadd{In Section~\ref{sec:convergence}, we establish convergence guarantees for RGD on a quotient manifold under non-degeneracy constraint, which we show (in Section~\ref{sec:noiseless_non_deg_results}), is less restrictive than the affine rigidity constraint. Specifically, we show that RGD converges locally linearly to a non-degenerate alignment of views, and derive the radius and rate of convergence via convexity argument \citea{udriste2013convex}. We also provide an exact recovery and a noise stability analysis of RGD when initialized with the solution of SPEC \cite{chaudhury2015global}. Through a simulation, we demonstrate that the alignment error decreases as RGD refines the initial spectral alignment, with a decay that is consistent with a linear order of convergence.}}

\section{A Quadratic Program in Orthogonal Groups for Aligning Views}
\label{sec:setup}
In this section we borrow and build upon the patch framework setup described in \citea{chaudhury2015global}, elucidate the structure of the objects underlying a framework, define an alignment, optimal alignment and perfect alignment of views, and the realization of a framework due to a perfect alignment. 

Suppose that a sequence of $m$ point clouds in $\mathbb{R}^d$ is available where each cloud represents a local view of the dataset $(\mathbf{x}_k)_1^n$.  Let $\Gamma$ be a graph of $m+n$ vertices, where the $k$th vertex represents the $k$th point for $k \in [1,n]$ and the $(n+i)$th vertex represents the $i$th view for $i \in [1,m]$. An edge $(k,i) \in E(\Gamma)$ means that the $k$th point has a local representation due to $i$th view, given by $\mathbf{x}_{k,i} \in \mathbb{R}^d$. In particular, $\Gamma$ is bipartite. The tuple $\Theta = (\Gamma, (\mathbf{x}_{k,i}))$ is called the patch framework.

Given a patch framework $\Theta$, the task is to align the overlapping views, precisely, to find orthogonal matrices $(\mathbf{S}_i)_{1}^{m} \subseteq \mathbb{O}(d)$ and translation vectors $(\mathbf{t}_i)_{1}^{m} \subseteq \mathbb{R}^{d}$ so that the local representations of the $k$th point across the rigidly transformed views are close in the $2$-norm. This naturally leads to the following problem,
\begin{equation}
    \mathcal{A}_0 \coloneqq \textstyle\min_{\substack{(\mathbf{S}_i)_1^m \subseteq \mathbb{O}(d)\\(\mathbf{t}_i)_1^m \subseteq \mathbb{R}^d}}\textstyle\sum_{\substack{(k,i)\in E(\Gamma)\\(k,j)\in E(\Gamma)}}\left\|(\mathbf{S}_i^T\mathbf{x}_{k,i}+\mathbf{t}_i)-(\mathbf{S}_j^T\mathbf{x}_{k,j}+\mathbf{t}_j)\right\|_2^2. \label{eq:A_0}
\end{equation}
An equivalent problem is to find $(\mathbf{x}_k)_1^n$, $(\mathbf{S}_i)_{1}^{m}$ and $(\mathbf{t}_i)_{1}^{m}$ that minimize
\begin{equation}
    \mathcal{A}_1 \coloneqq \textstyle\min_{\substack{(\mathbf{S}_i)_1^m \subseteq \mathbb{O}(d),\\ (\mathbf{t}_i)_1^m, (\mathbf{x}_k)_1^n \subseteq \mathbb{R}^d}}\textstyle\sum_{\substack{(k,i)\in E(\Gamma)}}\left\|\mathbf{x}_k-(\mathbf{S}_i^T\mathbf{x}_{k,i}+\mathbf{t}_i)\right\|_2^2. \label{eq:A_1}
\end{equation}
Here, $\mathcal{A}_0$ and $\mathcal{A}_1$ are the optimal alignment errors. 
\begin{prop}
\label{prop:A_0_2_A_1}
$\mathcal{A}_0 = 2\mathcal{A}_1$.
\end{prop}

Note that translating the optimal $\mathbf{x}_k$ and $\mathbf{t}_i$ by the same amount does not change the objective and still leads to an optimal solution. To avoid that, we add a centering constraint $\sum_1^n \mathbf{x}_k = 0$. A more concise representation of $\mathcal{A}_1$ is then derived as follows. Define $ \mathbf{H} \coloneqq [\mathbf{x}_1, \mathbf{x}_2, \ldots, \mathbf{x}_n, \mathbf{t}_1, \mathbf{t}_2, \ldots, \mathbf{t}_m] \in \mathbb{R}^{d \times (n+m)}$ and
$\mathbf{e}_{ki} \coloneqq \mathbf{e}^{n+m}_k - \mathbf{e}^{n+m}_{n+i}$, then $\mathbf{x}_k - \mathbf{t}_i = \mathbf{H}\mathbf{e}_{ki}$. Let $\mathbf{S} = [\mathbf{S}_i]_1^m \in \mathbb{O}(d)^m \subseteq \mathbb{R}^{md \times d}$. Then
\begin{align}
    \mathcal{A}_1 &= \textstyle\min_{\substack{\mathbf{S} \in \mathbb{O}(d)^m\\\mathbf{H}\mathbf{1}^{n+m}_n = 0}}\textstyle\sum_{(k,i) \in E(\Gamma)}\left\|\mathbf{H}\mathbf{e}_{ki}-\mathbf{S}^T(\mathbf{e}^m_i \otimes \mathbf{I}_d)\mathbf{x}_{k,i}\right\|_2^2\\
    &= \textstyle\min_{\substack{\mathbf{S} \in \mathbb{O}(d)^m\\\mathbf{H}\mathbf{1}^{n+m}_n = 0}}\Tr\left(\begin{bmatrix}\mathbf{H} & \mathbf{S}^T\end{bmatrix}\begin{bmatrix}\boldsymbol{\mathcal{L}}_{\Gamma} & -\mathbf{B}^T\\-\mathbf{B} & \mathbf{D}\end{bmatrix}\begin{bmatrix}\mathbf{H}^T\\\mathbf{S}\end{bmatrix}\right) \text{, where} \label{eq:A_1_}
\end{align}
\begin{align}
    \boldsymbol{\mathcal{L}}_{\Gamma} &= \textstyle\sum_{(k,i) \in E(\Gamma)}\mathbf{e}_{ki}\mathbf{e}_{ki}^T\retainlabel{eq:L_Gamma}\\
    \mathbf{B} &= \textstyle\sum_{(k,i) \in E(\Gamma)} (\mathbf{e}^m_i \otimes \mathbf{I}_d)\mathbf{x}_{k,i}\mathbf{e}_{ki}^T \label{eq:B}\\
    \mathbf{D} &= \textstyle\sum_{(k,i) \in E(\Gamma)} (\mathbf{e}^m_i \otimes \mathbf{I}_d)\mathbf{x}_{k,i}\mathbf{x}_{k,i}^T(\mathbf{e}^m_i \otimes \mathbf{I}_d)^T. \retainlabel{eq:D}
\end{align}
\begin{rmk}
\label{rmk:L0DB}
Note that $\boldsymbol{\mathcal{L}}_{\Gamma}$ is the combinatorial Laplacian of the graph $\Gamma$. Due to the bipartite structure of $\Gamma$, $\boldsymbol{\mathcal{L}}_{\Gamma}^\dagger$ can be computed in $O(nm^2)$ \citea{HO2005917}, against the complexity of $O((m+n)^3)$ for the general case.  The matrix $\mathbf{D}$ is a block diagonal matrix where the $i$th block is $\sum_{(k,i) \in E(\Gamma)}\mathbf{x}_{k,i}\mathbf{x}_{k,i}^T$. The matrix $\mathbf{B} = [\mathbf{B}_i]_1^m$ is a vertical stack of $m$ matrices, one for each view, where each $\mathbf{B}_i \in \mathbb{R}^{d \times (n+m)}$. For a fixed $i$, $\mathbf{B}_i(:,n+i) = -\sum_{(k,i)\in E}\mathbf{x}_{k,i}$, if $(k,i) \in E(\Gamma)$ then $\mathbf{B}_i(:,k) = \mathbf{x}_{k,i}$ otherwise zero. Thus, $\mathbf{1}_{n+m} \in \ker(\mathbf{B}_i)$. Finally, the $j$th row of $\mathbf{B}_i$ contains information about the $j$th coordinates of the points in the $i$th view. We will use the following result later,
\end{rmk}
\begin{prop}
\label{prop:kerB}
$\ker (\boldsymbol{\mathcal{L}}_{\Gamma}) \subseteq \ker (\mathbf{B})$.
\end{prop}
\begin{assump}
\label{assump:connected_gamma}
The dimension of the kernel of $\boldsymbol{\mathcal{L}}_{\Gamma}$ equals the number of connected components in $\Gamma$. To keep the subsequent calculations simple, we assume that $\Gamma$ is connected and thus $\ker (\boldsymbol{\mathcal{L}}_{\Gamma})$ is the span of a single vector $\mathbf{1}_{n+m}$.
\end{assump}

It is clear from Eq.~(\ref{eq:A_1_}) that the the optimal $\mathbf{H}$ satisfies $\mathbf{H}^*\boldsymbol{\mathcal{L}}_{\Gamma} = \mathbf{S}^T\mathbf{B}$. Using the fact that $\boldsymbol{\mathcal{L}}_{\Gamma}^\dagger \boldsymbol{\mathcal{L}}_{\Gamma} = \boldsymbol{\mathcal{L}}_{\Gamma}\boldsymbol{\mathcal{L}}_{\Gamma}^\dagger = \mathbf{I}_{n+m} - (n+m)^{-1}\mathbf{1}_{n+m}\mathbf{1}_{n+m}^T$, the solution $\mathbf{H}^*$ is of the form $\mathbf{S}^T\mathbf{B}\boldsymbol{\mathcal{L}}_{\Gamma}^\dagger - \mathbf{h}\mathbf{1}_{n+m}^T$ for some translation vector $\mathbf{h} \in \mathbb{R}^d$. Since $\mathbf{H}\mathbf{1}^{n+m}_n = 0$, the optimal value of $\mathbf{h}$ is $\mathbf{S}^T\mathbf{B}\boldsymbol{\mathcal{L}}_{\Gamma}^\dagger \mathbf{1}^{n+m}_n/n$. Substituting back,
\begin{equation}
    \mathbf{H}^* = \mathbf{S}^T\mathbf{B}\boldsymbol{\mathcal{L}}_{\Gamma}^\dagger \left(\mathbf{I}_{n+m} - n^{-1}\mathbf{1}^{n+m}_n\mathbf{1}_{n+m}^T\right). \label{eq:opt_Z}
\end{equation}
Substituting Eq.~(\ref{eq:opt_Z}) back in Eq.~(\ref{eq:A_1_}), the problem reduces to
\begin{equation}
    \mathcal{A}_1 = \min_{\mathbf{S} \in \mathbb{O}(d)^m} F(\mathbf{S}) = \min_{\mathbf{S} \in \mathbb{O}(d)^m} \Tr(\mathbf{C}\mathbf{S}\mathbf{S}^T) \text{ where } \mathbf{C} = \mathbf{D} - \mathbf{B}\boldsymbol{\mathcal{L}}_{\Gamma}^{\dagger}\mathbf{B}^T. \label{eq:GPOP}
\end{equation}
Here $\mathbf{C}$ is named the \textit{patch-stress matrix} and is positive semidefinite \citea{chaudhury2015global}.
\begin{dfn}
The objective $F$ depends only on $\mathbf{S}$, so we define an \textit{alignment} of local views as an element of $\mathbb{O}(d)^m$. If $\mathbf{S}$ is a global minimum of $F$, it's called an optimal alignment. If $F(\mathbf{S}) = 0$, it's called a \textit{perfect alignment}. Since $\mathbf{C} \succeq 0$, every perfect alignment is optimal, but not every optimal alignment is perfect.
\end{dfn}

\begin{dfn}
\label{def:realization}
The consensus representation of the framework $\Theta$ due to an alignment $\mathbf{S}$ is given by $\Theta(\mathbf{S}) \coloneqq \mathbf{H}^*(:,1:n) = \mathbf{S}^T\mathbf{B}\boldsymbol{\mathcal{L}}_{\Gamma}^\dagger(:,1:n)\left(\mathbf{I}_n - n^{-1}\mathbf{1}_{n}\mathbf{1}_{n}^T\right)$. The one due to a perfect alignment is called a realization of the framework \citea{gortler2010affine}.
\end{dfn}

\section{Non-degeneracy and Uniqueness in the General Setting}
\label{sec:non_deg}
In this section, we derive the Hessian of $\widetilde{F}$, the function induced by $F$ on a certain quotient space $\mathbb{O}(d)^m/_{\sim}$. We refer the reader to \citeb[Chapter 3 and 5]{absil2009optimization} for the definitions of differential of a mapping, metric, gradient, connection and Hessian in the context of Riemannian manifolds, as well as to \citea{zhu2023rotation,luo2022nonconvex,dong2022analysis,zhao2015riemannian}, which also adopt a quotient manifold approach in their analyses. We also obtain the equations governing the non-singularity and positivity of the Hessian and consequently, a characterization of a non-degenerate alignment in the general (noisy) setting. \revadd{Moreover, we identify the vicinity of a non-degenerate alignment in which the Hessian is positive definite.}

\subsection{Preliminaries}
\label{subsec:prelims}
Recall that the problem under consideration is the minimization of $F(\mathbf{S}) = \Tr(\mathbf{C}\mathbf{S}\mathbf{S}^T)$ over $\mathbf{S} \in \mathbb{O}(d)^m$ where $\mathbf{C} \succeq 0$ is the patch-stress matrix defined in Eq.~(\ref{eq:GPOP}). Note that the objective is invariant to the action of $\mathbb{O}(d)$ i.e. for any $\mathbf{Q} \in \mathbb{O}(d)$, $F(\mathbf{S}) = F(\mathbf{S}\mathbf{Q})$.
\begin{assump}
If $\Gamma$ has $K$ connected components, then the objective is invariant to the action of $\mathbb{O}(d)^K$ on $\mathbb{O}(d)^m = \prod_1^K \mathbb{O}(d)^{m_j}$ where $\textstyle\sum_1^K m_j = m$ and where each $\mathbb{O}(d)$ acts independently on $\mathbb{O}(d)^{m_j}$ for $j \in [1,K]$. To keep the computations clean, we assume that the bipartite graph $\Gamma$ is connected (as in Assumption~\ref{assump:connected_gamma}) throughout.
\end{assump}
Subsequently, we define an equivalence relation on $\mathbb{O}(d)^m$; $\mathbf{S}_1 \sim \mathbf{S}_2$ iff $\mathbf{S}_1 = \mathbf{S}_2\mathbf{Q}$ for some $\mathbf{Q} \in \mathbb{O}(d)$. Given $\mathbf{S} \in \mathbb{O}(d)^m$, its equivalence class is $[\mathbf{S}] = \left\{\mathbf{S}\mathbf{Q}: \mathbf{Q} \in \mathbb{O}(d)\right\}$. Clearly, there exists a bijection between $\mathbb{O}(d)^m/_{\sim}$ and $\mathbb{O}(d)^{m-1}$, thus an element of $\mathbb{O}(d)^m/_{\sim}$ will be identified with an element of $\mathbb{O}(d)^{m-1}$. Define the projection, 
\begin{align}
    \pi: \mathbb{O}(d)^m &\mapsto \mathbb{O}(d)^m/_{\sim}\\
    \pi(\mathbf{S}_{1:m}) &= \mathbf{S}_{2:m}\mathbf{S}_1^T. \label{eq:pi}
\end{align}
Let $\widetilde{\mathbf{S}} \in \mathbb{O}(d)^m/_{\sim}$, then
\begin{equation}
    \pi^{-1}(\widetilde{\mathbf{S}}) = \left\{\begin{bmatrix}
    \mathbf{Q}\\\widetilde{\mathbf{S}}\mathbf{Q}
    \end{bmatrix}: \mathbf{Q} \in \mathbb{O}(d)\right\} = \{\mathbf{S} \in \mathbb{O}(d)^m: \mathbf{S}_{i+1}\mathbf{S}_1^T = \widetilde{\mathbf{S}}_{i}, i \in [1,m-1]\}. \label{eq:pi_inv_wtS}
\end{equation}

The Riemannian metric $g$ on $\mathbb{O}(d)^m$ is the canonical one given by
\begin{equation}
g(\mathbf{Z},\mathbf{W}) \coloneqq \Tr(\mathbf{Z}^T\mathbf{W}) = \textstyle\sum_1^m \Tr(\mathbf{Z}_i^T\mathbf{W}_i) \text{ where } \mathbf{Z},\mathbf{W} \in T_{\mathbf{S}}\mathbb{O}(d)^m \subseteq \mathbb{R}^{md \times d}. \label{eq:g_Z_W}
\end{equation}
By a simple extension of the $m=1$ case \citea{absil2009optimization}, it is easy to deduce the following result.
\begin{prop}
\label{prop:T_SOdm}
For $\mathbf{S} \in \pi^{-1}(\widetilde{\mathbf{S}})$, the tangent space to $\mathbb{O}(d)^m$ at $\mathbf{S}$ is given by
\begin{equation}
    T_{\mathbf{S}}\mathbb{O}(d)^m = \{[\mathbf{S}_i\boldsymbol{\Omega}_i]_1^m: \boldsymbol{\Omega}_i \in \Skew(d)\}. \label{eq:T_SOdm}
\end{equation}
The orthogonal projection of $\boldsymbol{\xi} = [\boldsymbol{\xi}_i]_1^m$, where $\boldsymbol{\xi}_i \in \mathbb{R}^{d \times d}$, onto $T_{\mathbf{S}}\mathbb{O}(d)^m$ is
\begin{align}
    P_{\mathbf{S}}\left(\boldsymbol{\xi}\right) = \argmin_{[\mathbf{S}_i\boldsymbol{\Omega}_i]_1^m,\boldsymbol{\Omega}_i \in \Skew(d)} \textstyle\sum_1^m\left\|\boldsymbol{\xi}_i - \mathbf{S}_i\boldsymbol{\Omega}_i\right\|_F^2 = [\mathbf{S}_i\Skew(\mathbf{S}_i^T\boldsymbol{\xi}_i)]_1^m. \label{eq:P_S_xi}
\end{align}
\end{prop}

Then, $\pi^{-1}(\widetilde{\mathbf{S}})$ admits a tangent space at $\mathbf{S} \in \pi^{-1}(\widetilde{\mathbf{S}})$ called the vertical space $\mathcal{V}_{\mathbf{S}}$ at $\mathbf{S}$. The horizontal space $\mathcal{H}_{\mathbf{S}}$ at $\mathbf{S}$ is the subspace of $T_{\mathbf{S}}\mathbb{O}(d)^m$ that is the orthogonal complement to the vertical space $\mathcal{V}_{\mathbf{S}}$.
\begin{prop}
\label{prop:V_S_H_S}
The vertical space $\mathcal{V}_{\mathbf{S}}$ at $\mathbf{S} \in \pi^{-1}(\widetilde{\mathbf{S}})$ is
 $$   \mathcal{V}_{\mathbf{S}} = \{\mathbf{S}\boldsymbol{\Omega}: \boldsymbol{\Omega} \in \Skew(d)\}.$$
The orthogonal projection of $\mathbf{Z} = [\mathbf{S}_i\boldsymbol{\Omega}_i]_1^m \in T_{\mathbf{S}}\mathbb{O}(d)^m$ onto $\mathcal{V}_{\mathbf{S}}$ is
\begin{equation}
    P^{v}_{\mathbf{S}}([\mathbf{S}_i\boldsymbol{\Omega}_i]_1^m) = \left[\mathbf{S}_i\argmin_{\boldsymbol{\Omega}\in \Skew(d)} \textstyle\sum_1^m\left\|\mathbf{S}_j(\boldsymbol{\Omega}_j-\boldsymbol{\Omega})\right\|_F^2\right]_1^m = \left[\mathbf{S}_i\left(m^{-1}\textstyle\sum_1^m\boldsymbol{\Omega}_i\right)\right]_1^m. \label{eq:P^v_S}
\end{equation}
The horizontal space at $\mathbf{S} \in \pi^{-1}(\widetilde{\mathbf{S}})$ is
$$\mathcal{H}_{\mathbf{S}} = \left\{[\mathbf{S}_i\boldsymbol{\Omega}_i]_1^m, \boldsymbol{\Omega}_i \in \Skew(d), \textstyle\sum_1^m \boldsymbol{\Omega}_i= 0\right\}.$$
The orthogonal projection of $\mathbf{Z} = [\mathbf{S}_i\boldsymbol{\Omega}_i]_1^m \in T_{\mathbf{S}}\mathbb{O}(d)^m$ to $\mathcal{H}_{\mathbf{S}}$ is
\begin{equation}
    P^{h}_{\mathbf{S}}([\mathbf{S}_i\boldsymbol{\Omega}_i]_1^m) = [\mathbf{S}_i\boldsymbol{\Omega}_i]_1^m - P^{v}_{\mathbf{S}}([\mathbf{S}_i\boldsymbol{\Omega}_i]_1^m) = \left[\mathbf{S}_i\left(\boldsymbol{\Omega}_i-m^{-1}\textstyle\sum_1^m\boldsymbol{\Omega}_i\right)\right]_1^m.\label{eq:P^h_S}
\end{equation}
\end{prop}

Note that $T_\mathbf{S}\mathbb{O}(d)^m$ is a vector space of dimension $md(d-1)/2$ and $\mathcal{V}_{\mathbf{S}}$ forms a $d(d-1)/2$ dimensional subspace of $T_{\mathbf{S}}\mathbb{O}(d)^m$. The dimension of $\mathcal{H}_{\mathbf{S}}$ and $T_{\widetilde{\mathbf{S}}}\mathbb{O}(d)^m/_{\sim}$ is $(m-1)d(d-1)/2$. In particular, $\mathcal{H}_{\mathbf{S}}$ can be identified with $T_{\widetilde{\mathbf{S}}}\mathbb{O}(d)^m/_{\sim}$. Let $\widetilde{\mathbf{S}} \in \mathbb{O}(d)^{m}/_{\sim}$ and $\widetilde{\mathbf{Z}} \in T_{\widetilde{\mathbf{S}}}\mathbb{O}(d)^{m}/_{\sim}$. Then the horizontal lift of $\widetilde{\mathbf{Z}}$ at $\mathbf{S} \in \pi^{-1}(\widetilde{\mathbf{S}})$ is defined as $\overline{\widetilde{\mathbf{Z}}} \in \mathcal{H}_{\mathbf{S}}$ such that for each $i \in [1,m-1]$,
\begin{equation}
    D\pi[\mathbf{S}]\left(\overline{\widetilde{\mathbf{Z}}}\right)_i = \widetilde{\mathbf{Z}}_i. \label{eq:hlift_def}
\end{equation}

\begin{prop}
\label{prop:hlift_char}
Let $\widetilde{\mathbf{S}} \in \mathbb{O}(d)^{m}/_{\sim}$ and $\widetilde{\mathbf{Z}} \in T_{\widetilde{\mathbf{S}}}\mathbb{O}(d)^{m}/_{\sim}$. Let $\mathbf{Z}$ be the horizontal lift of $\widetilde{\mathbf{Z}}$ at $\mathbf{S} \in \pi^{-1}(\widetilde{\mathbf{S}})$. If $(\widetilde{\boldsymbol{\Omega}}_i)_1^{m-1} \subseteq \Skew(d)$ are such that $\widetilde{\mathbf{Z}}_i = \widetilde{\mathbf{S}}_i\widetilde{\boldsymbol{\Omega}}_i$, and $(\boldsymbol{\Omega}_i)_1^m \subseteq \Skew(d)$ are such that $\mathbf{Z}_i=\mathbf{S}_i\boldsymbol{\Omega}_i$ and $\textstyle\sum_1^m \boldsymbol{\Omega}_i = 0$, then
\begin{align}
    \boldsymbol{\Omega}_1 &= -m^{-1}\mathbf{S}_1^T\left(\textstyle\sum_1^{m-1}\widetilde{\boldsymbol{\Omega}}_i\right)\mathbf{S}_1 \label{eq:hlift1}\\
    \boldsymbol{\Omega}_{i+1} &= \mathbf{S}_1^T\widetilde{\boldsymbol{\Omega}}_i\mathbf{S}_1 + \boldsymbol{\Omega}_1 \text{ for all } i \in [1,m-1]. \label{eq:hlifti}
\end{align}
Moreover, the linear system above has full rank, and thus the horizontal lift $\mathbf{Z}$ of $\widetilde{\mathbf{Z}}$ at $\mathbf{S} \in \mathbb{O}(d)^m$ is a unique element of $\mathcal{H}_{\mathbf{S}}$.
\end{prop}

\begin{prop}
\label{prop:g_tilde}
Let $\widetilde{\mathbf{Z}},\widetilde{\mathbf{W}} \in T_{\widetilde{\mathbf{S}}}\mathbb{O}(d)^{m}/_{\sim}$ and $\mathbf{Z}, \mathbf{W} \in T_{\mathbf{S}}\mathbb{O}(d)^m$ be their horizontal lifts at $\mathbf{S} \in \pi^{-1}(\widetilde{\mathbf{S}})$. Then
 $   \widetilde{g}(\widetilde{\mathbf{Z}},\widetilde{\mathbf{W}}) \coloneqq g(\mathbf{Z}, \mathbf{W})$
defines a Riemannian metric on $\mathbb{O}(d)^m/_{\sim}$.
\end{prop}
We note that $\mathcal{H}_{\mathbf{S}}$ with the canonical metric $g$, is isometric to $T_{\widetilde{\mathbf{S}}}\mathbb{O}(d)^{m}/_{\sim}$ (equivalently $T_\mathbf{S}\mathbb{O}(d)^{m-1}$) when equipped with the above metric $\widetilde{g}$. \revadd{The following relation holds between the Frobenius norm of an element of $T_{\widetilde{\mathbf{S}}}\mathbb{O}(d)^{m}/_{\sim}$ and of its horizontal lift.}
\begin{prop}
\label{prop:hlift_frob_ineq}
\revadd{Let $\widetilde{\mathbf{Z}} \in T_{\widetilde{\mathbf{S}}}\mathbb{O}(d)^{m}/_{\sim}$ and $\mathbf{Z} \in T_{\mathbf{S}}\mathbb{O}(d)^m$ be the horizontal lift of $\mathbf{\widetilde{Z}}$ at $\mathbf{S} \in \pi^{-1}(\widetilde{\mathbf{S}})$. Then $\left\|\mathbf{Z}\right\|_F \leq \left\|\widetilde{\mathbf{Z}}\right\|_F \leq \sqrt{m+1}\left\|\mathbf{Z}\right\|_F$.}
\end{prop}

Now, coming back to the alignment error $F(\mathbf{S}) = \Tr(\mathbf{C}\mathbf{S}\mathbf{S}^T)$ defined on $\mathbb{O}(d)^m$. It induces the following function $\widetilde{F}(\widetilde{\mathbf{S}}) = \Tr\left(\mathbf{C}\begin{bsmallmatrix}
    \mathbf{I}_d\\\widetilde{\mathbf{S}}
    \end{bsmallmatrix}\begin{bsmallmatrix}
    \mathbf{I}_d\\\widetilde{\mathbf{S}}
    \end{bsmallmatrix}^T\right)$ on $\mathbb{O}(d)^{m}/_{\sim}$ (again identified with $\mathbb{O}(d)^{m-1}$).
In particular, $F = \widetilde{F} \circ \pi$ and $\widetilde{F} \circ \pi(\mathbf{S}_1) = \widetilde{F} \circ \pi(\mathbf{S}_2)$ if $\mathbf{S}_1 \sim \mathbf{S}_2$.

\begin{prop}
\label{prop:gradFS}
The horizontal lift of $\grad \widetilde{F}(\widetilde{\mathbf{S}})$ at $\mathbf{S} \in \pi^{-1}(\widetilde{\mathbf{S}})$ is
\begin{equation}
    \overline{\grad \widetilde{F}(\widetilde{\mathbf{S}})} = \grad F(\mathbf{S}) = [\mathbf{S}_i\boldsymbol{\Omega}_i]_1^{m} \label{eq:gradFS}
\end{equation}
where $\boldsymbol{\Omega}_i \coloneqq \mathbf{S}_i^T[\mathbf{C}\mathbf{S}]_i - [\mathbf{C}\mathbf{S}]_i^T\mathbf{S}_i \in \Skew(d)$ and $\textstyle\sum_1^m \boldsymbol{\Omega}_i = \mathbf{S}^T\mathbf{C}\mathbf{S} - \mathbf{S}^T\mathbf{C}\mathbf{S} = 0$ (which validates that $\overline{\grad \widetilde{F}(\widetilde{\mathbf{S}})}$ is indeed in $\mathcal{H}_{\mathbf{S}}$, see Proposition~\ref{prop:V_S_H_S}). Consequently, the set of critical points of $\widetilde{F}$ is given by $\widetilde{\mathcal{C}} = \{\widetilde{\mathbf{S}} \in \mathbb{O}(d)^m/_{\sim}: \grad \widetilde{F}(\widetilde{\mathbf{S}}) = 0\}$, equivalently,
\begin{equation}
    \widetilde{\mathcal{C}} = \{\widetilde{\mathbf{S}} \in \mathbb{O}(d)^m/_{\sim}: \mathbf{S}_i^T[\mathbf{C}\mathbf{S}]_i = [\mathbf{C}\mathbf{S}]_i^T\mathbf{S}_i, \text{ for all } i \in [1,m], \mathbf{S} \in \pi^{-1}(\widetilde{\mathbf{S}})\}, \label{eq:crit_pts}
\end{equation}
and that of $F$ is $\mathcal{C} = \{\mathbf{S} \in \mathbb{O}(d)^m: \grad F(\mathbf{S}) = 0\}$, equivalently,
\begin{equation}
     \mathcal{C} = \{\mathbf{S} \in \mathbb{O}(d)^m: \mathbf{S}_i^T[\mathbf{C}\mathbf{S}]_i = [\mathbf{C}\mathbf{S}]_i^T\mathbf{S}_i, \text{ for all } i \in [1,m]\}. \label{eq:crit_pts2}
\end{equation}
\end{prop}
From Eq.~(\ref{eq:crit_pts}) and (\ref{eq:crit_pts2}), it is easy to see that if $\mathbf{S} \in \mathcal{C}$ then $\pi(\mathbf{S}) \in \widetilde{C}$. Similarly, if $\widetilde{\mathbf{S}} \in \widetilde{C}$ then $\mathbf{S} \in \mathcal{C}$ for all $\mathbf{S} \in \pi^{-1}(\widetilde{\mathbf{S}})$.

\begin{prop}
\label{prop:DgradFSZ}
Let \revadd{$\widetilde{\mathbf{S}} \in \mathbb{O}(d)^{m}/_{\sim}$} then for every $\mathbf{S} \in \pi^{-1}(\widetilde{\mathbf{S}})$, $\mathbf{Z} \in T_\mathbf{S}\mathbb{O}(d)^m$,
$$D\grad F(\mathbf{S})[\mathbf{Z}] = \left[\mathbf{S}_i(\mathbf{S}_i^T[\mathbf{C}\mathbf{Z}]_i - [\mathbf{C}\mathbf{Z}]_i^T\mathbf{S}_i - [\mathbf{C}\mathbf{S}]_i^T\mathbf{Z}_i + \mathbf{Z}_i^T\mathbf{S}_i[\mathbf{C}\mathbf{S}]_i^T\mathbf{S}_i)\right]_1^m.$$
\revadd{Moreover, if $\widetilde{\mathbf{S}} \in \widetilde{\mathcal{C}}$ then} 
$$D\grad F(\mathbf{S})[\mathbf{Z}] = \left[\mathbf{S}_i(\mathbf{S}_i^T[\mathbf{C}\mathbf{Z}]_i - [\mathbf{C}\mathbf{Z}]_i^T\mathbf{S}_i - [\mathbf{C}\mathbf{S}]_i^T\mathbf{Z}_i + \mathbf{Z}_i^T[\mathbf{C}\mathbf{S}]_i)\right]_1^m.$$
\end{prop}

Let $\nabla$ be the Levi-Civita connection (also known as the Riemannian connection) on $\mathbb{O}(d)^m$ and $\widetilde{\nabla}$ be the induced connection on $\mathbb{O}(d)^m /_{\sim}$. Then, from Proposition~\ref{prop:DgradFSZ} and the definition of the Riemannian Hessian operator \citeb[Section 5.5]{absil2009optimization}, we obtain
\begin{prop}
\label{prop:HessFSZ}
Let \revdel{$\widetilde{\mathbf{S}} \in \widetilde{\mathcal{C}}$}\revadd{$\widetilde{\mathbf{S}} \in \mathbb{O}(d)^m/_{\sim}$} and $\widetilde{\mathbf{Z}} \in T_{\widetilde{\mathbf{S}}}\mathbb{O}(d)^{m}/_{\sim}$. Let $\mathbf{Z}$ be the horizontal lift of $\widetilde{\mathbf{Z}}$ at $\mathbf{S} \in \pi^{-1}(\widetilde{\mathbf{S}})$. Then the horizontal lift of $\Hess \widetilde{F}(\widetilde{\mathbf{S}})[\widetilde{\mathbf{Z}}]$ at $\mathbf{S}$ is
\begin{equation}
    \overline{\Hess \widetilde{F}(\widetilde{\mathbf{S}})[\widetilde{\mathbf{Z}}]} = [\mathbf{S}_i\widehat{\boldsymbol{\Omega}}_i]_1^m \label{eq:eq3}
\end{equation}
where \revadd{$\widehat{\boldsymbol{\Omega}}_i = \Skew(\boldsymbol{\xi}_i) - m^{-1}\sum_1^m\Skew(\boldsymbol{\xi}_i)$ and $$\boldsymbol{\xi}_i = \mathbf{S}_i^T[\mathbf{C}\mathbf{Z}]_i - [\mathbf{C}\mathbf{Z}]_i^T\mathbf{S}_i - [\mathbf{C}\mathbf{S}]_i^T\mathbf{Z}_i + \mathbf{Z}_i^T\mathbf{S}_i[\mathbf{C}\mathbf{S}]_i^T\mathbf{S}_i,$$ and if $\widetilde{\mathbf{S}} \in \widetilde{\mathcal{C}}$ then }
$$\widehat{\boldsymbol{\Omega}}_i = \mathbf{S}_i^T[\mathbf{C}\mathbf{Z}]_i - [\mathbf{C}\mathbf{Z}]_i^T\mathbf{S}_i - [\mathbf{C}\mathbf{S}]_i^T\mathbf{Z}_i + \mathbf{Z}_i^T[\mathbf{C}\mathbf{S}]_i.$$ 
\end{prop}

Now we obtain a compact representation for $\widehat{\boldsymbol{\Omega}}_i$. We first define certain matrices, then use them to obtain an expression for $\widehat{\boldsymbol{\Omega}}_i$ and then describe their structure. Recall that $\mathbf{C} = \mathbf{D} - \mathbf{B}\boldsymbol{\mathcal{L}}_{\Gamma}^\dagger \mathbf{B}^T$ (see Eq.~(\ref{eq:GPOP}) and Remark~\ref{rmk:L0DB}) and for convenience define
\begin{equation}
    \mathbf{B}(\mathbf{S}) \coloneqq \blockdiag((\mathbf{S}_i)_1^m)^T\ \mathbf{B}  = [\mathbf{S}_i^T\mathbf{B}_i]_1^m \label{eq:BofS}
\end{equation}
and
$$\mathbf{D}(\mathbf{S}) \coloneqq \blockdiag((\mathbf{S}_i)_1^m)^T\ \mathbf{D}\ \blockdiag((\mathbf{S}_i)_1^m) = \blockdiag((\mathbf{S}_i^T\mathbf{D}_{ii}\mathbf{S}_i)_1^m).$$ Using these matrices, we also define,
\begin{align}
    \mathbf{C}(\mathbf{S}) &\coloneqq \blockdiag((\mathbf{S}_i)_1^m)^T\ \mathbf{C}\ \blockdiag((\mathbf{S}_i)_1^m) = \mathbf{D}(\mathbf{S}) - \mathbf{B}(\mathbf{S})\boldsymbol{\mathcal{L}}_{\Gamma}^\dagger \mathbf{B}(\mathbf{S})^T\label{eq:C_of_S}\\
    \widehat{\mathbf{C}}(\mathbf{S}) &\coloneqq \blockdiag(([\mathbf{C}(\mathbf{S})\mathbf{I}_d^m]_i)_1^m) \label{eq:C_hat}\\
    \revadd{\mathbf{L}(\mathbf{S})} &\coloneqq \revadd{\mathbf{C}(\mathbf{S}) - \widehat{\mathbf{C}}(\mathbf{S})}. \label{eq:L_of_S}
\end{align}
\begin{rmk}
\label{rmk:StildeCtilde}
From the definition of $\mathcal{C}$ (see Eq.~(\ref{eq:crit_pts2})), $\mathbf{S} \in \mathcal{C}$ iff $\widehat{\mathbf{C}}(\mathbf{S}) \in \Sym(md)$ i.e. for each $i \in [1,m]$, $[\mathbf{C}(\mathbf{S})\mathbf{I}_d^m]_i = [\mathbf{C}(\mathbf{S})\mathbf{I}_d^m]_i^T$ or equivalently,
\begin{equation}
     \textstyle\sum_{j=1}^{m}\mathbf{C}(\mathbf{S})_{ij} = \mathbf{S}_i^T[\mathbf{C}\mathbf{S}]_i = [\mathbf{C}\mathbf{S}]_i^T\mathbf{S}_i = \textstyle\sum_{j=1}^{m}\mathbf{C}(\mathbf{S})_{ij}^T. \label{eq:eq2}
\end{equation}
\end{rmk}
\begin{rmk}
\label{rmk:C_S_structure}
Since $\blockdiag((\mathbf{S}_i)_1^m)$ is an orthogonal matrix, $\mathbf{C}(\mathbf{S})$ is unitarily equivalent to $\mathbf{C}$. Thus, $\mathbf{C}(\mathbf{S}) \in \Sym(md)$, $\mathbf{C}(\mathbf{S}) \succeq 0$, $\rank(\mathbf{C}(\mathbf{S})) = \rank(\mathbf{C})$ and the $(i,j)$th block of size $d$ in $\mathbf{C}(\mathbf{S})$ is $\mathbf{C}(\mathbf{S})_{ij} = \delta_{ij}\mathbf{D}(\mathbf{S})_{ii} - \mathbf{B}(\mathbf{S})_i\boldsymbol{\mathcal{L}}_{\Gamma}^\dagger \mathbf{B}(\mathbf{S})_j^T$ where $\mathbf{B}(\mathbf{S})_i$ is the $i$th row block of $\mathbf{B}(\mathbf{S})$ of dimension $d \times (m+n)$.
\end{rmk}
\begin{rmk}
\label{rmk:C_hat_L_structure}
\revadd{For a general $\widetilde{\mathbf{S}} \in \mathbb{O}(d)^m/_{\sim}$ and $\mathbf{S} \in \pi^{-1}(\mathbf{S})$, $\mathbf{L}(\mathbf{S})$ may not be symmetric and the $(i,j)$th block of size $d$ in $\mathbf{L}(\mathbf{S})$ is $\delta_{ij}\mathbf{B}(\mathbf{S})_i\boldsymbol{\mathcal{L}}_{\Gamma}^\dagger \mathbf{B}(\mathbf{S})^T\mathbf{I}^m_d - \mathbf{B}(\mathbf{S})_i\boldsymbol{\mathcal{L}}_{\Gamma}^\dagger \mathbf{B}(\mathbf{S})_j^T$. In particular, $\mathbf{L}(\mathbf{S})$ does not depend on $\mathbf{D}(\mathbf{S})$. Since the sum $\sum_{j=1}^{m}\mathbf{L}(\mathbf{S})_{ij} = 0$, for each $k \in [1,d]$, the vector $\mathbf{1}_m \otimes \mathbf{e}^d_k$ lies in the kernel of $\mathbf{L}(\mathbf{S})$ and thus the rank of $\mathbf{L}(\mathbf{S})$ is at most $(m-1)d$. Equivalently, if $\boldsymbol{\Omega} = [\boldsymbol{\Omega}_0]_1^m$ for some $\boldsymbol{\Omega}_0 \in \Skew(d)$ then $\mathbf{L}(\mathbf{S})\boldsymbol{\Omega} = 0$. Finally, if $\mathbf{Q} \in \mathbb{O}(d)$ then $\mathbf{L}(\mathbf{S}\mathbf{Q}) = (\mathbf{I}_m \otimes \mathbf{Q})^T\mathbf{L}(\mathbf{S})(\mathbf{I}_m \otimes \mathbf{Q})$ (follows from Eq.~(\ref{eq:C_of_S}, \ref{eq:C_hat}, \ref{eq:eq2})) and since $\mathbf{I}_m \otimes \mathbf{Q}$ is orthogonal, $L(\mathbf{S}\mathbf{Q})$ is unitarily equivalent to $\mathbf{L}(\mathbf{S})$.} Now, if $\widetilde{\mathbf{S}} \in \widetilde{\mathcal{C}}$ then for every $\mathbf{S} \in \pi^{-1}(\mathbf{S})$, $\mathbf{L}(\mathbf{S}) \in \Sym(md)$ and $\textstyle\sum_{j=1}^{m}\mathbf{L}(\mathbf{S})_{ij} = \textstyle\sum_{j=1}^{m}\mathbf{L}(\mathbf{S})_{ji} = 0$ for all $i \in [1,m]$ (see Eq.~(\ref{eq:eq2})).
\end{rmk}

We proceed to derive equations that will be used to determine the conditions under which the Hessian is non-singular and positive definite.
\begin{prop}
\label{prop:Omega_hat_compact}
Consider the same setup as in Proposition~\ref{prop:HessFSZ}. Then
\revadd{
\begin{equation}
    \widehat{\boldsymbol{\Omega}}_i = \left\{\begin{matrix} [\mathbf{L}(\mathbf{S})\boldsymbol{\Omega}]_i - [\mathbf{L}(\mathbf{S})\boldsymbol{\Omega}]_i^T, & \text{if } \widetilde{\mathbf{S}} \in \widetilde{\mathcal{C}}\\
    \frac{1}{2}\left\{[(\mathbf{L}(\mathbf{S}) + \mathbf{L}(\mathbf{S})^T)\boldsymbol{\Omega}]_i - [(\mathbf{L}(\mathbf{S}) + \mathbf{L}(\mathbf{S})^T)\boldsymbol{\Omega}]_i^T\right\}, & \text{otherwise}.\end{matrix}\right. \label{eq:eq4}
\end{equation}
}where $\boldsymbol{\Omega}=[\boldsymbol{\Omega}_i]_1^m$ and $\boldsymbol{\Omega}_i \in \Skew(d)$ is such that $\mathbf{Z}_i = \mathbf{S}_i\boldsymbol{\Omega}_i$ and $\textstyle\sum_1^m\boldsymbol{\Omega}_i = 0$.
\end{prop}

Combining the above with Proposition~\ref{prop:HessFSZ}, one can obtain a slightly compact representation of the horizontal lift of the Hessian, and the following result.
\begin{prop}
\label{prop:HessFSZZ}
Let \revadd{$\widetilde{\mathbf{S}} \in \mathbf{O}(d)^m/_{\sim}$}\revdel{$\widetilde{\mathbf{S}} \in \widetilde{\mathcal{C}}$} and $\widetilde{\mathbf{Z}} \in T_{\widetilde{\mathbf{S}}}\mathbb{O}(d)^{m}/_{\sim}$. Let $\mathbf{Z}$ be the horizontal lift of $\widetilde{\mathbf{Z}}$ at $\mathbf{S} \in \pi^{-1}(\widetilde{\mathbf{S}})$. Then
\begin{equation}
    \widetilde{g}(\Hess \widetilde{F}(\widetilde{\mathbf{S}})[\widetilde{\mathbf{Z}}],\widetilde{\mathbf{Z}}) = \revadd{\Tr(\boldsymbol{\Omega}^T(\mathbf{L}(\mathbf{S})+\mathbf{L}(\mathbf{S})^T)\boldsymbol{\Omega})} =  2\Tr(\boldsymbol{\Omega}^T\mathbf{L}(\mathbf{S})\boldsymbol{\Omega}) \label{eq:TrOmegaTLSOmega}
\end{equation}
where $\boldsymbol{\Omega}=[\boldsymbol{\Omega}_i]_1^m$ and $\boldsymbol{\Omega}_i \in \Skew(d)$ is such that $\mathbf{Z}_i = \mathbf{S}_i\boldsymbol{\Omega}_i$ and $\textstyle\sum_1^m\boldsymbol{\Omega}_i = 0$.
\end{prop}
Then the non-singularity and the positive definiteness of the Hessian amounts to the right side of Eq.~(\ref{eq:TrOmegaTLSOmega}) being non-zero and positive, respectively, for every non-zero $\boldsymbol{\Omega}$. Although $\mathbf{C}(\mathbf{S})$, and thus $\mathbf{L}(\mathbf{S})$, can be calculated from the patch framework $\Theta$ and the alignment $\mathbf{S}$, it is not obvious how to test the above practically. The main issue is that $\boldsymbol{\Omega}$ in Eq.~(\ref{eq:TrOmegaTLSOmega}) is not unconstrained, and in fact has a specific structure.

For the above reason, we are going to manipulate Eq.~(\ref{eq:TrOmegaTLSOmega}), utilizing the structure of $\boldsymbol{\Omega}$. The aim is to obtain an expression of the form $\boldsymbol{\omega}^T \mathbb{L}(\mathbf{S})\boldsymbol{\omega}$ where (i) the vector $\boldsymbol{\omega}$ is essentially unconstrained and (ii) $\boldsymbol{\Omega}$ and $\mathbf{L}(\mathbf{S})$ are related to $\boldsymbol{\omega}$ and $\mathbb{L}(\mathbf{S})$, respectively, through permutation matrices and vectorization operations (and thus the two pairs carry the same information). To achieve that, we first define certain matrices, then rewrite Eq.~(\ref{eq:TrOmegaTLSOmega}) in terms of those matrices and then describe their structure.


To this end, for $\boldsymbol{\Omega} = [\boldsymbol{\Omega}_i]_1^m$ where $\boldsymbol{\Omega}_i \in \Skew(d)$, let $\{\boldsymbol{\Omega}_{i}(r,s): 1 \leq r < s \leq d\}$ be the elements in the upper triangular region of $\boldsymbol{\Omega}_i$. For a fixed pair $(r,s)$ such that $1 \leq r < s \leq d$, define the column vector $\boldsymbol{\omega}_{r,s} \coloneqq [\boldsymbol{\Omega}_{i}(r,s)]_{i=1}^{m} \in \mathbb{R}^m$, a vertical stack of the $(r,s)$th element of each $\boldsymbol{\Omega}_i$. Then there exists a permutation matrix $\mathbf{P}$ such that
\begin{align}
    \mathbf{P}\boldsymbol{\Omega} &= \begin{bmatrix}
    \mathbf{0}_{m} & \boldsymbol{\omega}_{1,2} &  \ldots & \boldsymbol{\omega}_{1,d-1} & \boldsymbol{\omega}_{1,d}\\
    -\boldsymbol{\omega}_{1,2} & \mathbf{0}_{m}  & \ldots & \boldsymbol{\omega}_{2,d-1} & \boldsymbol{\omega}_{2,d}\\
        \vdots & \vdots & \vdots & \vdots & \vdots\\
    -\boldsymbol{\omega}_{1,d-1} & -\boldsymbol{\omega}_{2,d-1} & \ldots & \mathbf{0}_{m} & \boldsymbol{\omega}_{d-1,d}\\
    -\boldsymbol{\omega}_{1,d} & -\boldsymbol{\omega}_{2,d}  & \ldots & -\boldsymbol{\omega}_{d-1,d} & \mathbf{0}_{m}
    \end{bmatrix}. \label{eq:P0Omega}
\end{align}
In words, for $1 \leq r < s \leq d$, the $(r,s)$th block of $\mathbf{P}\boldsymbol{\Omega}$ is a vertical stack of the $(r,s)$th element of each $\boldsymbol{\Omega}_{i}$. For $r = s$, this is just a zero vector and for $r > s$, this is $-\boldsymbol{\omega}_{s,r}$.

Then, we collect the (strictly) upper triangular elements of $\mathbf{P}\boldsymbol{\Omega}$ in the column-major order in the vector $\boldsymbol{\omega}$. Note that $\mathbf{P}\boldsymbol{\Omega}$ can be fully described by $\boldsymbol{\omega}$. In particular, there exist a block matrix $\overline{\mathbf{P}}$ of size $d(d-1)/2 \times d^2$ siuch that $\vecz (\mathbf{P}\boldsymbol{\Omega}) = \overline{\mathbf{P}}^T \boldsymbol{\omega}$. The blocks of $\overline{\mathbf{P}}$ when indexed using tuples $(r,s)$ and $(p,q)$ where $1 \leq r < s \leq d$ and $p,q \in [1,d]$, are as follows: $\mathbf{0}_{m \times m}$ when $p=q$, $\delta_{pr}\delta_{qs}\mathbf{I}_m$ when $p < q$, and $-\delta_{ps}\delta_{qr}\mathbf{I}_m$ when $p > q$. Finally, we define $\mathcal{B}(\mathbf{S}) \coloneqq \mathbf{P}\mathbf{B}(\mathbf{S})$,
\begin{align}
    \boldsymbol{\mathcal{L}}(\mathbf{S}) &\coloneqq \mathbf{P}\mathbf{L}(\mathbf{S})\mathbf{P}^T \label{eq:mathcal_L}\\
    \mathbb{L}(\mathbf{S}) &\coloneqq \overline{\mathbf{P}}(\mathbf{I}_d \otimes \boldsymbol{\mathcal{L}}(\mathbf{S}))\overline{\mathbf{P}}^T \label{eq:mathbb_L}
\end{align}

\begin{prop}
\label{prop:Omega^TLSOmega2}
Consider the same setup as in Proposition~\ref{prop:HessFSZZ}. Then
\begin{equation}
    \Tr(\boldsymbol{\Omega}^T\mathbf{L}(\mathbf{S})\boldsymbol{\Omega}) = \boldsymbol{\omega}^T\mathbb{L}(\mathbf{S})\boldsymbol{\omega} \label{eq:omega^TmbbLomega}
\end{equation}
\end{prop}

The following remarks reveal the structure of $\mathcal{B}(\mathbf{S})$, $\boldsymbol{\mathcal{L}}(\mathbf{S})$ and $\mathbb{L}(\mathbf{S})$.
\begin{rmk}
\label{rmk:mathcalBS}
For $p \in [1,d]$, the $p$th row-block of $\mathcal{B}(\mathbf{S})$, $\mathcal{B}(\mathbf{S})_p$, is of size $m \times (n+m)$, and can be viewed as a vertical stack of the $p$th rows of $\mathbf{B}(\mathbf{S})_i$, $i \in [1,m]$. In particular, $\mathcal{B}(\mathbf{S})_p$ depends only on the $p$th coordinate of the local views (see Remark~\ref{rmk:L0DB}).
\end{rmk}
\begin{rmk}
\label{rmk:mathcalLS}
\revadd{For a general $\widetilde{\mathbf{S}} \in \mathbb{O}(d)^m/_{\sim}$ and $\mathbf{S} \in \pi^{-1}(\widetilde{\mathbf{S}})$, the matrix $\mathbf{\mathcal{L}}(\mathbf{S})$ may not be symmetric and for $p,q \in [1,d]$, $\boldsymbol{\mathcal{L}}(\mathbf{S})_{p,q} = \diag (\mathcal{B}(\mathbf{S})_p\boldsymbol{\mathcal{L}}_{\Gamma}^\dagger \mathcal{B}(\mathbf{S})_{q}^T\mathbf{1}_m) - \mathcal{B}(\mathbf{S})_p\boldsymbol{\mathcal{L}}_{\Gamma}^\dagger \mathcal{B}(\mathbf{S})_q^T$. Thus $\boldsymbol{\mathcal{L}}(\mathbf{S})_{p,q}$ depends on $\Gamma$ through $\boldsymbol{\mathcal{L}}_{\Gamma}^\dagger$ and the $p$th and $q$th coordinates of the rigidly transformed local views $\mathbf{B}(\mathbf{S})$. Since constant vectors are in the kernel of $\boldsymbol{\mathcal{L}}(\mathbf{S})_{p,q}$, therefore for each $p \in [1,d]$, the vector $\mathbf{e}^d_p \otimes \mathbf{1}_m$ lies in the kernel of $\boldsymbol{\mathcal{L}}(\mathbf{S})$. Thus the rank of $\boldsymbol{\mathcal{L}}(\mathbf{S})$ is atmost $(m-1)d$. Fianally, if $\mathbf{Q} \in \mathbb{O}(d)$ then, from Remark~\ref{rmk:C_hat_L_structure} and Eq.~(\ref{eq:mathcal_L}), it follows that $\boldsymbol{\mathcal{L}}(\mathbf{S}\mathbf{Q})$ is unitarily equivalent to $\boldsymbol{\mathcal{L}}(\mathbf{S})$. Now, if  $\widetilde{\mathbf{S}} \in \widetilde{\mathcal{C}}$ then $\mathbf{\mathcal{L}}(\mathbf{S})$ is symmetric and in particular $\diag (\mathcal{B}(\mathbf{S})_p\boldsymbol{\mathcal{L}}_{\Gamma}^\dagger \mathcal{B}(\mathbf{S})_{q}^T\mathbf{1}_m) = \diag (\mathcal{B}(\mathbf{S})_q\boldsymbol{\mathcal{L}}_{\Gamma}^\dagger \mathcal{B}(\mathbf{S})_{p}^T\mathbf{1}_m)$ which results in $\boldsymbol{\mathcal{L}}(\mathbf{S})_{p,q}^T = \boldsymbol{\mathcal{L}}(\mathbf{S})_{q,p}$ for all $p,q \in [1,d]$.
}\revdel{If $\widetilde{\mathbf{S}} \in \widetilde{\mathcal{C}}$ then for every $\mathbf{S} \in \pi^{-1}(\widetilde{\mathbf{S}})$, the following hold for $\mathbf{L}(\mathbf{S})$. First, since $\mathbf{L}(\mathbf{S})$ is symmetric therefore $\boldsymbol{\mathcal{L}}(\mathbf{S})_{q,p} = \boldsymbol{\mathcal{L}}(\mathbf{S})_{p,q}^T$. In particular, for $p,q \in [1,d]$, $\boldsymbol{\mathcal{L}}(\mathbf{S})_{p,q} = -\diag (\mathcal{B}(\mathbf{S})_p\boldsymbol{\mathcal{L}}_{\Gamma}^\dagger \mathcal{B}(\mathbf{S})_{q}^T\mathbf{1}_m) + \mathcal{B}(\mathbf{S})_p\boldsymbol{\mathcal{L}}_{\Gamma}^\dagger \mathcal{B}(\mathbf{S})_q^T$. Thus $\boldsymbol{\mathcal{L}}(\mathbf{S})_{p,q}$ depends on $\Gamma$ through $\boldsymbol{\mathcal{L}}_{\Gamma}^\dagger$ and the $p$th and $q$th coordinates of the rigidly transformed local views $\mathbf{B}(\mathbf{S})$. Moreover, $\boldsymbol{\mathcal{L}}(\mathbf{S})_{p,q}$ is a Laplacian-like matrix and the constant vectors are in its kernel. Moreover, for each $p \in [1,d]$, the vector $\mathbf{e}^d_p \otimes \mathbf{1}_m$ lies in the kernel of $\boldsymbol{\mathcal{L}}(\mathbf{S})$, thus the rank of $\boldsymbol{\mathcal{L}}(\mathbf{S})$ is atmost $(m-1)d$. Finally, if $\mathbf{Q} \in \mathbb{O}(d)$ then, from Remark~\ref{rmk:C_hat_L_structure} and Eq.~(\ref{eq:mathcal_L}), it follows that $\boldsymbol{\mathcal{L}}(\mathbf{S}\mathbf{Q})$ is unitarily equivalent to $\boldsymbol{\mathcal{L}}(\mathbf{S})$.}
\end{rmk}

\revadd{
The following characterization of the matrices $\mathcal{L}(\mathbf{S})_{p,p}$ will be useful for understanding the geometrical aspects of non-degenerate alignments in the case of $d=2$.}
\begin{rmk}
\label{rmk:mathcalLpp}
\revadd{For a general $\widetilde{\mathbf{S}} \in \mathbb{O}(d)^m/_{\sim}$ and $\mathbf{S} \in \pi^{-1}(\widetilde{\mathbf{S}})$, the symmetric Laplacian-like matrix $\mathcal{L}(\mathbf{S})_{p,p} = \diag (\mathcal{B}(\mathbf{S})_p\boldsymbol{\mathcal{L}}_{\Gamma}^\dagger \mathcal{B}(\mathbf{S})_{p}^T\mathbf{1}_m) - \mathcal{B}(\mathbf{S})_p\boldsymbol{\mathcal{L}}_{\Gamma}^\dagger \mathcal{B}(\mathbf{S})_p^T$ as defined in above remark can be constructed in the following way. Consider transforming all the views by $\mathbf{S}$ and projecting them to the $p$th coordinate, $p \in [1,d]$, i.e. replacing $\mathbf{x}_{k,i} \in \mathbb{R}^d$ to $\mathbf{S}_i^T\mathbf{x}_{k,i}(p) \in \mathbb{R}$. Now, to construct $\mathcal{L}(\mathbf{S})_{p,p}$, we (i) construct the positive semidefinite kernel $\mathbf{B}\boldsymbol{\mathcal{L}}_{\Gamma}^{\dagger}\mathbf{B}^T$ as in (Eq.~\ref{eq:GPOP}) using the one-dimensional patch framework $(\Gamma, (\mathbf{S}_i^T\mathbf{x}_{k,i}(p)))$, and then (ii) construct an unnormalized graph Laplacian matrix \citea{belkin2003laplacian} using the constructed kernel. Consequently, $\boldsymbol{\mathcal{L}}(\mathbf{S})_{p,p} \succeq 0$ and $\rank(\boldsymbol{\mathcal{L}}(\mathbf{S})_{p,p}) \leq m-1$.}
\end{rmk}
\begin{rmk}
\label{rmk:mathbb_L_structure}
\revadd{For $\widetilde{\mathbf{S}} \in \mathbb{O}(d)^m/_{\sim}$ and $\mathbf{S} \in \pi^{-1}(\widetilde{\mathbf{S}})$}, the following hold for $\mathbb{L}(\mathbf{S})$. First, $\mathbb{L}(\mathbf{S})$ is a block matrix of size $d(d-1)/2$ where each block is of size $m$. Indexing the rows and columns of $\mathbb{L}$ by tuples of the form $(r,s)$ where $1 \leq r < s \leq d$ we have,
\begin{align}
    \mathbb{L}(\mathbf{S})_{(r_1,s_1),(r_2,s_2)} 
    &= \left\{\begin{matrix}
    \boldsymbol{\mathcal{L}}(\mathbf{S})_{r,r}+\boldsymbol{\mathcal{L}}(\mathbf{S})_{s,s}, & r_1=r_2=r, s_1=s_2=s\\
    0, & \{r_1,s_1\} \cap \{r_2,s_2\} = \emptyset\\
    \boldsymbol{\mathcal{L}}(\mathbf{S})_{s_1,s_2}, & r_1 = r_2, s_1 \neq s_2\\
    \boldsymbol{\mathcal{L}}(\mathbf{S})_{r_1,r_2}, & s_1 = s_2, r_1 \neq r_2\\
    -\boldsymbol{\mathcal{L}}(\mathbf{S})_{r_1,s_2}, & s_1 = r_2\\
    -\boldsymbol{\mathcal{L}}(\mathbf{S})_{s_1,r_2}, & s_2 = r_1.\\
    \end{matrix}\right.
\end{align}
\revadd{For $d=2$, $3$ and $4$, $\mathbb{L}(\mathbf{S})$ is given by $\boldsymbol{\mathcal{L}}(\mathbf{S})_{1,1} + \boldsymbol{\mathcal{L}}(\mathbf{S})_{2,2}$,
\begin{equation}
    \begin{bmatrix}
    \boldsymbol{\mathcal{L}}(\mathbf{S})_{1,1} + \boldsymbol{\mathcal{L}}(\mathbf{S})_{2,2} & \boldsymbol{\mathcal{L}}(\mathbf{S})_{2,3} & -\boldsymbol{\mathcal{L}}(\mathbf{S})_{1,3}\\
    \boldsymbol{\mathcal{L}}(\mathbf{S})_{3,2} & \boldsymbol{\mathcal{L}}(\mathbf{S})_{1,1}+\boldsymbol{\mathcal{L}}(\mathbf{S})_{3,3} & \boldsymbol{\mathcal{L}}(\mathbf{S})_{1,2}\\
    -\boldsymbol{\mathcal{L}}(\mathbf{S})_{3,1} & \boldsymbol{\mathcal{L}}(\mathbf{S})_{2,1} & \boldsymbol{\mathcal{L}}(\mathbf{S})_{2,2}+\boldsymbol{\mathcal{L}}(\mathbf{S})_{3,3}
    \end{bmatrix}
\end{equation}
and (for brevity, here $\boldsymbol{\mathcal{L}} \equiv \boldsymbol{\mathcal{L}}(\mathbf{S})$)
\begin{equation}
    \begin{bmatrix}
    \boldsymbol{\mathcal{L}}_{1,1}+\boldsymbol{\mathcal{L}}_{2,2} & \boldsymbol{\mathcal{L}}_{2,3} & \boldsymbol{\mathcal{L}}_{2,4} & -\boldsymbol{\mathcal{L}}_{1,3} & -\boldsymbol{\mathcal{L}}_{1,4} & 0\\
    \boldsymbol{\mathcal{L}}_{3,2} & \boldsymbol{\mathcal{L}}_{1,1}+\boldsymbol{\mathcal{L}}_{3,3} & \boldsymbol{\mathcal{L}}_{3,4} & \boldsymbol{\mathcal{L}}_{1,2} & 0 & -\boldsymbol{\mathcal{L}}_{1,4}\\
    \boldsymbol{\mathcal{L}}_{4,2} & \boldsymbol{\mathcal{L}}_{4,3} & \boldsymbol{\mathcal{L}}_{1,1}+\boldsymbol{\mathcal{L}}_{4,4} & 0 & \boldsymbol{\mathcal{L}}_{1,2} & \boldsymbol{\mathcal{L}}_{1,3}\\
    -\boldsymbol{\mathcal{L}}_{3,1} & \boldsymbol{\mathcal{L}}_{2,1} & 0 & \boldsymbol{\mathcal{L}}_{2,2}+\boldsymbol{\mathcal{L}}_{3,3} & \boldsymbol{\mathcal{L}}_{3,4} & -\boldsymbol{\mathcal{L}}_{2,4}\\
    -\boldsymbol{\mathcal{L}}_{4,1} & 0 & \boldsymbol{\mathcal{L}}_{2,1} & \boldsymbol{\mathcal{L}}_{4,3} & \boldsymbol{\mathcal{L}}_{2,2}+\boldsymbol{\mathcal{L}}_{4,4} & \boldsymbol{\mathcal{L}}_{2,3}\\
    0 & -\boldsymbol{\mathcal{L}}_{4,1} & \boldsymbol{\mathcal{L}}_{3,1} & -\boldsymbol{\mathcal{L}}_{4,2} & \boldsymbol{\mathcal{L}}_{3,2} & \boldsymbol{\mathcal{L}}_{3,3}+\boldsymbol{\mathcal{L}}_{4,4}
    \end{bmatrix},
\end{equation}
respectively.
Here $\mathcal{L}(\mathbf{S})_{p,q}$ depends only on the $p$th and $q$th coordinates of the points in the local views, the structure of $\Gamma$ and the alignment $\mathbf{S}$.} Additionally, the set of vectors of the form $\boldsymbol{\omega} = [\boldsymbol{\omega}_{r,s}]_{1 \leq r < s \leq d}$ where each $\boldsymbol{\omega}_{r,s}$ is a constant vector, lie in the kernel of $\mathbb{L}(\mathbf{S})$. Therefore, the rank of $\mathbb{L}(\mathbf{S})$ is at most $(m-1)d(d-1)/2$.
\revadd{Now, if $\widetilde{\mathbf{S}} \in \widetilde{\mathcal{C}}$, then $\mathbb{L}(\mathbf{S})$ is also symmetric and in particular $\mathbb{L}(\mathbf{S})_{(r_1,s_1),(r_2,s_2)}^T = \mathbb{L}(\mathbf{S})_{(r_2,s_2),(r_1,s_1)}$}.
\end{rmk}

\subsection{Non-degenerate Alignment in the General Setting}
\label{subsec:non_deg_gen_setting}
As argued in Section~\ref{sec:setup}, since $F(\mathbf{S}) = F(\mathbf{S}\mathbf{Q})$ for all $\mathbf{Q} \in \mathbb{O}(d)$, every alignment $\mathbf{S}$ is degenerate in this sense. With a slight abuse of notation, we define a non-degenerate alignment as,
\begin{dfn}
\label{def:non_deg_alignment0}
An alignment $\mathbf{S} \in \mathbb{O}(d)^m$ is non-degenerate if $\pi(\mathbf{S})$ is a non-degenerate local minimum of $\widetilde{F}$.
\end{dfn}
With the above definition, to characterize the non-degenerate alignments, it suffices to characterize the non-degenerate local minima of $\widetilde{F}$. We accomplish the same in the following theorem. Note that we have not made any assumption about the affine non-degeneracy of the points and the noise in the local views.

\begin{thm}{\textbf{(Condition for $\widetilde{\mathbf{S}}$ to be a non-degenerate local minimum of $\widetilde{F}$)}}. Let $\widetilde{\mathbf{S}} \in \widetilde{\mathcal{C}}$ and $\mathbf{S} \in \pi^{-1}(\widetilde{\mathbf{S}})$. Then the following are equivalent.
\begin{enumerate}[leftmargin=*]
    \item $\widetilde{\mathbf{S}}$ is a non-degenerate local minimum of $\widetilde{F}$.
    \item $\widetilde{g}(\Hess \widetilde{F}(\widetilde{\mathbf{S}})[\widetilde{\mathbf{Z}}],\widetilde{\mathbf{Z}}) > 0$ for all $\widetilde{\mathbf{Z}} \in T_{\widetilde{\mathbf{S}}}\mathbb{O}(d)^m/_{\sim}$ such that $\widetilde{\mathbf{Z}} \neq 0$.
    \item $\Tr(\boldsymbol{\Omega}^T\mathbf{L}(\mathbf{S})\boldsymbol{\Omega}) > 0$ for all $\boldsymbol{\Omega} = [\boldsymbol{\Omega}_i]_1^m$ where $\boldsymbol{\Omega}_i \in \Skew(d)$, $\textstyle\sum_1^m \boldsymbol{\Omega}_i = 0$ and not all $\boldsymbol{\Omega}_i$ equal zero.
    \item $\Tr(\boldsymbol{\Omega}^T\mathbf{L}(\mathbf{S})\boldsymbol{\Omega}) > 0$ for all $\boldsymbol{\Omega} = [\boldsymbol{\Omega}_i]_1^m$ where $\boldsymbol{\Omega}_i \in \Skew(d)$ and not all $\boldsymbol{\Omega}_i$ are equal.
    \item $\boldsymbol{\omega}^T\mathbb{L}(\mathbf{S})\boldsymbol{\omega} > 0$ for all $\boldsymbol{\omega} = [\boldsymbol{\omega}_{r,s}]_{1 \leq r < s \leq d}$ where not all $\boldsymbol{\omega}_{r,s}$ are constant vectors.
    \item $\mathbb{L}(\mathbf{S})$ is positive semi-definite and of rank $(m-1)d(d-1)/2$.
    \item \revadd{$\lambda_{d(d-1)/2+1}(\mathbb{L}(\mathbf{S})) > 0$.}
\end{enumerate}
\label{thm:non_deg_loc_min}
\end{thm}
\begin{rmk}
\label{rmk:non_deg_loc_min}
Given the patch framework $\Theta$ and the alignment $\mathbf{S} \in \mathcal{C}$, one can compute the matrix $\mathbb{L}(\mathbf{S})$ in polynomial time in $m$, $n$ and $d$, and then check the non-degeneracy of the alignment $\mathbf{S}$ by testing the last condition in the above theorem (which again requires polynomial time in $m$ and $d$).
\end{rmk}

Although, $\widetilde{\mathbf{S}}$ is a non-degenerate local minimum of $\widetilde{F}$ if any of the equivalent conditions in the above theorem hold for every $\mathbf{S} \in \pi^{-1}(\widetilde{\mathbf{S}})$, the following result shows that if a conditions hold for one $\mathbf{S} \in \pi^{-1}(\widetilde{\mathbf{S}})$ then it holds for all other elements as well i.e. for all $\mathbf{S}\mathbf{Q}$ where $\mathbf{Q} \in \mathbb{O}(d)$ is arbitrary.
\begin{prop}
\label{prop:one_all1}
Let $\mathbf{S} \in \pi^{-1}(\widetilde{\mathbf{S}})$ and $\mathbf{Q} \in \mathbb{O}(d)$. Suppose \revadd{a} condition in Theorem~\ref{thm:non_deg_loc_min} holds for $\mathbf{S}$ then it holds for $\mathbf{S}\mathbf{Q}$ also. Consequently, an alignment $\mathbf{S}$ is non-degenerate if $\mathbf{S} \in \mathcal{C}$ (see Eq.~(\ref{eq:crit_pts2})) and it satisfies any of the (equivalent) conditions 3-7 in Theorem~\ref{thm:non_deg_loc_min}.
\end{prop}

\revadd{Following Remark~\ref{rmk:mathbb_L_structure} and Remark~\ref{rmk:mathcalLpp}, in the case of $d=2$ and an arbitrary number of views, the non-degeneracy of an alignment can be checked by investigating the the second smallest eigenvalues of certain laplacian-like matrices.
\begin{cor}
\label{cor:non_deg_d_2}
If $d=2$ then $\mathbf{S}$ is a non-degenerate alignment iff $\mathbf{\mathcal{L}}(\mathbf{S})_{1,1}+\mathbf{\mathcal{L}}(\mathbf{S})_{2,2}$ has rank of $m-1$ or equivalently, $\lambda_2(\mathbf{\mathcal{L}}(\mathbf{S})_{1,1} + \mathbf{\mathcal{L}}(\mathbf{S})_{2,2}) > 0$. In particular, $\mathbf{S}$ is a non-degenerate alignment if either $\mathbf{\mathcal{L}}(\mathbf{S})_{1,1}$ or $\mathbf{\mathcal{L}}(\mathbf{S})_{2,2}$ has a rank of $m-1$, equivalently $\max \{\lambda_2(\mathbf{\mathcal{L}}(\mathbf{S})_{1,1}), \lambda_2(\mathbf{\mathcal{L}}(\mathbf{S})_{2,2})\} > 0$.
\end{cor}
}
\begin{figure}[H]
    \centering
     \includegraphics[width=0.15\textwidth,keepaspectratio]{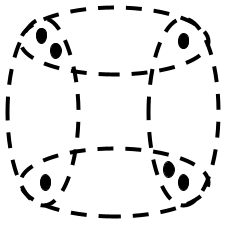}
    \caption{The dotted lines represent views and the filled points represent points on the overlaps. Here $d=2$ and all the pair of views are perfectly aligned (same is the case for the rest of the figures). It will be clear from Proposition~\ref{prop:noiseless_setting1} in Section~\ref{sec:noiseless_non_deg_results} that $\mathbf{L}(\mathbf{S}) \succeq 0$, and thus $\mathbb{L}(\mathbf{S}) \succeq 0$. Through simple calculations one can deduce that the rank of $\mathbb{L}(\mathbf{S})$ is $3$ (which equals $(m-1)d(d-1)/2$) while the rank of $\mathbf{L}(\mathbf{S})$ is $ 3 < 6 = (m-1)d$.}
    \label{fig:suff_cond_views_non_deg}
\end{figure}
A sufficient condition for a non-degenerate alignment in arbitrary dimensions is as follows. It is not a necessary condition though as demonstrated in Figure~\ref{fig:suff_cond_views_non_deg}.
\begin{cor}
\label{cor:suff_non_deg_loc_min}
If $\mathbf{L}(\mathbf{S}) \succeq 0$ and of rank $(m-1)d$, then $\mathbf{S}$ is a non-degenerate alignment. The same holds for $\boldsymbol{\mathcal{L}}(\mathbf{S})$ as it unitarily equivalent to $\mathbf{L}(\mathbf{S})$.
\end{cor}

\revadd{In the following we characterize the neighborhood of a non-degenerate alignment in which the Hessian remains non-singular and positive definite. This will be useful in characterizing the radius of convergence of Riemannian gradient descent in Section~\ref{sec:convergence}.}
\begin{prop}
\label{prop:HessVicinity}
\revadd{Let $\mathbf{S}$ be a non-degenerate alignment. For brevity, define $c_1 \coloneqq \max_{1}^{m}\sigma_{\max}(\mathbf{C}_{k,:})$, $c_2(\mathbf{S}) \coloneqq \max_{1}^{m}\sigma_{\max}([\mathbf{C}\mathbf{S}]_i)$, $c_3(\mathbf{S}) \coloneqq \sigma_{\max}(\mathbf{L}(\mathbf{S}))$, $\lambda_{-}(\mathbf{S}) \coloneqq \lambda_{d(d-1)/2+1}(\mathbb{L}(\mathbf{S})) > 0$ (follows from Theorem~\ref{thm:non_deg_loc_min}) and $\lambda_{+}(\mathbf{S}) \coloneqq \lambda_{md(d-1)/2}(\mathbb{L}(\mathbf{S})) > 0$. Note that the above functions are invariant under the action of $\mathbb{O}(d)$ i.e. they have the same value at $\mathbf{S}\mathbf{Q}$ for all $\mathbf{Q} \in \mathbb{O}(d)$. Let $\zeta \in (0,1)$ be fixed and define}
\begin{equation}
\label{eq:delta}
\revadd{\delta(\mathbf{S}) \coloneqq \lambda_{-}(\mathbf{S})/ 2(c_1 + c_2(\mathbf{S}) + 2 c_3(\mathbf{S})).}
\end{equation}
\revadd{If $\mathbf{O} \in \mathbb{O}(d)^m$ satisfies $\min_{\mathbf{Q}\in\mathbb{O}(d)}\left\|\mathbf{O}-\mathbf{S}\mathbf{Q}\right\|_F \leq \zeta\delta(\mathbf{S})$, then for all $\widetilde{\mathbf{Z}} \in T_{\widetilde{\mathbf{O}}}\mathbb{O}(d)^m/_{\sim}$,}
\begin{equation}
    \revadd{(1-\zeta)\lambda_{-}(\mathbf{S})\widetilde{g}(\widetilde{\mathbf{Z}}, \widetilde{\mathbf{Z}}) \leq \widetilde{g}(\Hess \widetilde{F}(\widetilde{\mathbf{O}})[\widetilde{\mathbf{Z}}],\widetilde{\mathbf{Z}}) \leq (\lambda_{+}(\mathbf{S})+\zeta \lambda_{-}(\mathbf{S})) \widetilde{g}(\widetilde{\mathbf{Z}}, \widetilde{\mathbf{Z}}).}
\end{equation}
\end{prop}
We end this subsection by deriving a necessary and sufficient condition for an alignment of two views to be non-degenerate. First we need the following definitions (note that the objects in these definitions are related but not identical to $\mathbf{B}_i$ (see Eq.~(\ref{eq:B}), Remark~\ref{rmk:L0DB}) and $\mathbf{B}(\mathbf{S})_i$ (see Eq.~(\ref{eq:BofS}), Remark~\ref{rmk:C_S_structure})),
\begin{dfn}
\label{def:Bij}
Let $i,j \in [1,m]$ be indices of two views. Define $\mathbf{B}_{i,j}$ to be a matrix whose columns are $\mathbf{x}_{k,i}$ (in the increasing order of $k$) where $(k,i),(k,j) \in E(\Gamma)$. Generally, $\mathbf{B}_{i,j} \neq \mathbf{B}_{j,i}$. Also, define $\overline{\mathbf{B}}_{i,j} = \mathbf{B}_{i,j}\left(\mathbf{I}_{n'} - (1/n')\mathbf{1}_{n'}\mathbf{1}_{n'}^T\right)$ where $n' = |\{k: (k,i),(k,j) \in E(\Gamma)\}|$ is the number of points on the overlap of the $i$th view and the $j$th view, equivalently the number of columns in $\mathbf{B}_{i,j}$.
\end{dfn}
\begin{dfn}
\label{def:BSicapj}
In a similar manner as above, let $i,j \in [1,m]$ be indices of two views and let $\mathbf{S}$ be an alignment. Define $\mathbf{B}(\mathbf{S})_{i,j}$ to be a matrix whose columns are $\mathbf{S}_i^T\mathbf{x}_{k,i}+\mathbf{t}_i$ (in increasing order of $k$) where $(k,i),(k,j) \in E(\Gamma)$ and where $\mathbf{t}_i$ is obtained using Eq.~(\ref{eq:opt_Z}). Also, define $\overline{\mathbf{B}(\mathbf{S})}_{i,j} = \mathbf{B}(\mathbf{S})_{i,j}\left(\mathbf{I}_{n'} - (1/n')\mathbf{1}_{n'}\mathbf{1}_{n'}^T\right)$.
\end{dfn}
\begin{rmk}
\label{rmk:BS_ijB_ij}
Let $i,j \in [1,m]$ be indices of two views and $\mathbf{S}$ be an alignment. Let  $n' = |\{k:(k,i),(k,j) \in E(\Gamma)\}|$ be the number of points on the overlap of the two views. Then $\mathbf{B}(\mathbf{S})_{i,j} = \mathbf{S}_i^T \mathbf{B}_{i,j} + \mathbf{t}_i\mathbf{1}_{n'}^T$ where $\mathbf{t}_i$ is obtained using Eq.~(\ref{eq:opt_Z}). Thus, we have $\rank (\overline{\mathbf{B}}_{i,j}) = \rank (\overline{\mathbf{B}(\mathbf{S})}_{i,j})$.
Moreover, $\mathbf{B}(\mathbf{S})_{i,j}\left(\mathbf{I}_{n'} - (1/n'){\mathbf{1}_{n'}\mathbf{1}_{n'}^T}\right)\mathbf{B}(\mathbf{S})_{j,i}^T$ equals $\mathbf{S}_i^T\mathbf{B}_{i,j}\left(\mathbf{I}_{n'} - (1/n'){\mathbf{1}_{n'}\mathbf{1}_{n'}^T}\right)\mathbf{B}_{j,i}^T\mathbf{S}_j$ which in turn equals $\mathbf{S}_1^T\overline{\mathbf{B}}_{i,j}\overline{\mathbf{B}}_{j,i}^T\mathbf{S}_2$.
Consequently, the $\rank (\overline{\mathbf{B}(\mathbf{S})}_{i,j}\overline{\mathbf{B}(\mathbf{S})}_{j,i}^T) = \rank (\overline{\mathbf{B}}_{i,j}\overline{\mathbf{B}}_{j,i}^T)$.
\end{rmk}
\begin{figure}[H]
    \centering
    \begin{tabular}{ccc}
    \begin{subfigure}[b]{0.175\textwidth}
         \centering
         \includegraphics[width=0.9\textwidth,keepaspectratio]{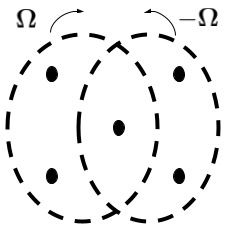}
         \caption{}
         \label{fig:nec_suff_cond_loc_rigid_two_views_1}
     \end{subfigure}
     &
     \begin{subfigure}[b]{0.175\textwidth}
         \centering
         \includegraphics[width=0.9\textwidth,keepaspectratio]{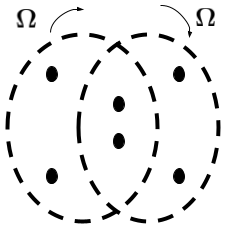}
         \caption{}
         \label{fig:nec_suff_cond_loc_rigid_two_views}
     \end{subfigure}
     &
     \begin{subfigure}[b]{0.175\textwidth}
         \centering
         \includegraphics[width=0.9\textwidth,keepaspectratio]{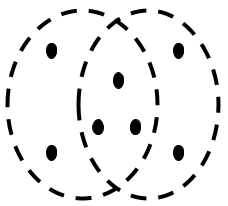}
         \caption{}
         \label{fig:nec_suff_cond_glob_rigid_two_views}
     \end{subfigure}
     \end{tabular}
    \caption{(a) $\rank (\overline{\mathbf{B}}_{1,2}\overline{\mathbf{B}}_{2,1}^T) = 0$ and the two views can be rotated by a different amount while still being perfectly aligned. (b) $\rank (\overline{\mathbf{B}}_{1,2}\overline{\mathbf{B}}_{2,1}^T) = 1$ and in order for the views to be perfectly aligned, every infinitesimal rotation of the two views must be identical. However the perfect alignment of the views is not unique because the second view can be flipped (a non-infinitesimal rotation) to obtain another perfect alignment of the views. (c) $\rank (\overline{\mathbf{B}}_{1,2}\overline{\mathbf{B}}_{2,1}^T) = 2$ and the perfect alignment is unique.}
    \label{fig:geom_intuit}
\end{figure}
\begin{thm}
\label{thm:non_deg_two_views_gen_setting}
Consider $m=2$ and let $\mathbf{S} \in \mathbb{O}(d)^2$. Then $\mathbf{S}$ is a non-degenerate alignment iff all of the following hold: (see Figures~\ref{fig:nec_suff_cond_loc_rigid_two_views_1} and \ref{fig:nec_suff_cond_loc_rigid_two_views} for intuition when $d=2$)
\begin{enumerate}[leftmargin=*]
    \itemsep0em 
    \item $\overline{\mathbf{B}(\mathbf{S})}_{1,2}\overline{\mathbf{B}(\mathbf{S})}_{2,1}^T$ is symmetric.
    \item $\mathrm{Tr}(\boldsymbol{\Omega}^T \overline{\mathbf{B}(\mathbf{S})}_{1,2}\overline{\mathbf{B}(\mathbf{S})}_{2,1}^T\boldsymbol{\Omega}) \geq 0$ for all $\boldsymbol{\Omega} \in \Skew(d)$.
    \item $\rank\left(\overline{\mathbf{B}(\mathbf{S})}_{1,2}\overline{\mathbf{B}(\mathbf{S})}_{2,1}^T\right) \geq d-1$ (equivalently $\rank (\overline{\mathbf{B}}_{1,2}\overline{\mathbf{B}}_{2,1}^T) \geq d-1$).
\end{enumerate}
\end{thm}

\subsection{Unique Optimal Alignment in the General Setting}
\label{subsec:uniq_gen_setting}
Since $F(\mathbf{S}) = F(\mathbf{S}\mathbf{Q})$ for all $\mathbf{Q} \in \mathbb{O}(d)$, if $\mathbf{S}$ is an optimal alignment (a global minimum) then so is $\mathbf{S}\mathbf{Q}$. In this sense, no optimal alignment is unique. With a slight abuse of convention we define a unique optimal alignment below.
\begin{dfn}
\label{def:uniq_alignment}
An alignment $\mathbf{S} \in \mathbb{O}(d)^m$ is a unique optimal alignment if $\pi(\mathbf{S})$ is the unique global minimum of $\widetilde{F}$ i.e. if $\mathbf{O} \in \mathbb{O}(d)^m$ is also an optimal alignment then $\pi(\mathbf{O}) = \pi(\mathbf{S})$, equivalently, $\mathbf{O} = \mathbf{S}\mathbf{Q}$ for some $\mathbf{Q} \in \mathbb{O}(d)$.
\end{dfn}

\begin{thm}
\label{thm:uniq_two_views_gen_setting}
Let $m=2$ and $\mathbf{S}$ be an optimal alignment. Then $\mathbf{S}$ is unique iff $\rank (\overline{\mathbf{B}}_{1,2}\overline{\mathbf{B}}_{2,1}^T) = d$ (see Figures~\ref{fig:nec_suff_cond_loc_rigid_two_views} and \ref{fig:nec_suff_cond_glob_rigid_two_views} for $d=2$, and \citea{schonemann1966generalized} for the proof).
\end{thm}

\section{Non-degeneracy and Uniqueness in the Noiseless Regime}
\label{sec:noiseless_non_deg_results}
\begin{figure}[H]
    \centering
    \resizebox{0.65\textwidth}{!}{%
    \begin{tikzpicture}[
        NodeA/.style={rectangle, draw=black!60, fill=green!5, very thick, minimum size=7mm, text width=3cm},
        NodeB/.style={rectangle, draw=black!60, fill=red!5, very thick, minimum size=5mm, text width=3cm},
        NodeC/.style={rectangle, draw=black!60, fill=blue!5, very thick, minimum size=7mm, text width=3cm},
        NodeD/.style={rectangle, draw=black!60, fill=blue!5, very thick, minimum size=7mm, text width=5cm},
        every text node part/.style={align=center}
        ]
        \node[NodeB]        (irigid)  {infinitesimally rigid $\Theta(\mathbf{S})$};
        \node[NodeC]      (rankR) [above=0.85cm of irigid] {$\rank(\boldsymbol{\mathcal{R}}(\mathbf{S})) \geq nd - d(d+1)/2$};
        \node[NodeB]      (lrigid)       [below=of irigid] {locally rigid $\Theta(\mathbf{S})$};
        \node[NodeB]      (grigid)   [below=of lrigid]        {globally rigid $\Theta(\mathbf{S})$};
        \node[NodeB]        (arigid)       [below=of grigid] {affinely rigid $\Theta(\mathbf{S})$};

        \node[NodeA]        (nondegS)       [right=3cm of irigid] {non-degenerate $\mathbf{S}$};
        \node[NodeA]        (strictS)       [right=3cm of lrigid] {$\pi(\mathbf{S})$ is a strict minimum of $\widetilde{F}$};
        \node[NodeC]        (rankLbb)       [above=of nondegS] {$\rank(\boldsymbol{\mathbb{L}}(\mathbf{S})) = (m-1)d(d-1)/2$};
        \node[NodeA]        (uniqueS)       [right=3cm of grigid] {unique $\mathbf{S}$};
        \node[NodeC]        (rankC)       [right=3cm of arigid] {$\rank(\mathbf{C}) = (m-1)d$};
        \coordinate[below right=2cm and 0.5cm of rankLbb]  (rankLbb0) ;
        \coordinate[above right=2cm and 0.5cm of rankC]  (rankC0) ;

        \path (rankLbb) -- node (rankLbbiffnondegS) {Proposition~\ref{prop:noiseless_setting1}} (nondegS);
        \path (nondegS) -- node (nondegSimpstrictS) {Trivial} (strictS);
        \path (lrigid.20) -- node (lrigidimpirigid) {Generic $\Theta(\mathbf{S})$} (irigid.330);
        
        \draw[implies-implies, double, line width=0.4mm] (irigid.east) -- node [midway,above] {Theorem~\ref{thm:inf_rigid}} (nondegS.west);
        \draw[-implies,double, line width=0.4mm] (irigid.230) --node[midway,left] {\citea{toth2017handbook}} (lrigid.145);
        \draw[-implies,double, line width=0.4mm] (lrigid.20) -- (lrigidimpirigid) -- (irigid.330);
        \draw[implies-implies, double, line width=0.4mm] (lrigid.east) -- node[midway,above] {Proposition~\ref{prop:non_deg_views}} (strictS.west);
        \draw[-implies,double, line width=0.4mm] (grigid.north) -- node [midway,left] {\citea{gortler2010affine}} (lrigid.south);
        \draw[implies-implies,double, line width=0.4mm] (irigid.north) -- node [midway,left] {\citea{toth2017handbook}} (rankR.south);
        \draw[implies-implies,double, line width=0.4mm] (grigid.east)--node[midway,above]{Theorem~\ref{thm:glob_rigid}} (uniqueS.west);
        \draw[-implies,double, line width=0.4mm] (arigid.north) -- node [midway,left] {\citea{gortler2010affine}} (grigid.south);

        \draw[-,double, line width=0.4mm] (rankC.east) -| (rankC0);
        \draw[-,double, line width=0.4mm] (rankLbb0) --node[above,midway,sloped] {Corollary~\ref{cor:noiseless_setting1}} (rankC0);
        \draw[-implies,double, line width=0.4mm] (rankLbb0) |- (rankLbb.east) ;

        \draw[implies-implies,double, line width=0.4mm] (rankC.west) -- node [midway, above] {\citea{chaudhury2015global,zha2009spectral}} (arigid.east);

        \draw[implies-,double, line width=0.4mm] (rankLbb.south)--(rankLbbiffnondegS);
        \draw[-implies,double, line width=0.4mm] (rankLbbiffnondegS)--(nondegS.north);
        \draw[-implies,double, line width=0.4mm] (nondegS.south) --(nondegSimpstrictS)-- (strictS.north);
        
    \end{tikzpicture}
    }
    \caption{\revadd{The implications between the type of a perfect alignment $\mathbf{S}$ and the rigidity of the resulting realization $\Theta(\mathbf{S})$.}}
    \label{fig:rigidity_flow}
\end{figure}
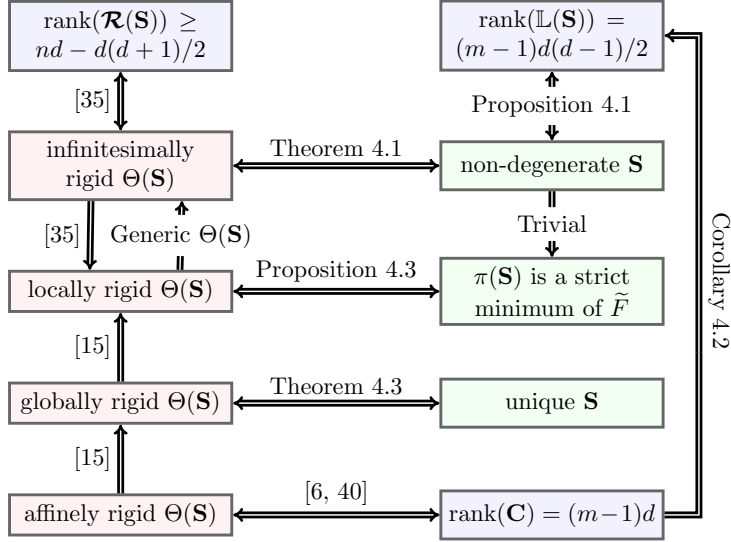

We start by deriving some important consequences of the noiseless setting. Under a mild assumption on the structure of the local views, we show that the non-degeneracy and uniqueness of a perfect alignment is equivalent to certain notions of rigidity of the resulting realization (Figure~\ref{fig:rigidity_flow}). We then provide necessary and sufficient conditions on the overlapping structure of the views for the non-degeneracy and uniqueness of a perfect alignment. \revadd{Consequently, we obtain conditions on a perfect alignment for the resulting realization to be infinitesimally/locally/globally rigid. These should be contrasted with the ones in \citea{zha2009spectral}, for the affine rigidity of a realization.}
\subsection{Consequences of Noiseless Setting}
\label{subsec:noiseless_conseq}
As discussed in \citea{chaudhury2015global}, in the noiseless case, the patch-stress matrix $\mathbf{C}$ is constructed from $\Gamma$ and clean measurements. In particular, there exists a perfect alignment $\mathbf{S}$ such that $F(\mathbf{S}) = 0$.
\begin{prop}
\label{prop:noiseless_setting1}
Let $\mathbf{S}$ be a perfect alignment. Then $\widehat{\mathbf{C}}(\mathbf{S}) = 0$, $\mathbf{L}(\mathbf{S}) = \mathbf{C}(\mathbf{S})$ (see Eq.~(\ref{eq:L_of_S})) and $\mathbb{L}(\mathbf{S}) = \overline{\mathbf{P}}(\mathbf{I}_d \otimes  (\mathbf{P}\mathbf{C}(\mathbf{S})\mathbf{P}^T))\overline{\mathbf{P}}^T$ (see Eq.~(\ref{eq:mathbb_L})). Consequently, it is easy to deduce from Remark~\ref{rmk:C_S_structure} that $\mathbf{L}(\mathbf{S}) \succeq 0$ and $\mathbb{L}(\mathbf{S}) \succeq 0$. \revadd{It follows from Theorem~\ref{thm:non_deg_loc_min} that $\mathbf{S}$ is non-degenerate if and only if $\rank(\mathbb{L}(\mathbf{S})) = (m-1)d(d-1)/2$.}
\end{prop}

\begin{rmk}
\label{rmk:C1Sp}
\revadd{
Following Remark~\ref{rmk:mathcalLpp} and due to the above proposition, it is easy to deduce that for a perfect alignment $\mathbf{S}$, $\mathbf{\mathcal{L}}(\mathbf{S})_{pp}$ is exactly the patch-stress matrix (see Eq.~(\ref{eq:GPOP})) of the one-dimensional patch framework  $(\Gamma, (\mathbf{S}_i^T\mathbf{x}_{k,i}(p)))$.}
\end{rmk}

\revadd{As a direct corollary of Proposition~\ref{prop:HessVicinity}, we obtain a bound on the neighborhood of a non-degenerate perfect alignment where the Hessian is positive~definite.}
\begin{cor}
\label{cor:HessVicinity}
\revadd{Let $\mathbf{S}$ be a non-degenerate perfect alignment.
As in Proposition~\ref{prop:HessVicinity}, define $c_1 = \max_{1}^{m}\sigma_{\max}(\mathbf{C}_{k,:})$, $c_3 \coloneqq \sigma_{\max}(\mathbf{C})$, $\lambda_{0_-}(\mathbf{S}) \coloneqq \lambda_{d(d-1)/2+1}(\mathbb{L}(\mathbf{S}))$ and $\lambda_{0_+}(\mathbf{S}) \coloneqq \lambda_{md(d-1)/2}(\mathbb{L}(\mathbf{S}))$ (from Proposition~\ref{prop:noiseless_setting1} and Eq.~(\ref{eq:omega^TmbbLomega}) it is easy to deduce that $0 < \lambda_{0_-}(\mathbf{S}) \leq \lambda_{0_+}(\mathbf{S}) \leq 2\lambda_{md}(\mathbf{C})$). Let $\zeta \in (0,1)$ be fixed and define}
\begin{equation}
\label{eq:delta}
\revadd{\delta_0(\mathbf{S}) \coloneqq \lambda_{0_-}(\mathbf{S})/ 2(c_1 + 2 c_{3}).}
\end{equation}
\revadd{If $\mathbf{O} \in \mathbb{O}(d)^m$ satisfies $\min_{\mathbf{Q}\in\mathbb{O}(d)}\left\|\mathbf{O}-\mathbf{S}\mathbf{Q}\right\|_F \leq \zeta\delta_0(\mathbf{S})$, then for all $\widetilde{\mathbf{Z}} \in T_{\widetilde{\mathbf{O}}}\mathbb{O}(d)^m/_{\sim}$,}
\begin{equation}
    \revadd{(1-\zeta)\lambda_{0_-}(\mathbf{S})\widetilde{g}(\widetilde{\mathbf{Z}}, \widetilde{\mathbf{Z}}) \leq \widetilde{g}(\Hess \widetilde{F}(\widetilde{\mathbf{O}})[\widetilde{\mathbf{Z}}],\widetilde{\mathbf{Z}}) \leq (\lambda_{0_+}(\mathbf{S}) + \zeta\lambda_{0_-}(\mathbf{S})) \widetilde{g}(\widetilde{\mathbf{Z}}, \widetilde{\mathbf{Z}}).}
\end{equation}
\end{cor}

Finally, \revadd{using Proposition~\ref{prop:one_all1} and Proposition~\ref{prop:noiseless_setting1},} we provide a simplified characterization of a non-degenerate perfect alignment that will be useful in proving the subsequent results. First, similar to \citea{zha2009spectral}, we define a certificate of $\mathbf{L}(\mathbf{S})$.
\begin{dfn}
\label{def:LScertificate}
An $\boldsymbol{\Omega} \in \Skew(d)^m$ is said to be a certificate of $\mathbf{L}(\mathbf{S})$ if $\mathbf{L}(\mathbf{S})\boldsymbol{\Omega} = 0$. It is a trivial certificate if $\boldsymbol{\Omega}_i = \boldsymbol{\Omega}_0$ for all $i \in [1,m]$ and for some $\boldsymbol{\Omega}_0 \in \Skew(d)$.
\end{dfn}

\begin{prop}
\label{prop:non_deg_triv_cert}
If $\mathbf{S}$ is a perfect alignment then $\mathbf{S}$ is non-degenerate iff every certificate of $\mathbf{L}(\mathbf{S})$ is trivial.
\end{prop}

\subsection{Rigidity of a Realization}
\label{subsec:loc_glob_rigid}
In the following, we reveal the relation between non-degenerate and unique perfect alignment with the various notions of the rigidity of the resulting realization. \revadd{These are summarized in Figure~\ref{fig:rigidity_flow}.}
\revadd{Throughout the rest of this work, we assume the following.}
\begin{assump}
\label{assump:non_deg_views}
\revadd{Each view is affinely non-degenerate i.e. has at least $d+1$ points whose affine span has a rank of $d$.}
\end{assump}

\revadd{Consequently, the perfect alignment of the local views can be uniquely determined by their realization. This can be easily inferred from the following result.}
\begin{prop}
\label{prop:non_deg_views}
\revadd{Let $\mathbf{B}_{i,i}$, $i \in [1,m]$, be as in Definition~\ref{def:Bij}. Define $\varrho = (\sum_1^m 1/\sigma_{\min}(\mathbf{B}_{i,i}\mathbf{B}_{i,i}^T)^{2})^{1/2}$. Then for perfect alignments $\mathbf{S}$ and $\mathbf{O}$, and the corresponding realizations $\Theta(\mathbf{S})$ and $\Theta(\mathbf{O})$, $\left\|\mathbf{S} - \mathbf{O}\right\|_F \leq \varrho \left\|\Theta(\mathbf{S}) - \Theta(\mathbf{O})\right\|_F$. In particular, $\Theta(\mathbf{S}) = \Theta(\mathbf{O})$ if and only if $\mathbf{S} = \mathbf{O}$.}
\end{prop}

Now we define various notions of the rigidity of a realization $\Theta(\mathbf{S})$. Although phrased differently, the definitions are the same as those in \citea{toth2017handbook, gortler2010affine, chaudhury2015global}.
\begin{dfn}
\label{def:inf_rigid}
\revadd{Let $\mathbf{S}$ be a perfect alignment. Then $\Theta(\mathbf{S}) = (\mathbf{x}_k(\mathbf{S}))_1^n$ is infinitesimally rigid if there does not exist a perturbation $(\mathbf{p}_k)_1^n \subseteq \mathbb{R}^d$ satisfying:
\begin{enumerate}[leftmargin=*]
    \item $(\mathbf{p}_k)_1^n$ is not a trivial perturbation (it is a trivial perturbation if there exist $\boldsymbol{\Omega} \in \Skew(d)$ and $\mathbf{t} \in \mathbb{R}^d$ such that $\mathbf{p}_k = \boldsymbol{\Omega}\mathbf{x}_k + \mathbf{t}$),
    \item and for all $(k_1,i), (k_2,i) \in E(\Gamma)$, $(\mathbf{x}_{k_1}(\mathbf{S})-\mathbf{x}_{k_2}(\mathbf{S}))^T(\mathbf{p}_{k_1}-\mathbf{p}_{k_2}) = 0$.
\end{enumerate}
}
\begin{rmk}
\label{rmk:inf_rigid}
\revadd{The above two conditions can be described in terms of the rank of the so-called rigidity matrix $\boldsymbol{\mathcal{R}}(\mathbf{S})$ \citea{toth2017handbook}. The rigidity matrix has a row for each triplet $(k_1, k_2, i)$ satisfying $(k_1,i), (k_2,i) \in E(\Gamma)$, and the $k_1$th and $k_2$th the blocks of the row are $(\mathbf{x}_{k_1}(\mathbf{S})-\mathbf{x}_{k_2}(\mathbf{S}))^T$ and $(\mathbf{x}_{k_2}(\mathbf{S})-\mathbf{x}_{k_1}(\mathbf{S}))^T$, respectively. Overall, the sparse matrix $\boldsymbol{\mathcal{R}}(\mathbf{S})$ has $\sum_{1}^{m}{n_i \choose 2}$ rows and $nd$ columns, and the realization $\Theta(\mathbf{S}) = (\mathbf{x}_k(\mathbf{S}))_1^n$ is infinitesimally rigid if and only if $\rank(\boldsymbol{\mathcal{R}}(\mathbf{S})) \geq nd - d(d+1)/2$.}
\end{rmk}
\end{dfn}
For the following definitions, we use the facts due to Definition~\ref{def:realization}: (i) $\mathbf{Q}^T\Theta(\mathbf{S}) = \Theta(\mathbf{S}\mathbf{Q})$ for any $\mathbf{Q} \in \mathbb{R}^{d \times d}$ and (ii) $\mathbf{0}_{d} = \textstyle\argmin_{\mathbf{t} \in \mathbb{R}^d} \left\|\Theta(\mathbf{O}) - \mathbf{Q}^T\Theta(\mathbf{S}) - \mathbf{t}\mathbf{1}_n^T\right\|_F$.
\begin{dfn}
\label{def:loc_rigid}
Let $\mathbf{S}$ be a perfect alignment. Then $\Theta(\mathbf{S})$ is locally rigid if there exists $\epsilon > 0$ such that for any other perfect alignment $\mathbf{O} \in \mathbb{O}(d)^m$ with $\left\|\Theta(\mathbf{O})-\Theta(\mathbf{S})\right\|_F < \epsilon$, we have $\Theta(\mathbf{O})$ to be a rigid transformation of $\Theta(\mathbf{S})$ or equivalently $\Theta(\mathbf{O})  = \Theta(\mathbf{S}\mathbf{Q})$ for some $\mathbf{Q} \in \mathbb{O}(d)$.
\end{dfn}
\begin{dfn}
\label{def:glob_rigid}
Let $\mathbf{S}$ be a perfect alignment. Then $\Theta(\mathbf{S})$ is globally rigid if for any other perfect alignment $\mathbf{O} \in \mathbb{O}(d)^m$ we have $\Theta(\mathbf{O}) = \Theta(\mathbf{S}\mathbf{Q})$ for some $\mathbf{Q} \in \mathbb{O}(d)$.
\end{dfn}
\begin{dfn}
\label{def:affine_rigid}
\revadd{Let $\mathbf{S}$ be a perfect alignment. Then $\Theta(\mathbf{S})$ is affinely rigid if for any realization $\mathbf{Y} \in \mathbb{R}^{d \times n}$ satisfying: for each $i \in [1,m]$ there exist an affine transform $\mathbf{A}_i$ such that $\mathbf{Y}_k = \mathbf{A}_i(\mathbf{x}_{k,i})$, we have $\mathbf{Y} = \mathbf{A}\Theta(\mathbf{S})$ for some global affine transform~$\mathbf{A}$.}
\end{dfn}
\revadd{From the above definitions, it is easy to see that an affinely rigid realization is globally rigid which in turn is locally rigid.} Examples of realizations that are not locally rigid, locally rigid but not globally rigid and \revadd{globally rigid but not affinely rigid} are provided in Figure~\ref{fig:nec_cond_loc_rigid_of_views}, Figure~\ref{fig:suff_cond_views_non_deg}, Figure~\ref{fig:G_star_1} and \revadd{\cite[Figure~3]{gortler2010affine}} respectively.
\revadd{Follows our first result connecting type of a perfect alignment with the rigidity of the resulting realization.}
\begin{thm}
\label{thm:inf_rigid}
\revadd{Let $\mathbf{S}$ be a perfect alignment. The realization $\Theta(\mathbf{S})$ is infinitesimally rigid if and only if the alignment $\mathbf{S}$ is non-degenerate.}
\end{thm}
\revadd{Due to the proof of the above theorem, Proposition~\ref{prop:noiseless_setting1} and Remark~\ref{rmk:inf_rigid}, one can derive non-trivial perturbations of a non-infinitesimally rigid realization $\Theta(\mathbf{S})$ by using the non-trivial vectors in the null space of $\boldsymbol{\mathcal{R}}(\mathbf{S})$ or $\boldsymbol{\mathbb{L}}(\mathbf{S})$. Furthermore, in the noiseless setting, one can test if a perfect alignment $\mathbf{S}$ is non-degenerate - either by checking if $\rank(\boldsymbol{\mathcal{R}}(\mathbf{S})) \geq nd-d(d+1)/2$ or if $\rank(\boldsymbol{\mathbb{L}}(\mathbf{S})) \geq md(d-1)/2$.}

\revadd{It is well known that an infinitesimally rigid realization is also locally rigid \citea{toth2017handbook}, and the converse holds for generic realizations ($\Theta(\mathbf{S}) = (\mathbf{x}_k(\mathbf{S}))_1^n$ is generic if the coordinates do not satisfy any non-zero algebraic equation with rational coefficients). Here, we provide a result which elucidates a more clear picture.}
\begin{thm}
\label{thm:loc_rigid}
Let $\mathbf{S}$ be a perfect alignment.
Then the realization $\Theta(\mathbf{S})$ is locally rigid iff $\pi(\mathbf{S})$ is a strict global minimum of $\widetilde{F}$. Consequently, if $\mathbf{S}$ is a non-degenerate perfect alignment then $\Theta(\mathbf{S})$ is locally rigid, and the converse holds if $\Theta(\mathbf{S})$ is generic.
\end{thm}

\revadd{Moreover, using Corollary~\ref{cor:suff_non_deg_loc_min}, Proposition~\ref{prop:noiseless_setting1}, we obtain a sufficient condition for any realization of a patch framework to be locally rigid.
\begin{cor}
\label{cor:noiseless_setting1}
If $\mathbf{C}$ is of rank $(m-1)d$ then every perfect alignment $\mathbf{S}$ of $F$ is non-degenerate and the realization $\Theta(\mathbf{S})$ is locally rigid.
\end{cor}
}

Using the Definition~\ref{def:uniq_alignment} and~\ref{def:glob_rigid}, it is easy to deduce that a unique perfect alignment results in a globally rigid realization and vice versa. Then, using the fact that affine rigidity implies global rigidity \cite{gortler2010affine,toth2017handbook}, and a realization is affinely rigid if and only if the rank of $\mathbf{C}$ is $(m-1)d$ \citea{chaudhury2015global}, it follows from Corollary~\ref{cor:noiseless_setting1} that the unique perfect alignment underlying an affinely rigid realization is also non-degenerate.
\begin{thm}
\label{thm:glob_rigid}
Let $\mathbf{S}$ be a perfect alignment.
Then $\Theta(\mathbf{S})$ is globally rigid iff $\mathbf{S}$ is unique.
\end{thm}
\begin{prop}
\label{prop:affine_rigid}
\revadd{Let $\mathbf{S}$ be a perfect alignment.
If the realization $\Theta(\mathbf{S})$ is affinely rigid then $\mathbf{S}$ is a non-degenerate and unique perfect alignment. The converse does not hold due to the counterexample in \cite[Figure~3]{gortler2010affine}.}
\end{prop}
\revadd{
Finally, combining the affine rigidity rank condition ($\rank(\mathbf{C}) = (m-1)d$) with Corollary~\ref{cor:non_deg_d_2}, Remark~\ref{rmk:C1Sp} and Theorem~\ref{thm:loc_rigid}, we are also able to connect the local rigidity of a realization in two dimensions with affine rigidity of its projection in one dimension.
\begin{cor}
\label{cor:local_affine_rigid_in_d_2}
Let $d=2$ and $\mathbf{S}$ be a perfect alignment. Then $\Theta(\mathbf{S})$ is locally rigid if its projection in at least one of the two dimensions is affinely rigid.
\end{cor}}
\subsection{Conditions on Overlapping Views for a Non-degenerate Perfect Alignment}
\label{subsec:non_deg_noiseless_setting}
\revadd{We now focus on deriving the necessary and sufficient conditions on the overlapping structure of the views for a perfect alignment to be non-degenerate. These are inspired by the affine rigidity criteria discussed in~\cite{zha2009spectral}. The main difference is that we impose relatively weaker rank constraint on the overlaps.
In fact, as indicated by our previous results (Figure~\ref{fig:rigidity_flow}), the conditions presented here can be viewed as those ensuring infinitesimal and generic local rigidity of a realization.} To begin with,

\begin{dfn}
\label{def:BSAcapB}
Let $\mathbf{S}$ be a perfect alignment. Let $A$ and $B$ be non-empty disjoint subsets of $[1,m]$. Define $\mathbf{B}(\mathbf{S})_{A,B}$ to be a matrix whose columns are $\mathbf{S}_i^T\mathbf{x}_{k,i}+\mathbf{t}_i$ (in the increasing order of $k$) where $(k,i),(k,j) \in E(\Gamma)$ for some $i \in A$ and $j \in B$, and where $\mathbf{t}_i$ is obtained using Eq.~(\ref{eq:opt_Z}). Also define $\overline{\mathbf{B}(\mathbf{S})}_{A,B} = \mathbf{B}(\mathbf{S})_{A,B}\left(\mathbf{I}_{n'} - (1/n')\mathbf{1}_{n'}\mathbf{1}_{n'}^T\right)$ where $n' = |\{k:(k,i),(k,j) \in E(\Gamma) \text{ for some } (i,j) \in A \times B\}|$.
For brevity, we denote $\mathbf{B}(\mathbf{S})_{\{i\},\{j\}}$ and $\overline{\mathbf{B}(\mathbf{S})}_{\{i\},\{j\}}$ by $\mathbf{B}(\mathbf{S})_{i,j}$ and $\overline{\mathbf{B}(\mathbf{S})}_{i,j}$ respectively, where $i \neq j$. Note that the notation is consistent with that of Definition~\ref{def:BSicapj}.
\end{dfn}

\begin{rmk}
\label{rmk:BS_ijB_ij_noiseless}
Since $\mathbf{S}$ is a perfect alignment $\mathbf{S}_i^T\mathbf{x}_{k,i}+\mathbf{t}_i = \mathbf{S}_j^T\mathbf{x}_{k,j}+\mathbf{t}_j$ for all $(k,i),(k,j) \in E(\Gamma)$, thus $\mathbf{B}(\mathbf{S})_{A,B}$ is well defined and $\mathbf{B}(\mathbf{S})_{A,B} = \mathbf{B}(\mathbf{S})_{B,A}$.
Let $i,j \in [1,m]$ then, since $\mathbf{B}(\mathbf{S})_{i,j} = \mathbf{B}(\mathbf{S})_{j,i}$, from Remark~\ref{rmk:BS_ijB_ij}, $\rank (\overline{\mathbf{B}}_{i,j}) = \rank (\overline{\mathbf{B}(\mathbf{S})}_{i,j}) = \rank (\overline{\mathbf{B}(\mathbf{S})}_{j,i}) = \rank (\overline{\mathbf{B}}_{j,i})$ and $\rank (\overline{\mathbf{B}}_{i,j}) = \rank (\overline{\mathbf{B}(\mathbf{S})}_{i,j}\overline{\mathbf{B}(\mathbf{S})}_{j,i}^T) = \rank (\overline{\mathbf{B}}_{i,j}\overline{\mathbf{B}}_{j,i}^T)$.
\end{rmk}

Due to the above two remarks and Theorem~\ref{thm:non_deg_two_views_gen_setting}, a necessary and sufficient condition for a perfect alignment of two views to be non-degenerate is easily obtained.
\begin{thm}
\label{thm:nec_suff_cond_loc_rigid_two_views}
Consider $m=2$ and let $\mathbf{S}$ be a perfect alignment. Then $\mathbf{S}$ is non-degenerate iff $\rank(\overline{\mathbf{B}}_{1,2}) \geq d-1$.
\end{thm}
\begin{figure}[H]
    \centering
    \begin{tabular}{cc}
    \begin{subfigure}[b]{0.35\textwidth}
         \centering
         \includegraphics[width=0.9\textwidth,keepaspectratio]{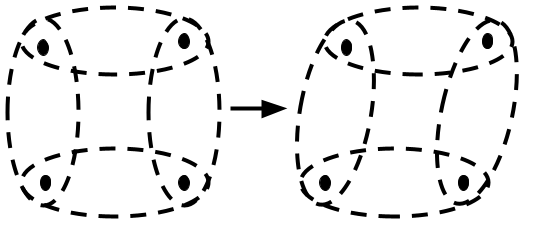}         \caption{Theorem~\ref{thm:nec_cond_loc_rigid_of_views}}
         \label{fig:nec_cond_loc_rigid_of_views}
     \end{subfigure}
     & 
     \begin{subfigure}[b]{0.175\textwidth}
         \centering
         \includegraphics[width=0.9\textwidth,keepaspectratio]{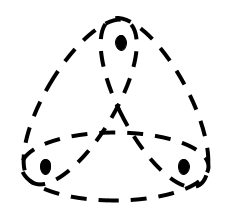}
         \caption{Theorem~\ref{thm:G_star_1}}
         \label{fig:G_star_1}
     \end{subfigure}
     \end{tabular}
    \caption{Counterexamples for the converse of various Theorems. (a) For every pair of nonempty partitions $A$ and $B$ of $[1,4]$, $\rank(\overline{\mathbf{B}(\mathbf{S})}_{A,B}) \geq 1$ but $\mathbf{S}$ is degenerate. (b) $\mathbf{S}$ is non-degenerate but $|\mathbb{G}^*(\mathbf{S})| = 3$. These views should be considered affinely non-degenerate, satisfying Assumption~\ref{assump:non_deg_views}. For clarity, only the points in the overlapping regions are shown.}
    \label{fig:counterex}
\end{figure}
A necessary condition for a perfect alignment of $m \geq 3$ views to be non-degenerate is as follows. The converse of the theorem does not hold, as demonstrated in Figure~\ref{fig:nec_cond_loc_rigid_of_views}.
\begin{thm}
\label{thm:nec_cond_loc_rigid_of_views}
Let $\mathbf{S}$ be a perfect alignment. If $\mathbf{S}$ is non-degenerate then the $\rank(\overline{\mathbf{B}(\mathbf{S})}_{A,B})$ is at least $d-1$ for all non-empty partitions $A$ and $B$ of $[1,m]$ i.e. for all $A,B \subseteq [1,m]$, $A,B \neq \emptyset$, $A \cap B = \emptyset$ and $A \cup B = [1,m]$.
\end{thm}

Now we derive a sufficient condition for a perfect alignment of $m \geq 3$ views to be non-degenerate. As in \citea{zha2009spectral}, we construct a graph $\mathbb{G}$ with $m$ vertices where each vertex corresponds to a view and an edge exists between the $i$th and $j$th vertices iff $\rank(\overline{\mathbf{B}}_{i,j}) \geq d-1$. The Theorem~\ref{thm:nec_suff_cond_loc_rigid_two_views} and the following propositions will play a crucial role in our next set of results,

\begin{lem}
\label{lem:subproblem_cert}
Let $\mathbf{S}$ be a perfect alignment and $\boldsymbol{\Omega}$ be a certificate of $\mathbf{L}(\mathbf{S})$. Consider removing the $i$th view and the points that lie exclusively in it. Then $\mathbf{S}_{-i} = [\mathbf{S}_j]_{j \in [1,m] \setminus \{i\}}$ is a perfect alignment of the remaining views and $[\boldsymbol{\Omega}_j]_{j \in [1,m] \setminus \{i\}}$ is a certificate of $\mathbf{L}_{-i}(\mathbf{S}_{-i})$, the matrix in Eq.~(\ref{eq:L_of_S}) associated with the remaining views.
\end{lem}

\begin{prop}
\label{prop:same_conn_comp_non_deg}
Let $\mathbf{S}$ be a perfect alignment. Let $\boldsymbol{\Omega}$ be a certificate of $\mathbf{L}(\mathbf{S})$. If $i$th and $j$th view lie in the same connected component of $\mathbb{G}$ then $\boldsymbol{\Omega}_i = \boldsymbol{\Omega}_j$.
\end{prop}

Similar to \citea{zha2009spectral}, consider the following coarsening procedure on $\mathbb{G}$ given a perfect alignment $\mathbf{S}$: (i) transform all the views using $\mathbf{S}$ (and $\mathbf{t}$ computed using Eq.~\ref{eq:opt_Z}), (ii) merge the views that lie in the same connected component of $\mathbb{G}$ and replace them with a single view, (iii) then construct the graph (in the same manner as $\mathbb{G}$) associated with the new set of views, (iv) repeat the procedure from (ii). Let the final graph over the remaining views be $\mathbb{G}^*(\mathbf{S})$, then the following result holds (the corollary follows trivially and the converse of the theorem may not hold, as shown in Figure~\ref{fig:G_star_1}).
\begin{thm}
\label{thm:G_star_1}
    A perfect alignment $\mathbf{S}$ is non-degenerate if $|\mathbb{G}^*(\mathbf{S})| = 1$.
\end{thm}
\begin{cor}
\label{cor:suff_cond_views_non_deg}
Every perfect alignment is non-degenerate if $\mathbb{G}$ is connected.
\end{cor}
\begin{rmk}
\label{rmk:tree_structure}
\revadd{From Theorem~\ref{thm:nec_cond_loc_rigid_of_views}, it is easy to see that the converse of the above corollary holds when a graph over views, in which two views are connected if they are overlapping (i.e. they share at least one common point), is a tree.}
\end{rmk}

\revadd{It is important to note that we have derived a necessary and sufficient rank based condition (Proposition~\ref{prop:noiseless_setting1}) that can be tested in polynomial time to assess the non-degeneracy of a given perfect alignment. While the above results offer a geometric interpretation of this condition, a complete understanding in the form of a single condition on the overlapping structure of the views, that is both necessary and sufficient, has yet to be established.}

\subsection{Conditions on Overlapping Views for a Unique Perfect Alignment}
\label{subsec:uniq_noiseless_setting}
\revdel{We previously showed that the global rigidity of a realization is equivalent to the uniqueness of the corresponding perfect alignment (Theorem~\ref{thm:glob_rigid}). Thus}Here, we focus on deriving necessary and sufficient conditions on the overlapping structure of the views for a perfect alignment to be unique, equivalently, for the resulting realization to be globally rigid. From Remark~\ref{rmk:BS_ijB_ij_noiseless} and Theorem~\ref{thm:uniq_two_views_gen_setting},
\begin{thm}
\label{thm:nec_suff_cond_glob_rigid_two_views}
Consider $m=2$ and let $\mathbf{S}$ be a perfect alignment. Then $\mathbf{S}$ is unique (Definition~\ref{def:uniq_alignment}) iff $\rank (\overline{\mathbf{B}}_{1,2}) = d$.
\end{thm}

A necessary condition for a perfect alignment of $m \geq 3$ views to be unique is, 
\begin{thm}
\label{thm:nec_cond_glob_rigid_views}
If $\mathbf{S}$ is a unique perfect alignment then $\rank(\overline{\mathbf{B}(\mathbf{S})}_{A,B}) = d$ for all non-empty partitions $A$ and $B$ of $[1,m]$.
\end{thm}

\revadd{It was shown in \citea{zha2009spectral} that the above rank condition holds for affinely rigid realization $\Theta(\mathbf{S})$. In contrast, our requirement only requires the perfect alignment $\mathbf{S}$ to be unique, which is equivalent to the global rigidity of $\Theta(\mathbf{S})$ (Figure~\ref{fig:rigidity_flow}). Moreover, we conjecture that the converse of the above theorem holds, in which case we would obtain a characterization of a unique perfect alignment and an exponential-time algorithm to test it, aligning with the NP-hardness of testing global rigidity~\cite{saxe1979embeddability}.}


Now we derive a sufficient condition for a perfect alignment of $m \geq 3$ to be unique. As in the previous section, we construct a graph $\overline{\mathbb{G}}$ with $m$ vertices, one for each view. An edge exists between the $i$th and $j$th vertices iff $\rank(\overline{\mathbf{B}}_{i,j}) = d$. We need Lemma~\ref{lem:subproblem_cert}, Theorem~\ref{thm:nec_suff_cond_glob_rigid_two_views} and the following proposition for our next result,
\begin{prop}
\label{prop:same_conn_comp_uniq}
Let $\mathbf{S}$ and $\mathbf{S}'$ be perfect alignments. If $i$th and $j$th view lie in the same connected component of $\overline{\mathbb{G}}$ then $\mathbf{S}'_i = \mathbf{S}_i\mathbf{Q}$ and $\mathbf{S}'_j = \mathbf{S}_j\mathbf{Q}$ for some $\mathbf{Q} \in \mathbb{O}(d)$.
\end{prop}
\begin{figure}[H]
    \centering
     \includegraphics[width=0.15\textwidth,keepaspectratio]{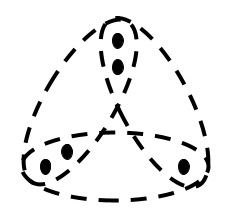}
    \caption{$\mathbf{S}$ is unique but $|\overline{\mathbb{G}}^*(\mathbf{S})| = 3$, thus the converse of Theorem~\ref{thm:overline_G_star_1} may not hold.}
    \label{fig:overline_G_star_1}
\end{figure}

Consider the same coarsening procedure as in Theorem~\ref{thm:G_star_1}, except that $\mathbb{G}$ and $\mathbb{G}^*(\mathbf{S})$ are replaced by $\overline{\mathbb{G}}$ and $\overline{\mathbb{G}}^*(\mathbf{S})$, respectively. Then the following holds (the corollary follows trivially and a counterexample for the converse is shown in Figure~\ref{fig:overline_G_star_1}).
\begin{thm}
\label{thm:overline_G_star_1}
    A perfect alignment $\mathbf{S}$ is unique if $|\overline{\mathbb{G}}^*(\mathbf{S})| = 1$.
\end{thm}
\begin{cor}
\label{cor:suff_cond_views_uniq}
Every perfect alignment is unique if $\overline{\mathbb{G}}$ is connected.
\end{cor}
\revadd{We note that the above result is weaker than the one in \citea{zha2009spectral} where the authors showed that the sufficient condition leads to affinely rigid realization $\Theta(\mathbf{S})$, while we show global rigidity. Nevertheless, we keep the result since the the proving technique is different than the one in \citea{zha2009spectral}.}

\section{Linear Convergence of RGD}
\label{sec:convergence}
In this section, we describe the RGD algorithm for solving the alignment problem in Eq.~(\ref{eq:GPOP}). Using the theory of Morse functions, we show that if the sequence of iterates generated by RGD converges to a non-degenerate alignment then the convergence is linear. Moreover, We obtain an estimate of the radius and the rate of convergence. We also present an exact recovery and noise stability analysis of RGD, initialized using the output of the spectral algorithm (SPEC) \citea{chaudhury2015global}.

\subsection{RGD Algorithm}
\label{subsec:rgd_algo}
A standard way to find a local minimum of Eq.~(\ref{eq:GPOP}) is to use RGD with a suitable initial point, step size and retraction strategy. In this work we use retraction based on the exponential map on $\mathbb{O}(d)$ \citea{van1996matrix, udriste2013convex}. Define,
\begin{align}
    R_{\EXP }: \cup_{\mathbf{S} \in \mathbb{O}(d)^m}(\{\mathbf{S}\} \times T_\mathbf{S}\mathbb{O}(d)^m) &\mapsto \mathbb{O}(d)^m\\
    R_{\EXP }\left([\mathbf{S}_i]_1^m, [\boldsymbol{\xi}_i]_1^m\right) &= [\mathbf{S}_i\exp(\mathbf{S}_i^T\boldsymbol{\xi}_i)]_1^m. \label{eq:R_PF}
\end{align}
where $\exp (\mathbf{A})$ denotes the matrix exponential of $\mathbf{A}$ \citea{van1996matrix, absil2009optimization}. Then the following lemma provides a consistent definition of a retraction on the quotient manifold $\mathbb{O}(d)^m/_{\sim}$.
\begin{lem}
\label{lem:retraction}
Let $\widetilde{\mathbf{S}} \in \mathbb{O}(d)^{m}/_{\sim}$ and $\mathbf{S}^a, \mathbf{S}^b \in \pi^{-1}(\widetilde{\mathbf{S}})$. If $\mathbf{Z}^a \in T_{\mathbf{S}^a}\mathbb{O}(d)^m$ and $\mathbf{Z}^b \in T_{\mathbf{S}^b}\mathbb{O}(d)^m$ are the horizontal lifts of $\widetilde{\mathbf{Z}} \in T_{\widetilde{\mathbf{S}}}\mathbb{O}(d)^{m}/_{\sim}$ then $\pi(R_{\EXP }(\mathbf{S}^a, \mathbf{Z}^a)) = \pi(R_{\EXP }(\mathbf{S}^b, \mathbf{Z}^b))$ (see Eq.~\ref{eq:pi}). As a result, the retraction
\begin{align}
    \widetilde{R}_{\EXP }: \cup_{\widetilde{\mathbf{S}} \in \mathbb{O}(d)^m/_{\sim} }(\{\widetilde{\mathbf{S}}\} \times T_{\widetilde{\mathbf{S}}}\mathbb{O}(d)^m/_{\sim}) &\mapsto \mathbb{O}(d)^m/_{\sim}\\
    \widetilde{R}_{\EXP }\left(\widetilde{\mathbf{S}}, \widetilde{\mathbf{Z}}\right) &=  \pi(R_{\EXP }(\mathbf{S}, \mathbf{Z}))\label{eq:Rtilde_PF}
\end{align}
is well defined for any $\mathbf{S} \in \pi^{-1}(\widetilde{\mathbf{S}})$ and $\mathbf{Z}$ being the horizontal lift of $\widetilde{\mathbf{Z}}$ at $\mathbf{S}$.
\end{lem}

The step direction will always be the horizontal lift of $-\grad \widetilde{F}(\widetilde{\mathbf{S}})$ at some $\mathbf{S} \in \pi^{-1}(\widetilde{\mathbf{S}})$. Consequently, due to Proposition~\ref{prop:gradFS}, the step direction is $\boldsymbol{\xi} = -\grad F(\mathbf{S}) = [[\mathbf{C}\mathbf{S}]_i - \mathbf{S}_i[\mathbf{C}\mathbf{S}]_i^T\mathbf{S}_i]_1^m$, the projection of the antigradient $-\nabla F(\mathbf{S})$ onto $T_\mathbf{S}\mathbb{O}(d)^m$. The step size $\alpha$ is calculated using the Armijo-type rule with parameters $\beta,\gamma \in (0,1)$ (here $g$ is the canonical metric on $\mathbb{O}(d)^m$ as in Eq.~(\ref{eq:g_Z_W})),
\begin{equation}
    \alpha = \max_{l \geq 0}\{\beta^l\ \vertbar\ F(R_{\EXP }(\mathbf{S}, -\beta^l\grad F(\mathbf{S}))) - F(\mathbf{S}) \leq -\gamma \beta^l g(\nabla F(\mathbf{S}),  \grad F(\mathbf{S})) \}. \label{eq:armijo_step}
\end{equation}
Since $F$ extends to a continuously differentiable non-negative function on $\mathbb{R}^{md \times d}$ containing $\mathbb{O}(d)^m$, it follows from \citeb[Proposition 2.8]{schneider2015convergence} that $\alpha$ is well-defined.
\begin{algorithm}[H]
\caption{Riemannian gradient descent for solving GPOP \label{algo:rgd}}
\begin{algorithmic}[1]
\REQUIRE $\widetilde{\mathbf{S}}^0 \in \mathbb{O}(d)^{m-1}$, $\Gamma$, $\{\mathbf{x}_{k,i}: (k,i) \in E(\Gamma)\}$, $\beta, \gamma \in (0,1)$
\STATE Construct $\mathbf{C}$ as in Eq.~(\ref{eq:GPOP}).
\REPEAT
    \STATE set $\mathbf{S}^k = [\mathbf{I}_d; \widetilde{\mathbf{S}}^k] \in \pi^{-1}(\widetilde{\mathbf{S}}^k) \subset \mathbb{O}(d)^m$ (Eq.~\ref{eq:pi_inv_wtS}).
    \STATE calculate the descent direction $-\grad F(\mathbf{S}^k)$ at $\mathbf{S}^k$ using Eq.~(\ref{eq:gradFS}).
    \STATE calculate the step size $\alpha_k$ according to the Armijo-type rule (see Eq.~(\ref{eq:armijo_step})).
    \STATE set $\widetilde{\mathbf{S}}^{k+1} = \pi(R_{\EXP}(\mathbf{S}^k, -\alpha_k \grad F(\mathbf{S}^k)))$ using Eq.~(\ref{eq:R_PF}, \ref{eq:pi}).
    \STATE $k \leftarrow k + 1$.
\UNTIL{convergence.}
\end{algorithmic}
\end{algorithm}
\subsection{Local linear Convergence of RGD}
\label{subsec:loc_lin_conv}
We proceed to show the local linear convergence of Algorithm~\ref{algo:rgd} to a non-degenerate alignment using the convergence analysis framework presented in \citea{schneider2015convergence} and as used in \citea{liu2019quadratic}.
To this end, we note that $\widetilde{F}$ and $F$ are real-analytic functions bounded from below by zero. $\mathbb{O}(d)^m/_{\sim}$  (whose elements are identified with $\mathbb{O}(d)^{m-1}$ here) and $\mathbb{O}(d)^m$ are compact submanifolds of $\mathbb{R}^{(m-1)d \times d}$ and $\mathbb{R}^{md \times d}$, respectively.
However, since $F(\mathbf{S}\mathbf{Q}) = F(\mathbf{S})$ for all $\mathbf{Q} \in \mathbb{O}(d)$, therefore every critical point of $F$ is degenerate and in particular $F$ is not a Morse-function \citea{cohen_iga_norbury_2006}. Nevertheless, if $\mathbf{S}^*$ is a non-degenerate alignment then $\widetilde{\mathbf{S}}^* = \pi(\mathbf{S}^*)$ is a non-degenerate critical point of $\widetilde{F}$. As a result $\widetilde{F}$ is a Morse function at $\widetilde{\mathbf{S}}^*$ and, due to \citeb[Proposition 4.2]{hu2018convergence}, the Lojasiewicz gradient inequality is satisfied.
\begin{prop}
Let $\mathbf{S}^*$ be a non-degenerate alignment and define $\widetilde{\mathbf{S}}^* = \pi(\mathbf{S}^*)$. Then there exist $\delta, \eta > 0$ such that $$|\widetilde{F}(\widetilde{\mathbf{S}}) - \widetilde{F}(\widetilde{\mathbf{S}}^*)| \leq \eta \left\| \grad \widetilde{F}(\widetilde{\mathbf{S}})\right\|_F^2$$
holds for every $\widetilde{\mathbf{S}} \in \mathbb{O}(d)^m/_{\sim}$ satisfying $\left\|\widetilde{\mathbf{S}}-\widetilde{\mathbf{S}}^*\right\|_F < \delta$.
\end{prop}

Moreover, the iterates $\{\widetilde{\mathbf{S}}^k\}_{k \geq 0}$ generated by Algorithm~\ref{algo:rgd} satisfy the (\textbf{A1}) sufficient descent, (\textbf{A2}) stationarity and (\textbf{A3}) safeguard assumptions below. The proofs are in the appendix and we make use of the following results to prove them.
\begin{prop}
\label{prop:liu_pf}
For all $\mathbf{S}_i \in \mathbb{O}(d)$ and $\mathbf{Z}_i \in T_{\mathbf{S}_i}\mathbb{O}(d)$ satisfying $\left\|\mathbf{Z}_i\right\|_F \leq 1$,
$$\left\|\mathbf{S}_i\exp (\mathbf{S}_i^T\mathbf{Z}_i) - (\mathbf{S}_i + \mathbf{Z}_i)\right\|_F \leq (e-1)\left\|\mathbf{Z}_i\right\|_F^2.$$
\end{prop}
\begin{prop}
\label{prop:second_order_boundedness_of_Rtilde}
For $\widetilde{\mathbf{S}} \in \mathbb{O}(d)^m/_{\sim}$ and $\widetilde{\mathbf{Z}} \in T_{\widetilde{\mathbf{S}}}\mathbb{O}(d)^m/_{\sim}$ satisfying  $\left\|\widetilde{\mathbf{Z}}\right\|_F \leq 1/2$,
\begin{enumerate}[leftmargin=*,label=(\alph*)]
    \item $\left\|R_\EXP(\mathbf{S}, \mathbf{Z}) - (\mathbf{S} + \mathbf{Z})\right\|_F \leq (e-1)\left\|\mathbf{Z}\right\|_F^2$ for any $\mathbf{S} \in \pi^{-1}(\widetilde{\mathbf{S}})$ and $\mathbf{Z} \in T_{\mathbf{S}}\mathbb{O}(d)^m$, the horizontal lift of $\widetilde{\mathbf{Z}}$ at $\mathbf{S}$. 
    \item $\left\|\widetilde{R}_\EXP(\widetilde{\mathbf{S}}, \widetilde{\mathbf{Z}}) - (\widetilde{\mathbf{S}} + \widetilde{\mathbf{Z}})\right\|_F \leq (e-1)\left\|\widetilde{\mathbf{Z}}\right\|_F^2$.
\end{enumerate}
\end{prop}
\begin{prop}
\label{prop:alpha_grad}
$\lim \alpha_k \left\|\grad \widetilde{F}(\widetilde{\mathbf{S}}^k)\right\|_F = 0$ and $\lim \alpha_k \left\|\grad F(\mathbf{S}^k)\right\|_F = 0$.
\end{prop}

\noindent \textbf{(A1)}. \textit{(Sufficient Descent)} There exist $\kappa_0 > 0$ and $k_1 \in \mathbb{N}$ such that, the inequality $\widetilde{F}(\widetilde{\mathbf{S}}^{k+1}) - \widetilde{F}(\widetilde{\mathbf{S}}^k) \leq - \kappa_0 \left\|\grad \widetilde{F}(\widetilde{\mathbf{S}}^k)\right\|_F \cdot \left\|\widetilde{\mathbf{S}}^{k+1}-\widetilde{\mathbf{S}}^k\right\|_F$ holds for all $k \geq k_1$.
\smallskip

\noindent \textbf{(A2)}. \textit{(Stationarity)} There exist $k_2 \in \mathbb{N}$ such that for all $k \geq k_2$, if $\left\|\grad \widetilde{F}(\widetilde{\mathbf{S}}^k)\right\|_F = 0$ then $\widetilde{\mathbf{S}}^{k+1} = \widetilde{\mathbf{S}}^k$. The sequence $\{\widetilde{\mathbf{S}}^{k}\}_{k \geq 0}$ satisfies this trivially.
\smallskip

\noindent \textbf{(A3)}. \textit{(Safeguard)} There exist a constant $\mu > 0$ and $k_3 \in \mathbb{N}$ such that the inequality $\left\|\grad \widetilde{F}(\widetilde{\mathbf{S}}^k)\right\|_F \leq \mu \left\|\widetilde{\mathbf{S}}^{k+1}-\widetilde{\mathbf{S}}^k\right\|_F$ holds for all $k \geq k_3$.

Combined with Theorem 2.3 in \citea{schneider2015convergence} and the fact that $\mathbb{O}(d)^m/_{\sim}$ is compact (thus every sequence on it has a cluster point), we obtain the following result.
\begin{thm}
\label{thm:rgd_conv}
Let $\mathbf{S}^*$ be a non-degenerate alignment and $\widetilde{\mathbf{S}}^* = \pi(\mathbf{S}^*)$. If the sequence $\{\widetilde{\mathbf{S}}^k\}_{k \geq 0}$ due to Algorithm~\ref{algo:rgd} converges to $\widetilde{\mathbf{S}}^*$ then the convergence is linear.
\end{thm}

Finally, we obtain an estimate of the radius and rate of linear convergence. The proof follows directly from \citeb[Chapter 7, Theorem~4.2]{udriste2013convex} (here $d_{\widetilde{g}}$ is the geodesic distance induced by the metric $\widetilde{g}$ on $\mathbb{O}(d)^m/_\sim$ as defined in Proposition~\ref{prop:g_tilde}).
\begin{thm}
\label{thm:rgd_conv2}
Let $\mathbf{S}^*$ be a non-degenerate alignment
and $\zeta \in (0,1)$ be fixed. Let $\lambda_{-}(\mathbf{S}^*)$, $\lambda_{+}(\mathbf{S}^*)$ and $\delta(\mathbf{S}^*)$ be as defined in Proposition~\ref{prop:HessVicinity}. If  the initialization $\widetilde{\mathbf{S}}^0 = \pi(\mathbf{S}^0)$ of Algorithm~\ref{algo:rgd} and and the subsequent iterates $\widetilde{\mathbf{S}}^k = \pi(\mathbf{S}^k)$ generated by it satisfy $\min_{\mathbf{Q}\in \mathbb{O}(d)}\left\|\mathbf{S}^k-\mathbf{S}^*\mathbf{Q}\right\|_F < \min\left\{2,\frac{2}{\pi}\zeta\delta(\mathbf{S}^*)\right\}$, then the sequence $\{\widetilde{\mathbf{S}}^k\}_{k \geq 0}$ converges to $\widetilde{\mathbf{S}}^* = \pi(\mathbf{S}^*)$ linearly. Moreover,
\begin{align}
    \widetilde{F}(\widetilde{\mathbf{S}}^k) - \widetilde{F}(\widetilde{\mathbf{S}}^*) &\leq q^{k}(\widetilde{F}(\widetilde{\mathbf{S}}^0)-\widetilde{F}(\widetilde{\mathbf{S}}^*))\label{eq:alignment_err_ratio}\\
    d_{\widetilde{g}}(\widetilde{\mathbf{S}}^k, \widetilde{\mathbf{S}}^*) &\leq C q^{(k-1)/2}
\end{align}
where $C > 0$ is a constant, $q = 1 - 2\gamma (1-\gamma) r(1+r)\in (0,1)$ and $r = \frac{(1-\zeta)\lambda_{-}(\mathbf{S}^*)}{\lambda_{+}(\mathbf{S}^*) +\zeta\lambda_{-}(\mathbf{S}^*)}$.
\end{thm}

\subsection{Exact Recovery and Noise Stability}
\label{subsec:noise_stability}
A direct consequence of Theorem~\ref{thm:rgd_conv2} and Corollary~\ref{cor:HessVicinity} to the noiseless setting is that RGD converges locally linearly to a perfect alignment under a condition weaker than affine rigidity.
\begin{thm}
\label{thm:exact_recovery}
Suppose $\mathbf{S}^*$ is a perfect alignment. Then Algorithm~\ref{algo:rgd} converges locally linearly to $\widetilde{\mathbf{S}}^* = \pi(\mathbf{S}^*)$ if any of the following holds:
\begin{enumerate}
    \item $\rank(\mathbf{C}) = (m-1)d$.
    \item $\rank(\mathbb{L}(\mathbf{S}^*))  = (m-1)d(d-1)/2$.
\end{enumerate}
The first condition which characterizes affine rigidity implies the second that characterizes infinitesimal rigidity as well as generic local rigidity (Figure~\ref{fig:rigidity_flow}).
\end{thm}

While RGD achieves local linear convergence under less restrictive conditions, a key challenge lies in selecting an initial alignment that is sufficiently close to a non-degenerate perfect alignment $\mathbf{S}^*$. In \citea{chaudhury2015global}, the authors showed that under the affine rigidity constraints, SPEC recovers the perfect alignment, eliminating the need for RGD. However, under weaker non-degeneracy (equivalently, infinitesimal/local rigidity) constraints, it is yet to be established whether the output of SPEC recovers/remains close to a perfect alignment of noiseless views (similarly, to an optimal alignment of the noisy views). We aim to address this in our future work.

Nevertheless, under the bounded noise model and affine rigidity constraints, the spectral solution $\mathbf{S}_{spec}$ has been shown to approximate a perfect alignment $\mathbf{S}_0$ of the noiseless counterparts of the noisy views \citea{chaudhury2015global}. While $\mathbf{S}_0$ is generally not the optimal alignment $\mathbf{S}^*$ of the noisy views, one can expect them to be relatively close. Therefore, under affine rigidity constraints, we expect $\mathbf{S}_{spec}$ to be near the optimal alignment $\mathbf{S}^*$ of the noisy views and thus, refining $\mathbf{S}_{spec}$ using RGD could potentially yield $\mathbf{S}^*$. Here, we provide a noise stability analysis of RGD which support this idea.


We start with a set of noiseless views and inject them with bounded noise i.e.  $\mathbf{x}_{k,i} \leftarrow \mathbf{x}_{k,i} + \boldsymbol{\epsilon}_{k,i}$ where $\left\|\boldsymbol{\epsilon}_{k,i}\right\|_2 \leq \varepsilon$ for a fixed noise level $\varepsilon > 0$. Let $\mathbf{C}_0$ and $\mathbf{C}$ be the patch-stress matrices (Eq.~(\ref{eq:GPOP})) corresponding to the noiseless views and their noisy counterparts, respectively.
Then the following lemma establishes a quadratic growth condition at an optimal alignment in the noiseless setting. The subsequent lemma bounds the distance between the optimal alignments of noisy and noiseless views.
\begin{lem}
\label{lem:quadgrowth}
Let $\rank(\mathbf{C}_0) = (m-1)d$, and consequently $\mathbf{S}_0$ be a unique perfect alignment of the noiseless views (Figure~\ref{fig:rigidity_flow}). Then,
\begin{equation}
    \Tr(\mathbf{C}_0\mathbf{S}\mathbf{S}^T) \geq \frac{\lambda_{d+1}(\mathbf{C}_0)}{2} \min_{\mathbf{Q} \in \mathbb{O}(d)}\left\|\mathbf{S}- \mathbf{S}_0\mathbf{Q}\right\|_F^2.
\end{equation}
\end{lem}
\begin{lem}
\label{lem:distS_0Sstar}
Let $\rank(\mathbf{C}_0) = (m-1)d$ and $\mathbf{S}_0$ be a unique perfect alignment of the noiseless views. Let $\mathbf{S}^*$ be an optimal alignment of the noisy views. Then
\begin{equation}
    \min_{\mathbf{Q}\in\mathbb{O}(d)}\left\|\mathbf{S}^* - \mathbf{S}_0\mathbf{Q}\right\|_F \leq \frac{4 m \left\|\mathbf{C}-\mathbf{C}_0\right\|_F}{\lambda_{d+1}(\mathbf{C}_0)}.
\end{equation}
\end{lem}

Finally, we obtain a bound on the noise level for RGD, initialized with  the spectral alignment, to converge locally linearly to the optimal alignment of the noisy views. 
\begin{thm}
\label{thm:rgd_noise_stability}
Let $\rank(\mathbf{C}_0) = (m-1)d$ and $\mathbf{S}_0$ be a unique perfect alignment of the noiseless views. Let $\mathbf{S}^*$ be an optimal alignment of the noisy views and suppose $\mathbf{S}^*$ is non-degenerate. Let $\zeta \in (0,1)$ be fixed. Then Algorithm~\ref{algo:rgd}, initialized with $\pi(\mathbf{S}_{spec}(\mathbf{C}))$, converges locally linearly to $\pi(\mathbf{S}^*)$ if the noise level $\varepsilon$ satisfies
\begin{equation}
    4\sqrt{m}\left(\frac{\pi\sqrt{d(d+1)}}{\lambda_{d+1}(\mathbf{C})} + \frac{\sqrt{m}}{\lambda_{d+1}(\mathbf{C}_0)}\right)(K_1 \varepsilon + K_2\varepsilon^2) < \min\left\{2, \frac{2}{\pi}\zeta\delta(\mathbf{S}^*)\right\}
\end{equation}
and the subsequent iterates satisfy $\min_{\mathbf{Q}\in \mathbb{O}(d)}\left\|\mathbf{S}^k-\mathbf{S}^*\mathbf{Q}\right\|_F < \min\left\{2, \frac{2}{\pi}\zeta\delta(\mathbf{S}^*)\right\}$.
Here,
\begin{align}
    K_1 &= 2\sqrt{n |E(\Gamma)|}\left(4 \max_1^n\left\|\mathbf{x}_k^*\right\|_2\frac{\sqrt{n|E(\Gamma)|}}{\lambda_{2}(\mathbf{\mathcal{L}}_{\Gamma})} + 1\right)\\
    K_2 &= 2\sqrt{n |E(\Gamma)|}\left(2\frac{\sqrt{n|E(\Gamma)|}}{\lambda_{2}(\mathbf{\mathcal{L}}_{\Gamma})} + 1\right),
\end{align}
$\mathbf{x}_k^*$ is the realization of the noiseless views due to $\mathbf{S}_0$ and $\delta(\mathbf{S}^*)$ is defined in Proposition~\ref{prop:HessVicinity}.
\end{thm}
In the result above, we assumed $\lambda_{d+1}(\mathbf{C}) > 0$, as was the case in \citea{chaudhury2015global}, where the authors observed (and as we validate below) that $\lambda_{d+1}(\mathbf{C})$ increases with noise level $\varepsilon$.
\begin{figure}
    \centering
    \begin{tabular}{cc}
       \includegraphics[width=0.46\linewidth,keepaspectratio]{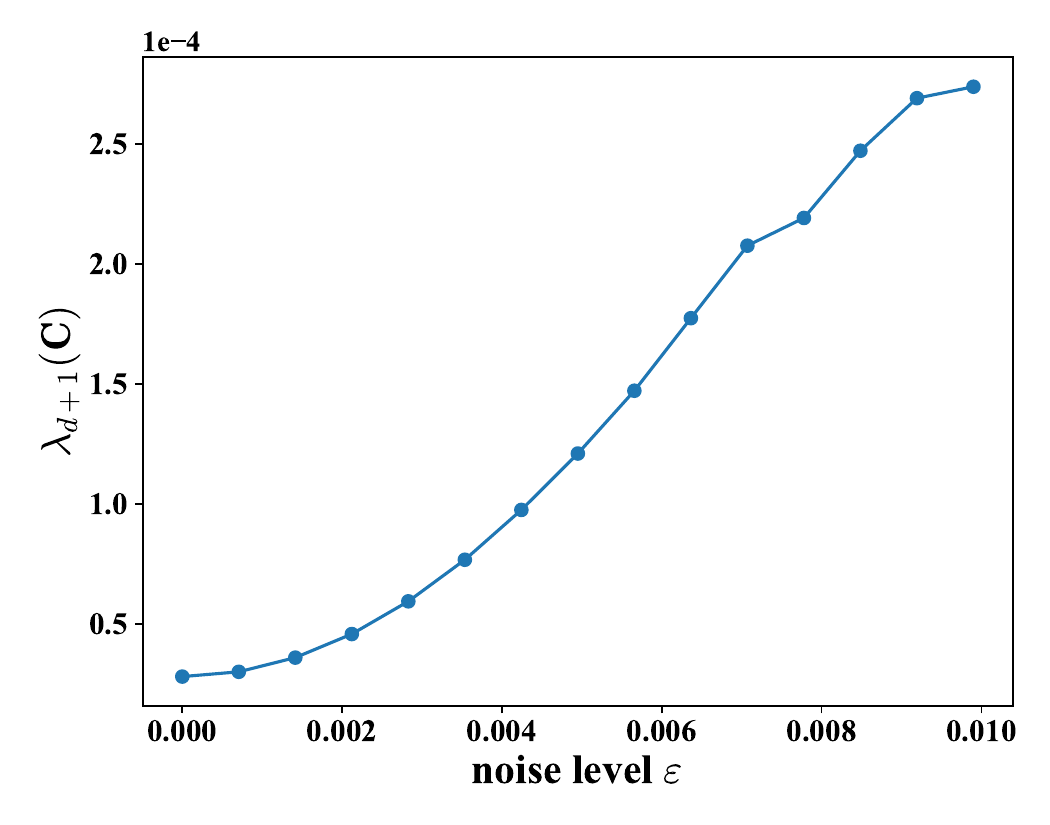}
      &  \includegraphics[width=0.46\linewidth,keepaspectratio]{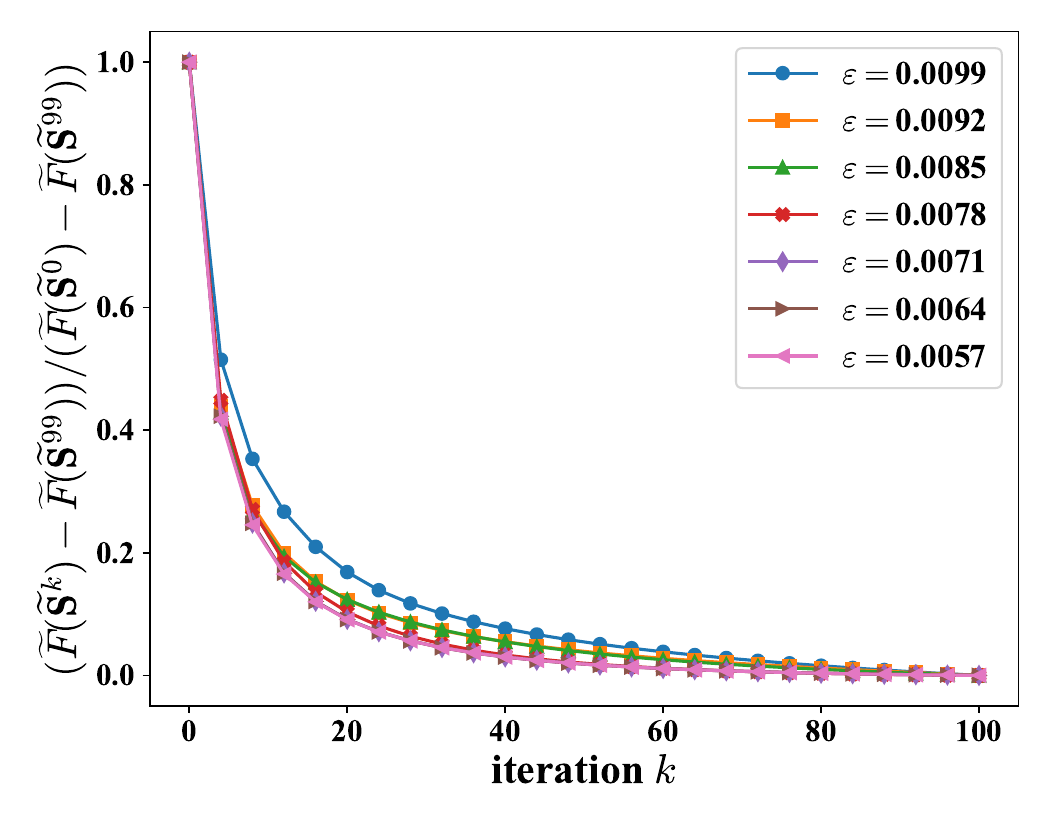}
    \end{tabular}
    \caption{(left) The eigenvalue $\lambda_{d+1}(\mathbf{C})$ against the noise levels. (right) The evolution of the ratio (Eq.~(\ref{eq:alignment_err_ratio})) due to the iterates generated by Algorithm~\ref{algo:rgd} when initialized with the output of SPEC \citea{chaudhury2015global}.}
    \label{fig:rgd_iterates}
\end{figure}

To demonstrate, we provide a simple simulation showing that when Algorithm~\ref{algo:rgd} is initialized with SPEC, it produce iterates with lower alignment error. We took about $n=5000$ points arranged in a unit square grid with a resolution of $70$ points per dimension and subsequently obtained $m=331$ overlapping views. We added random bounded noise in each view for a fixed noise level $\varepsilon$, obtained the corresponding patch-stress matrix $\mathbf{C}$, and computed the spectral alignment $\mathbf{S}_{spec}(\mathbf{C})$ of the noisy views. Finally, we refined it using Algorithm~\ref{algo:rgd} for $100$ iterations. Figure~\ref{fig:rgd_iterates} shows the eigenvalue $\lambda_{d+1}(\mathbf{C})$ against the noise level $\varepsilon$ and confirms the observation in \cite{chaudhury2015global} that the eigenvalue increases with the noise level. Figure~\ref{fig:rgd_iterates} also shows the ratio of
$\widetilde{F}(\widetilde{\mathbf{S}}^k)- \widetilde{F}(\widetilde{\mathbf{S}}^*)$ and $\widetilde{F}(\widetilde{\mathbf{S}}^0) - \widetilde{F}(\widetilde{\mathbf{S}}^*)$ for the alignments $\widetilde{\mathbf{S}}^k$ produced by RGD. The evolution of the ratio is consistent with the linear convergence predicted by Theorem~\ref{thm:rgd_conv2}.

\section{Discussion}
\label{sec:conc}
The following questions remain unanswered and we hope to address them in future work.
\begin{enumerate}[leftmargin=*]
    \item As we showed in Section~\ref{subsec:loc_glob_rigid}, in the case of noiseless views, a non-degenerate perfect alignment  characterizes the infinitesimal rigidity of the resulting realization and also its local rigidity if the realization is generic (Theorem~\ref{thm:inf_rigid} and \ref{thm:loc_rigid}). These results serve as a geometric interpretation of the non-degeneracy conditions in the noiseless case. Although a similar interpretation in the case of noisy views is still missing, it seems natural to expect that the non-degeneracy of an alignment would be associated with infinitesimal rigidity of the consensus representation of the framework (Definition~\ref{def:realization}). We conjecture that the two notions are equivalent. Notably, such a result would offer physically interpretable insights into the algebraic structure of $\mathbb{L}(\mathbf{S})$ (Remark~\ref{rmk:mathbb_L_structure}).
    \item In relation to the previous problem, and analogous to those presented in Section~\ref{subsec:non_deg_noiseless_setting} and Section~\ref{subsec:uniq_noiseless_setting}, necessary and sufficient conditions on the overlapping structure of $m > 2$ \textit{noisy} views for a non-degenerate alignment and for a unique optimal alignment are still unknown.
    \item The proof/counterexample of the converse of Theorem~\ref{thm:nec_cond_glob_rigid_views} is to be investigated.
    \item The relationship between the eigenvalues of $\mathbb{L}(\mathbf{S})$ and those of $\mathbf{C}(\mathbf{S})$ remains unclear in both the noisy and noiseless settings. Specifically, the connection between $\lambda_{d(d-1)/2+2}(\mathbb{L}(\mathbf{S}))$ and the eigenvalues of $\mathbf{C}$ should be explored as it would aid in determining the radius of convergence $\delta(\mathbf{S})$ for RGD (Theorem~\ref{thm:rgd_conv2}) in terms of the eigenvalues of $\mathbf{C}$ instead of $\mathbb{L}(\mathbf{S})$. It is important to note that a direct application of Ostrowski's theorem \cite{higham1998modifying} does not seem to provide the desired link.
    \item Requiring an alignment to be non-degenerate for the local linear convergence of RGD is, in a sense, a strong ask. The convergence to arbitrary critical points using a geometric approach based on \citeb[Section 6.2]{usevich2020approximate} is to be investigated.
    \item It remains unclear whether the solution of the spectral relaxation of Eq.~(\ref{eq:GPOP}) remains close to an optimal alignment under weaker infinitesimal or local rigidity constraints.
\end{enumerate}

\section*{Acknowledgements}
We thank Lijun Ding for providing references which assisted the noise stability analysis. DK was partly supported by a grant from Kavli Institute for Brain and Mind (UCSD). GM was partially funded by NSF CCF-2217058. AC was partially funded by NSF DMS-2012266 and a gift from Intel Research.

\appendix
\section{Notation and Proofs}
\label{supp:sec:all_proofs}
$[a,b]$ is the set $\{a,\ldots,b\}$ where $a,b \in \mathbb{Z}$.
$(a_i)_1^k$ is the sequence $a_1,\ldots,a_k$ where $a_i$ is either a scalar or a vector or a matrix.
$\mathbf{e}^p_q$ is a vector of zeros of length $p$ with $1$ at the $q$th location.
$\mathbf{1}^{p}_q$ is a vector of zeros of length $p$ whose first $q$ elements are $1$s.
$\mathbf{1}_p$ equals $\mathbf{1}^{p}_p$.
$\mathbf{0}_p$ and $\mathbf{0}_{m \times n}$ is a vector of length $p$ and a matrix of zeros with $m$ rows and $n$ columns, respectively.
$[\mathbf{A}_i]_1^n$ denotes a matrix obtained by vertically stacking the matrices $(\mathbf{A}_i)_1^m$.
$[\mathbf{A}]_1^n$ equals $[\mathbf{A}_i]_1^n$ where $\mathbf{A}_i = \mathbf{A}$ for all $i \in [1,n]$.
$\mathbf{I}_d$ and $\mathbf{I}^m_d$ denotes the identity matrix of size $d$ and $[\mathbf{I}_d]_1^m$, respectively.
$\mathbf{A}_i$ is the $i$th row block of $\mathbf{A}$ (the dimensions are contextual).
$\mathbf{A}_{ij}$ is the $(i,j)$th block of the block matrix $\mathbf{A}$ (the dimensions are contextual).
$\vecz (\mathbf{A})$ denotes the column-major vectorization of the matrix $\mathbf{A}$.
$\blockdiag((\mathbf{A}_i)_1^m)$ is a block diagonal matrix with $\mathbf{A}_i$ as the $i$th block.
$\diag (\mathbf{v})$ is a diagonal matrix with $\mathbf{v}(i)$ as the $i$th diagonal element.
$\mathbf{A}(i:j,:) (\mathbf{A}(:,i:j))$ denotes a stacking of $i$th to $j$th rows (columns) of $\mathbf{A}$.
$\mathbb{O}(n)$ is the set of orthogonal matrices of size $n$.
$\Sym (n)$ and $\Skew (n)$ is the set of symmetric and skew-symmetric matrices of size $n$, respectively.
$\Sym(\mathbf{A})$ and $\Skew(\mathbf{A})$ equals $(\mathbf{A}+\mathbf{A}^T)/2$ and $(\mathbf{A}-\mathbf{A}^T)/2$, respectively.
$\mathbf{A} \otimes \mathbf{B}$ denotes the Kronecker product of $\mathbf{A}$ and $\mathbf{B}$. \revadd{For a symmetric matrix $\mathbf{A}$, $\lambda_{i}(\mathbf{A})$ denotes the $i$th smallest eigenvalue of $\mathbf{A}$ and $\sigma_{\min}(\mathbf{A})$ ($\sigma_{\max}(\mathbf{A})$) denotes the smallest (largest) singular value of $\mathbf{A}$.}
\proofof{Proposition~\ref{prop:A_0_2_A_1}}
By differentiating the objective in Eq.~(\ref{eq:A_1}) with respect to $\mathbf{x}_k$, the optimal $\mathbf{x}_k^* \coloneqq \mathbf{x}_k((\mathbf{S}_i)_1^m,(\mathbf{t}_i)_1^m) = n_i^{-1}\textstyle\sum_{(k,i)\in E(\Gamma)}(\mathbf{S}_i^T\mathbf{x}_{k,i}+\mathbf{t}_i)$ (where $n_i = |\{i:(k,i)\in E(\Gamma)\}|$ is the number of points in the $i$th view) which is the consensus of all the views for the $k$th point. Since $\textstyle\sum_{(k,i)\in E(\Gamma)}(\mathbf{S}_i^T\mathbf{x}_{k,i}+\mathbf{t}_i-\mathbf{x}_k^*) = 0$, by adding and subtracting $\mathbf{x}_k^*$, we conclude that 
$$\textstyle\sum_{\substack{(k,i)\in E(\Gamma)\\(k,j)\in E(\Gamma)}}\left\|(\mathbf{S}_i^T\mathbf{x}_{k,i}+\mathbf{t}_i)-(\mathbf{S}_j^T\mathbf{x}_{k,j}+\mathbf{t}_j)\right\|^2_2 = 2\textstyle\sum_{(k,i)\in E(\Gamma)}\left\|(\mathbf{S}_i^T\mathbf{x}_{k,i}+\mathbf{t}_i)-\mathbf{x}_k^*\right\|^2_2$$
and the result follows. The above equality also shows that the minimizers
of Eq.~(\ref{eq:A_0}) and Eq.~(\ref{eq:A_1}) are the same.

\proofof{Proposition~\ref{prop:kerB}}
Let $\mathbf{u} \in \mathbb{R}^{n+m}$. Since $\boldsymbol{\mathcal{L}}_{\Gamma} \succeq 0$, $\boldsymbol{\mathcal{L}}_{\Gamma}\mathbf{u} = 0 \iff \mathbf{u}^T\boldsymbol{\mathcal{L}}_{\Gamma}\mathbf{u} = 0$. From Eq.~(\ref{eq:L_Gamma}), the latter holds iff $\mathbf{e}_{ki}^T\mathbf{u} = 0$ for all $(k,i) \in E(\Gamma)$. Thus, using Eq.~(\ref{eq:B}), $\boldsymbol{\mathcal{L}}_{\Gamma}\mathbf{u} = 0$ implies $\mathbf{B}\mathbf{u} = 0$. 

\proofoffirst{Proposition~\ref{prop:V_S_H_S}}
Using Eq.~(\ref{eq:pi_inv_wtS}),
\begin{align}
    \mathcal{V}_{\mathbf{S}} &= T_{\mathbf{S}}\pi^{-1}(\widetilde{\mathbf{S}}) = \{\mathbf{Z} \in \mathbb{R}^{md \times d}: \mathbf{Z}_i\mathbf{S}_1^T+\mathbf{S}_i\mathbf{Z}_1^T = 0, \mathbf{Z}_j\mathbf{S}_j^T+\mathbf{S}_j\mathbf{Z}_j^T=0, i \in [2,m], j \in [1,m]\}\\
    &= \{[\mathbf{S}_i\boldsymbol{\Omega}_i]_1^m: \boldsymbol{\Omega}_i \in \Skew(d), \boldsymbol{\Omega}_i+\boldsymbol{\Omega}_1^T = 0, i \in [1,m]\}\\
    &= \{\mathbf{S}\boldsymbol{\Omega}: \boldsymbol{\Omega} \in \Skew(d)\}
\end{align}

Since the objective in Eq.~(\ref{eq:P^v_S}) is strictly convex in $\boldsymbol{\Omega}$, it suffices to solve for the critical point and the Eq.~(\ref{eq:P^v_S}) follows immediately. Then, using Eq.~(\ref{eq:g_Z_W}), we obtain
\begin{align}
    \mathcal{H}_{\mathbf{S}} &= \mathcal{V}_{\mathbf{S}}^\perp = \{\mathbf{W} \in T_{\mathbf{S}}\mathbb{O}(d)^m: \Tr(\mathbf{Z}^T\mathbf{W}) = 0 \text{ for all } \mathbf{Z} \in \mathcal{V}_{\mathbf{S}}\}\\
    &= \left\{[\mathbf{S}_i\boldsymbol{\Omega}_i]_1^m: \boldsymbol{\Omega}_i \in \Skew(d) \text{ and }\Tr\left(\boldsymbol{\Omega}^T\left(\textstyle\sum_1^m\boldsymbol{\Omega}_i\right)\right) = 0 \text{ for all } \boldsymbol{\Omega}\in \Skew(d)\right\}.
\end{align}
The constraints on $\boldsymbol{\Omega}_i$ are equivalent to $\boldsymbol{\Omega}_i \in \Skew(d)$ and $ \textstyle\sum_1^m \boldsymbol{\Omega}_i \in \Sym(d)$, and subsequently to $\boldsymbol{\Omega}_i \in \Skew(d)$ and $\textstyle\sum_1^m \boldsymbol{\Omega}_i= 0$. Since $\mathcal{H}_{\mathbf{S}}$ is the orthogonal complement to $\mathcal{V}_{\mathbf{S}}$ in $T_{\mathbf{S}}\mathbb{O}(d)^m$, Eq.~(\ref{eq:P^h_S}) follows trivially.

\proofof{Proposition~\ref{prop:hlift_char}}
By Eq.~(\ref{eq:hlift_def}), we obtain, $\lim_{t \rightarrow 0}(\pi(\mathbf{S}+t\mathbf{Z})_i-\pi(\mathbf{S})_i)/t = \widetilde{\mathbf{Z}}_i$ for each $i \in [1,m-1]$ which further implies $\mathbf{S}_{i+1}\mathbf{Z}_1^T + \mathbf{Z}_{i+1}\mathbf{S}_1^T = \widetilde{\mathbf{Z}}_i$.
Since $T_{\widetilde{\mathbf{S}}}\mathbb{O}(d)^m/_{\sim}$ is identified with $T_{\widetilde{\mathbf{S}}}\mathbb{O}(d)^{m-1}$, therefore there exist $(\widetilde{\boldsymbol{\Omega}}_i)_1^{m-1} \subseteq \Skew(d)$, such that $\widetilde{\mathbf{Z}}_i = \widetilde{\mathbf{S}}_i\widetilde{\boldsymbol{\Omega}}_i$. Also, since $\mathbf{Z} \in \mathcal{H}_{\mathbf{S}}$, there exist $(\boldsymbol{\Omega}_i)_1^m \subseteq \Skew(d)$ such that $\textstyle\sum_1^m \boldsymbol{\Omega}_i = 0$ and $\mathbf{Z}_i=\mathbf{S}_i\boldsymbol{\Omega}_i$. Substituting $\widetilde{\mathbf{Z}}_i = \widetilde{\mathbf{S}}_i\widetilde{\boldsymbol{\Omega}}_i$ and $\mathbf{Z}_i=\mathbf{S}_i\boldsymbol{\Omega}_i$,
we obtain
\begin{equation}
    \mathbf{S}_{i+1}\boldsymbol{\Omega}_1^T\mathbf{S}_1^T + \mathbf{S}_{i+1}\boldsymbol{\Omega}_{i+1}\mathbf{S}_1^T = \widetilde{\mathbf{S}}_i\widetilde{\boldsymbol{\Omega}}_i \implies \boldsymbol{\Omega}_{i+1}-\boldsymbol{\Omega}_1 = \mathbf{S}_{1}^T\widetilde{\boldsymbol{\Omega}}_i\mathbf{S}_{1}, i \in [1,m-1], \label{supp:eq:eq1_}
\end{equation}
where we used the fact that $\widetilde{\mathbf{S}}_i = \mathbf{S}_{i+1}\mathbf{S}_1^T$. Observe that the linear system in $[\boldsymbol{\Omega}_i]_1^m$ is of full rank. By applying $\textstyle\sum_{i=1}^{m-1}$ and using $\textstyle\sum_1^m \mathbf{\Omega}_i = 0$ gives Eq.~(\ref{eq:hlift1}) and (\ref{eq:hlifti}).

\proofof{Proposition~\ref{prop:g_tilde}}
It suffices to show that $g(\mathbf{Z}, \mathbf{W})$ does not depend on the choice of $\mathbf{S} \in \pi^{-1}(\widetilde{\mathbf{S}})$. Let $\{\widetilde{\mathbf{U}}_i,\widetilde{\mathbf{V}}_i\}_1^{m-1},\{\mathbf{U}_i,\mathbf{V}_i\}_1^m$ be elements of $\Skew(d)$ such that $\widetilde{\mathbf{Z}}_i = \widetilde{\mathbf{S}}_i\widetilde{\mathbf{U}}_i$, $\widetilde{\mathbf{W}}_i = \widetilde{\mathbf{S}}_i\widetilde{\mathbf{V}}_i$, $\mathbf{Z}_i = \mathbf{S}_i\mathbf{U}_i$ and $\mathbf{W}_i = \mathbf{S}_i\mathbf{V}_i$. By the definition of $\mathcal{H}_{\mathbf{S}}$, $\textstyle\sum_1^m\mathbf{U}_i = \textstyle\sum_1^m \mathbf{V}_i = 0$ and the relation between $\mathbf{U}_i$ and $\widetilde{\mathbf{U}}_i$, and $\mathbf{V}_i$ and $\widetilde{\mathbf{V}}_i$ is given by Eq.~(\ref{eq:hlift1}, \ref{eq:hlifti}) (through Eq.~(\ref{supp:eq:eq1_})). Then we have
\begin{equation}
    g(\mathbf{Z}, \mathbf{W}) = \textstyle\sum_1^m \Tr(\mathbf{Z}_i^T\mathbf{W}_i) = \textstyle\sum_1^m \Tr(\mathbf{U}_i^T\mathbf{V}_i) = \textstyle\sum_1^m \Tr((\mathbf{U}_i-\mathbf{U}_1)^T\mathbf{V}_i) + \Tr(\mathbf{U}_1^T\mathbf{V}_i)
\end{equation}
Since $\textstyle\sum_1^m \mathbf{V}_i = 0$, the second term vanishes. The first reduces to 
\begin{align}
    &\textstyle\sum_1^m \Tr((\mathbf{U}_i-\mathbf{U}_1)^T\mathbf{V}_i) = \textstyle\sum_1^m \Tr((\mathbf{U}_i-\mathbf{U}_1)^T(\mathbf{V}_i-\mathbf{V}_1)) + \Tr((\mathbf{U}_i-\mathbf{U}_1)^T\mathbf{V}_1)\\
    &= \textstyle\sum_1^{m-1} \Tr((\mathbf{U}_{i+1}-\mathbf{U}_1)^T(\mathbf{V}_{i+1}-\mathbf{V}_1)) - \revadd{m}\Tr(\mathbf{U}_1^T\mathbf{V}_1)\\
    &=\textstyle\sum_1^{m-1}\Tr(\widetilde{\mathbf{U}}_{i}^T\widetilde{\mathbf{V}}_{i}) - m\Tr(\mathbf{U}_1^T\mathbf{V}_1)\\
    &= \textstyle\sum_1^{m-1}\Tr(\widetilde{\mathbf{U}}_{i}^T\widetilde{\mathbf{V}}_{i}) - \revdel{m^{-2}}\revadd{m^{-1}}\textstyle\sum_{i,j=1}^{m-1}\Tr(\widetilde{\mathbf{U}}_i^T\widetilde{\mathbf{V}}_j).
\end{align}
Since the above equation is independent of the choice of $\mathbf{S}$, the result follows.

\proofof{Proposition~\ref{prop:hlift_frob_ineq}} 
We have $\left\|\mathbf{Z}\right\|_F^2 = \sum_1^m \Tr(\mathbf{\Omega}_i\mathbf{\Omega}_i^T)$. Using Proposition~\ref{prop:hlift_char},
\begin{align}
    \left\|\widetilde{\mathbf{Z}}\right\|_F^2 &= \sum_1^{m-1} \Tr(\widetilde{\mathbf{\Omega}}_{i+1}\widetilde{\mathbf{\Omega}}_{i+1}^T) =\sum_1^{m-1} (\Tr(\mathbf{\Omega}_{i+1}\mathbf{\Omega}_{i+1}^T) -2\Tr(\mathbf{\Omega}_{i+1}\mathbf{\Omega}_1^T)) + (m-1)\Tr(\mathbf{\Omega}_1\mathbf{\Omega}_1^T)\\
    &=\sum_1^m \Tr(\mathbf{\Omega}_i\mathbf{\Omega}_i^T) + m \sum_1^m \Tr(\mathbf{\Omega}_1\mathbf{\Omega}_1^T) = \left\|\mathbf{Z}\right\|_F^2 + m\left\|\mathbf{\Omega}_1\right\|_F^2.
\end{align}

\proofof{Proposition~\ref{prop:gradFS}}
Using \citeb[Section 3.6.2]{absil2009optimization} and Proposition~\ref{prop:T_SOdm},
\begin{align}
    \overline{\grad \widetilde{F}(\widetilde{\mathbf{S}})} &= \grad F(\mathbf{S}) = \mathbf{P}_{\mathbf{S}}(\nabla F(\mathbf{S})) = \mathbf{P}_{\mathbf{S}}(2\mathbf{C}\mathbf{S}) =[2\mathbf{S}_i\text{skew}(\mathbf{S}_i^T[\mathbf{C}\mathbf{S}]_i])]_1^m\\
    &=[\mathbf{C}\mathbf{S}]_i - \mathbf{S}_i[\mathbf{C}\mathbf{S}]_i^T\mathbf{S}_i = \mathbf{S}_i (\mathbf{S}_i^T[\mathbf{C}\mathbf{S}]_i - [\mathbf{C}\mathbf{S}]_i^T\mathbf{S}_i)
\end{align}
Using the fact that $\mathbf{C}$ is symmetric, we conclude that $\sum_{1}^{m} \boldsymbol{\Omega}_i = \mathbf{S}^T\mathbf{C}\mathbf{S}-\mathbf{S}^T\mathbf{C}^T\mathbf{S} = 0$.

\proofof{Proposition~\ref{prop:DgradFSZ}}
\begin{align}
    D\grad F(\mathbf{S})[\mathbf{Z}]_i &= \lim_{t \rightarrow 0}(\grad F(\mathbf{S}+t\mathbf{Z})_i-\grad F(\mathbf{S})_i)/t\\
    &= \lim_{t \rightarrow 0} t^{-1}\left\{([\mathbf{C}(\mathbf{S}+t\mathbf{Z})]_i - (\mathbf{S}_i+t\mathbf{Z}_i)[\mathbf{C}(\mathbf{S}+t\mathbf{Z})]_i^T(\mathbf{S}_i+t\mathbf{Z}_i)) - ([\mathbf{C}\mathbf{S}]_i - \mathbf{S}_i[\mathbf{C}\mathbf{S}]_i^T\mathbf{S}_i)\right\}\\
    &=[\mathbf{C}\mathbf{Z}]_i - \mathbf{S}_i[\mathbf{C}\mathbf{Z}]_i^T\mathbf{S}_i - \mathbf{S}_i[\mathbf{C}\mathbf{S}]_i^T\mathbf{Z}_i -\mathbf{Z}_i[\mathbf{C}\mathbf{S}]_i^T\mathbf{S}_i\\
    &=\mathbf{S}_i(\mathbf{S}_i^T[\mathbf{C}\mathbf{Z}]_i - [\mathbf{C}\mathbf{Z}]_i^T\mathbf{S}_i - [\mathbf{C}\mathbf{S}]_i^T\mathbf{Z}_i - \mathbf{S}_i^T\mathbf{Z}_i[\mathbf{C}\mathbf{S}]_i^T\mathbf{S}_i)
\end{align}
\revadd{The result follows from the facts that $\mathbf{S}_i^T\mathbf{Z}_i + \mathbf{Z}_i^T\mathbf{S}_i = 0$ for $\mathbf{Z} \in T_{\mathbf{S}}\mathbb{O}(d)^m$ (see Proposition~\ref{prop:T_SOdm}), and $[\mathbf{C}\mathbf{S}]_i = \mathbf{S}_i[\mathbf{C}\mathbf{S}]_i^T\mathbf{S}_i$ for $\widetilde{\mathbf{S}} \in \widetilde{\mathcal{C}}$ (see the proof of Proposition~\ref{prop:gradFS})}.

\proofof{Proposition~\ref{prop:HessFSZ}}
Using \citeb[Chapter 5]{absil2009optimization},
$$\overline{\Hess \widetilde{F}(\widetilde{\mathbf{S}})[\widetilde{\mathbf{Z}}]} = \overline{\widetilde{\nabla}_{\widetilde{\mathbf{Z}}}\grad \widetilde{F}(\widetilde{\mathbf{S}})} = P^{h}_{\mathbf{S}}(\nabla_{\overline{\widetilde{\mathbf{Z}}}}\overline{\grad \widetilde{F}(\widetilde{\mathbf{S}})}) = P^{h}_{\mathbf{S}}(\nabla_{\mathbf{Z}}\grad F(\mathbf{S})).$$ Then from, Proposition~\ref{prop:T_SOdm}, \ref{prop:V_S_H_S} and \ref{prop:DgradFSZ}, the latter is reduced to
\begin{align}
    P^{h}_{\mathbf{S}}(\nabla_{\mathbf{Z}}\grad F(\mathbf{S})) &= P^{h}_{\mathbf{S}}(P_{\mathbf{S}}(D\grad F(\mathbf{S})[\mathbf{Z}])) = P^{h}_{\mathbf{S}}\left(\left[\mathbf{S}_i(\Skew(\boldsymbol{\xi}_i))\right]_1^m\right)\\
    &= [\mathbf{S}_i(\Skew(\boldsymbol{\xi}_i)-m^{-1}\sum_1^m\Skew(\boldsymbol{\xi}_i))]_1^m
\end{align}
For the case $\widetilde{\mathbf{S}} \in \widetilde{\mathcal{C}}$, we note that $\boldsymbol{\xi}_i = \Skew(\boldsymbol{\xi}_i) = \mathbf{S}_i^T[\mathbf{C}\mathbf{Z}]_i - [\mathbf{C}\mathbf{Z}]_i^T\mathbf{S}_i - [\mathbf{C}\mathbf{S}]_i^T\mathbf{Z}_i + \mathbf{Z}_i^T[\mathbf{C}\mathbf{S}]_i$ and $\sum_{1}^{m}\boldsymbol{\xi}_i = \mathbf{S}^T\mathbf{C}\mathbf{Z}-\mathbf{Z}^T\mathbf{C}\mathbf{S}-\mathbf{S}^T\mathbf{C}\mathbf{Z}+\mathbf{Z}^T\mathbf{C}\mathbf{S} = 0$.

\proofof{Proposition~\ref{prop:Omega_hat_compact}}
Since $\mathbf{Z} \in \mathcal{H}_{\mathbf{S}}$, using Proposition~\ref{prop:V_S_H_S}, there exist $\boldsymbol{\Omega} = [\boldsymbol{\Omega}_i]_1^m$ such that $\boldsymbol{\Omega}_i \in \Skew(d)$, $\textstyle\sum_1^m\boldsymbol{\Omega}_i = 0$ and $\mathbf{Z}_i = \mathbf{S}_i\boldsymbol{\Omega}_i$. Then,
\revadd{$\boldsymbol{\xi}_i = [\mathbf{C}(\mathbf{S})\boldsymbol{\Omega}]_i - [\mathbf{C}(\mathbf{S})\boldsymbol{\Omega}]_i^T - [\widehat{\mathbf{C}}(\mathbf{S})^T\boldsymbol{\Omega}]_i + [\widehat{\mathbf{C}}(\mathbf{S})\boldsymbol{\Omega}]_i^T$. The result follows from (i) $\mathbf{L}(\mathbf{S}) = \mathbf{C}(\mathbf{S}) - \widehat{\mathbf{C}}(\mathbf{S})$, (ii) $\mathbf{C}(\mathbf{S})$ is symmetric and (iii) for $\widetilde{\mathbf{S}}\in \widetilde{\mathcal{C}}$, $\widehat{\mathbf{C}}(\mathbf{S})$ is also symmetric (Remark~\ref{rmk:StildeCtilde}).}

\proofof{Proposition~\ref{prop:HessFSZZ}}
We obtain $\widetilde{g}(\Hess \widetilde{F}(\widetilde{\mathbf{S}})[\widetilde{\mathbf{Z}}],\widetilde{\mathbf{Z}}) = g(\overline{\Hess \widetilde{F}(\widetilde{\mathbf{S}})[\widetilde{\mathbf{Z}}]},\overline{\widetilde{\mathbf{Z}}})$ from Proposition~\ref{prop:g_tilde}. Due to Proposition \ref{prop:HessFSZ},
$$g(\overline{\Hess \widetilde{F}(\widetilde{\mathbf{S}})[\widetilde{\mathbf{Z}}]},\overline{\widetilde{\mathbf{Z}}}) = g(\overline{\Hess \widetilde{F}(\widetilde{\mathbf{S}})[\widetilde{\mathbf{Z}}]},\mathbf{Z}) = \textstyle\sum_{1}^{m}\Tr((\mathbf{S}_i\widehat{\boldsymbol{\Omega}}_i)^T\mathbf{S}_i\boldsymbol{\Omega}_i) = \textstyle\sum_{1}^{m}\Tr(\widehat{\boldsymbol{\Omega}}_i^T\boldsymbol{\Omega}_i).$$ \revadd{The result follows from the Proposition~\ref{prop:Omega_hat_compact} and the fact that $\Tr(\mathbf{A}) = \Tr(\mathbf{A}^T)$.}

\proofof{Proposition~\ref{prop:Omega^TLSOmega2}} Using $\vecz (\mathbf{A}\mathbf{X}\mathbf{B}) = (\mathbf{B}^T \otimes \mathbf{A})\vecz (\mathbf{X})$, we obtain
\begin{align}
\Tr(\boldsymbol{\Omega}^T&\mathbf{L}(\mathbf{S})\boldsymbol{\Omega}) = \Tr((\mathbf{P}\boldsymbol{\Omega})^T\mathbf{P}\mathbf{L}(\mathbf{S})\mathbf{P}^T(\mathbf{P}\boldsymbol{\Omega})) = \vecz (\mathbf{P}\boldsymbol{\Omega})^T\vecz (\mathcal{L}(\mathbf{S}) (\mathbf{P}\boldsymbol{\Omega}))\\
&= \vecz (\mathbf{P}\boldsymbol{\Omega})^T(\mathbf{I}_d \otimes \mathcal{L}(\mathbf{S}))\vecz (\mathbf{P}\boldsymbol{\Omega}) = \boldsymbol{\omega}^T\overline{\mathbf{P}}(\mathbf{I}_d \otimes \mathcal{L}(\mathbf{S}))\overline{\mathbf{P}}^T\boldsymbol{\omega} = \boldsymbol{\omega}^T\mathbb{L}(\mathbf{S})\boldsymbol{\omega}.
\end{align}

\proofof{Theorem~\ref{thm:non_deg_loc_min}}
$1$, $2$ and $3$ are equivalent by definitions and Proposition~\ref{prop:HessFSZZ}. For ($3 \iff 4$), by comparing dimensions, we note that the set of $\boldsymbol{\Omega} = [\boldsymbol{\Omega}_i]_1^m$ where $\boldsymbol{\Omega}_i \in \Skew(d)$ and $\textstyle\sum_1^m \boldsymbol{\Omega}_i = 0$, is the same as the set of $\boldsymbol{\Omega} = [\boldsymbol{\Omega}_i-\boldsymbol{\Omega}_0]_1^m$ where $\boldsymbol{\Omega}_i \in \Skew(d)$ and $\boldsymbol{\Omega}_0 = \frac{1}{m}\textstyle\sum_1^m\boldsymbol{\Omega}_i$. Using Remark~\ref{rmk:C_hat_L_structure}, we know that $\mathbf{L}(\mathbf{S})=\mathbf{L}(\mathbf{S})^T$ and $\mathbf{L}(\mathbf{S})[\boldsymbol{\Omega}_0]_1^m = 0$. Thus $(\boldsymbol{\Omega}-[\boldsymbol{\Omega}_0]_1^m)^T\mathbf{L}(\mathbf{S})(\boldsymbol{\Omega}-[\boldsymbol{\Omega}_0]_1^m) = \boldsymbol{\Omega}^T\mathbf{L}(\mathbf{S})\boldsymbol{\Omega}$ and the result follows. Subsequently, ($4 \iff 5$) follows directly from the definition of $\boldsymbol{\omega}$ and Proposition~\ref{prop:Omega^TLSOmega2}. Finally, ($5 \iff 6$) and ($6 \iff 7$) follow from Remark~\ref{rmk:mathbb_L_structure}.

\proofof{Proposition~\ref{prop:one_all1}} Suppose $\mathbf{L}(\mathbf{S})$ satisfies condition 4 and $\boldsymbol{\Omega}$ be as in the condition. Then using Remark~\ref{rmk:C_hat_L_structure}, 
\begin{align}
    \Tr(\boldsymbol{\Omega}^T\mathbf{L}(\mathbf{S}\mathbf{\mathbf{Q}})\boldsymbol{\Omega}) &= \Tr(\boldsymbol{\Omega}^T(\mathbf{I}_m \otimes \mathbf{Q})^T\mathbf{L}(\mathbf{S})(\mathbf{I}_m \otimes \mathbf{Q})\boldsymbol{\Omega}) = \Tr(\overline{\boldsymbol{\Omega}}^T \mathbf{L}(\mathbf{S})\overline{\boldsymbol{\Omega}})
\end{align}
which is positive because $\overline{\boldsymbol{\Omega}}_i = \mathbf{Q}\boldsymbol{\Omega}_i\mathbf{Q}^T \in \Skew(d)$ for all $i \in [1,m]$ and not all $\overline{\boldsymbol{\Omega}}_i$ are equal (if $\overline{\boldsymbol{\Omega}}_i = \overline{\boldsymbol{\Omega}}_j$ then $\boldsymbol{\Omega}_i = \boldsymbol{\Omega}_j$, a contradiction). Thus, condition 4 holds for~$\mathbf{S}\mathbf{Q}$ too. The result follows.

\proofof{Corollary~\ref{cor:suff_non_deg_loc_min}}
Since the rank of $\mathbf{L}(\mathbf{S})$ is $(m-1)d$, using Remark~\ref{rmk:C_hat_L_structure}, $\mathbf{L}(\mathbf{S})[\boldsymbol{\Omega}_i]_1^m = 0$ (equivalently, $\Tr(([\boldsymbol{\Omega}_i]_1^m)^T\mathbf{L}(\mathbf{S})[\boldsymbol{\Omega}_i]_1^m) = 0$) iff $\boldsymbol{\Omega}_i = \boldsymbol{\Omega}_0$ for all $i \in [1,m]$ and some $\boldsymbol{\Omega}_0 \in \Skew(d)$. Thus, condition~$4$ in Theorem~\ref{thm:non_deg_loc_min} is satisfied.

\proofof{Proposition~\ref{prop:HessVicinity}} For brevity, we define $\mathbf{D}_{\mathbf{S}} \coloneqq \blockdiag([\mathbf{S}_i]_1^m)$ and $\mathbf{Z}_i \coloneqq \mathbf{O}_i\boldsymbol{\Omega}_i$, i $\in [1,m]$, which satisfies $\left\|\mathbf{Z}_i\right\|_F^2 = \left\|\boldsymbol{\Omega}_i\right\|_F^2$. Then,
\begin{align}
    \Tr(\boldsymbol{\Omega}^T(\mathbf{L}(\mathbf{O})+\mathbf{L}(\mathbf{O})^T)\boldsymbol{\Omega}) &= 2\Tr(\boldsymbol{\Omega}^T\mathbf{C}(\mathbf{O})\boldsymbol{\Omega}) - \Tr(\boldsymbol{\Omega}^T(\widehat{\mathbf{C}}(\mathbf{O})+\widehat{\mathbf{C}}(\mathbf{O})^T)\boldsymbol{\Omega})\\
    &= 2\Tr(\boldsymbol{\Omega}^T\mathbf{C}(\mathbf{O})\boldsymbol{\Omega}) - \textstyle\sum_{i=1}^{m}\Tr\left(\boldsymbol{\Omega}_i^T \left(\textstyle\sum_{j=1}^{m}\mathbf{O}_i^T\mathbf{C}_{ij}\mathbf{O}_j + \mathbf{O}_j^T\mathbf{C}_{ji}\mathbf{O}_i\right)\boldsymbol{\Omega}_i\right)\\
    &= 2\Tr(\boldsymbol{\Omega}^T\mathbf{C}(\mathbf{O})\boldsymbol{\Omega}) - \textstyle\sum_{i=1}^{m}\Tr\left(\mathbf{Z}_i^T \left(\textstyle\sum_{j=1}^{m}\mathbf{C}_{ij}\mathbf{O}_j\mathbf{O}_i^T + \mathbf{O}_i\mathbf{O}_j^T\mathbf{C}_{ji}\right)\mathbf{Z}_i\right)\label{eq:HessO1}
\end{align}
Rewriting the first term,
\begin{align}
\textstyle\sum_{i=1}^{m}\Tr&\left(\mathbf{Z}_i^T\textstyle\sum_{j=1}^{m}\mathbf{C}_{ij}\mathbf{O}_j\mathbf{O}_i^T\mathbf{Z}_i\right) = \textstyle\sum_{i=1}^{m}\Tr\left(\mathbf{Z}_i^T\textstyle\sum_{j=1}^{m}\mathbf{C}_{ij}(\mathbf{O}_j\mathbf{O}_i^T-\mathbf{S}_j\mathbf{S}_i^T + \mathbf{S}_j\mathbf{S}_i^T)\mathbf{Z}_i\right)\\
&= \textstyle\sum_{i=1}^{m}\Tr\left(\mathbf{Z}_i^T\textstyle\sum_{j=1}^{m}\mathbf{C}_{ij}(\mathbf{O}_j\mathbf{O}_i^T-\mathbf{S}_j\mathbf{S}_i^T)\mathbf{Z}_i\right) + \Tr\left(\mathbf{Z}_i^T\mathbf{S}_i\mathbf{S}_i^T\textstyle\sum_{j=1}^{m}\mathbf{C}_{ij} \mathbf{S}_j\mathbf{S}_i^T\mathbf{Z}_i\right)\\
&= \textstyle\sum_{i=1}^{m} \Tr\left(\mathbf{Z}_i^T\textstyle\sum_{j=1}^{m}\mathbf{C}_{ij}(\mathbf{O}_j\mathbf{O}_i^T-\mathbf{S}_j\mathbf{S}_i^T)\mathbf{Z}_i\right) + \Tr\left(\mathbf{Z}^T\mathbf{D}_{\mathbf{S}}\widehat{\mathbf{C}}(\mathbf{S})\mathbf{D}_{\mathbf{S}}^T\mathbf{Z}\right) \label{eq:HessO2}
\end{align}
Using Cauchy-Schwarz inequality and the fact that $\left\|\mathbf{A}_1\mathbf{A}_2\right\|_F \leq \left\|\mathbf{A}_1\right\|_2\left\|\mathbf{A}_2\right\|_F$,
\begin{align}
\left|\textstyle\sum_{i=1}^{m} \Tr\right.&\left.\left(\mathbf{Z}_i^T\textstyle\sum_{j=1}^{m}\mathbf{C}_{ij}(\mathbf{O}_j\mathbf{O}_i^T-\mathbf{S}_j\mathbf{S}_i^T)\mathbf{Z}_i\right)\right| \leq \textstyle\sum_{i=1}^{m} \left\|\textstyle\sum_{j=1}^{m}\mathbf{C}_{ij}(\mathbf{O}_j\mathbf{O}_i^T-\mathbf{S}_j\mathbf{S}_i^T)\right\|_F\left\|\mathbf{Z}_i\right\|_F^2\\
&\leq \textstyle\sum_{i=1}^{m} \left\|\textstyle\sum_{j=1}^{m}\mathbf{C}_{ij}(\mathbf{O}_j\mathbf{O}_i^T-\mathbf{S}_j\mathbf{O}_i^T + \mathbf{S}_j\mathbf{O}_i^T-\mathbf{S}_j\mathbf{S}_i^T)\right\|_F\left\|\mathbf{Z}_i\right\|_F^2\\
&\leq  \textstyle\sum_{i=1}^{m}(\max_{k=1}^{m}\left\|\mathbf{C}_{k,:}\right\|_2 \left\|\mathbf{O}-\mathbf{S}\right\|_F + \left\|[\mathbf{C}\mathbf{S}]_i\right\|_2\left\|\mathbf{O}_i-\mathbf{S}_i\right\|_F)\left\|\boldsymbol{\Omega}_i\right\|_F^2\\
&\leq   \textstyle\sum_{i=1}^{m}(c_1 \left\|\mathbf{O}-\mathbf{S}\right\|_F + c_2(\mathbf{S})\left\|\mathbf{O}_i-\mathbf{S}_i\right\|_F)\left\|\boldsymbol{\Omega}_i\right\|_F^2\\
&\leq  (c_1 + c_2(\mathbf{S}))\left\|\mathbf{O}-\mathbf{S}\right\|_F\left\|\boldsymbol{\Omega}\right\|_F^2 \label{eq:HessO3}
\end{align}
Also, due to Eq.~(\ref{eq:C_of_S}),
$$\Tr(\boldsymbol{\Omega}^T\mathbf{C}(\mathbf{O})\boldsymbol{\Omega}) = \Tr(\boldsymbol{\Omega}^T\mathbf{D}_{\mathbf{O}}^T\mathbf{D}_{\mathbf{S}}\mathbf{C}(\mathbf{S})\mathbf{D}_{\mathbf{S}}^T\mathbf{D}_{\mathbf{O}}\boldsymbol{\Omega}) = \Tr(\mathbf{Z}^T\mathbf{D}_{\mathbf{S}}\mathbf{C}(\mathbf{S})\mathbf{D}_{\mathbf{S}}^T\mathbf{Z}).$$
Combining this with Eq.~(\ref{eq:HessO1}, \ref{eq:HessO2}, \ref{eq:HessO3}), we obtain
\begin{align}
    \Tr(\boldsymbol{\Omega}^T(\mathbf{L}(\mathbf{O})+\mathbf{L}(\mathbf{O})^T)\boldsymbol{\Omega}) &\geq 2\Tr(\mathbf{Z}^T\mathbf{D}_{\mathbf{S}}\mathbf{L}(\mathbf{S})\mathbf{D}_{\mathbf{S}}^T\mathbf{Z}) - 2(c_1 + c_2(\mathbf{S}))\left\|\mathbf{O}-\mathbf{S}\right\|_F\left\|\boldsymbol{\Omega}\right\|_F^2\\
    \Tr(\boldsymbol{\Omega}^T(\mathbf{L}(\mathbf{O})+\mathbf{L}(\mathbf{O})^T)\boldsymbol{\Omega}) &\leq 2\Tr(\mathbf{Z}^T\mathbf{D}_{\mathbf{S}}\mathbf{L}(\mathbf{S})\mathbf{D}_{\mathbf{S}}^T\mathbf{Z}) + 2(c_1 + c_2(\mathbf{S}))\left\|\mathbf{O}-\mathbf{S}\right\|_F\left\|\boldsymbol{\Omega}\right\|_F^2.
\end{align}
Moreover,
\begin{align}
\Tr(\mathbf{Z}^T\mathbf{D}_{\mathbf{S}}\mathbf{L}(\mathbf{S})\mathbf{D}_{\mathbf{S}}^T\mathbf{Z}) &= \Tr(\boldsymbol{\Omega}^T\mathbf{D}_{\mathbf{O}}^T\mathbf{D}_{\mathbf{S}}\mathbf{L}(\mathbf{S})\mathbf{D}_{\mathbf{S}}^T\mathbf{D}_{\mathbf{O}}\boldsymbol{\Omega})\\
&= \Tr(\boldsymbol{\Omega}^T\mathbf{L}(\mathbf{S})\boldsymbol{\Omega}) + 2\Tr(\boldsymbol{\Omega}^T(\mathbf{D}_{\mathbf{O}}^T\mathbf{D}_{\mathbf{S}}-\mathbf{I}_{md})\mathbf{L}(\mathbf{S})\mathbf{D}_{\mathbf{S}}^T\mathbf{D}_{\mathbf{O}}\boldsymbol{\Omega}),
\end{align}
where, for the first term,
\begin{equation}
    (\lambda_{-}(\mathbf{S})/2) \left\|\boldsymbol{\Omega}\right\|_F^2 \leq \Tr(\boldsymbol{\Omega}^T\mathbf{L}(\mathbf{S})\boldsymbol{\Omega}) =  \Tr(\boldsymbol{\omega}^T\mathbf{\mathbb{L}}(\mathbf{S})\boldsymbol{\omega}) \leq  (\lambda_{+}(\mathbf{S})/2) \left\|\boldsymbol{\Omega}\right\|_F^2.
\end{equation}
The fraction $1/2$ appears because $\left\|\boldsymbol{\omega}\right\|_F^2 = \left\|\boldsymbol{\Omega}\right\|_F^2/2$ as in Eq.~(\ref{eq:omega^TmbbLomega}). Then, for the second term, using $|\Tr(\mathbf{A}_1\mathbf{A}_2)| \leq \left\|\mathbf{A}_1\right\|_2\Tr(\mathbf{A}_2)$ and Cauchy-Schwarz inequality,
\begin{align}
\left|\Tr(\boldsymbol{\Omega}^T(\mathbf{D}_{\mathbf{O}}^T\mathbf{D}_{\mathbf{S}}-\mathbf{I}_{md})\mathbf{L}(\mathbf{S})\mathbf{D}_{\mathbf{S}}^T\mathbf{D}_{\mathbf{O}}\boldsymbol{\Omega})\right| &\leq \left\|\mathbf{L}(\mathbf{S})\right\|_2\left|\Tr(\boldsymbol{\Omega}^T(\mathbf{D}_{\mathbf{O}}^T\mathbf{D}_{\mathbf{S}}-\mathbf{I}_{md})\boldsymbol{\Omega})\right|\\
&\leq c_3(\mathbf{S})\left\|\mathbf{S}-\mathbf{O}\right\|_F\left\|\boldsymbol{\Omega}\right\|_F^2.
\end{align}
Overall,
\begin{align}
    \Tr(\boldsymbol{\Omega}^T(\mathbf{L}(\mathbf{O})+\mathbf{L}(\mathbf{O})^T)\boldsymbol{\Omega}) &\geq (\lambda_{-}(\mathbf{S}) - 2(c_1 + c_2(\mathbf{S}) + 2 c_3(\mathbf{S}))\left\|\mathbf{S}-\mathbf{O}\right\|_F)\left\|\boldsymbol{\Omega}\right\|_F^2\\
    \Tr(\boldsymbol{\Omega}^T(\mathbf{L}(\mathbf{O})+\mathbf{L}(\mathbf{O})^T)\boldsymbol{\Omega}) &\leq (\lambda_{+}(\mathbf{S}) + 2(c_1 + c_2(\mathbf{S}) + 2 c_3(\mathbf{S}))\left\|\mathbf{S}-\mathbf{O}\right\|_F)\left\|\boldsymbol{\Omega}\right\|_F^2
\end{align}
\revadd{Consequently, if $\left\|\mathbf{S}-\mathbf{O}\right\|_F < \zeta\delta(\mathbf{S})$ (as defined in the theorem statement) then}
\begin{equation}
\revadd{(1-\zeta)\lambda_{-}(\mathbf{S})\left\|\boldsymbol{\Omega}\right\|_F^2 \leq \Tr(\boldsymbol{\Omega}^T(\mathbf{L}(\mathbf{O})+\mathbf{L}(\mathbf{O})^T)\boldsymbol{\Omega}) \leq (\lambda_{+}(\mathbf{S}) + \zeta \lambda_{-}(\mathbf{S}))\left\|\boldsymbol{\Omega}\right\|_F^2.}
\end{equation}
\revadd{Finally, due to Proposition~\ref{prop:one_all1} and the fact that $\delta, \lambda_{-}$ and $\lambda_{+}$ are invariant under the action of $\mathbb{O}(d)$, we can replace $\mathbf{S}$ by $\mathbf{S}\mathbf{Q}$ for any $\mathbf{Q} \in \mathbb{O}(d)$, and the result follows.}

\proofof{Theorem~\ref{thm:non_deg_two_views_gen_setting}} The proof is divided into three parts specialized to the three conditions in the statement.

\noindent \underline{\textbf{Part 1}}. First note that $\mathbf{S} \in \mathcal{C}$ iff $\mathbf{S}_i^T[\mathbf{C}\mathbf{S}]_i = [\mathbf{C}\mathbf{S}]_i^T\mathbf{S}_i$ for $i = 1,2$ (see Eq.~(\ref{eq:crit_pts2})). Since
\begin{equation}
   \mathbf{S}_1^T[\mathbf{C}\mathbf{S}]_1 - [\mathbf{C}\mathbf{S}]_1^T\mathbf{S}_1 =  [\mathbf{C}\mathbf{S}]_2^T\mathbf{S}_2 - \mathbf{S}_2^T[\mathbf{C}\mathbf{S}]_2 = \mathbf{B}(\mathbf{S})_1\boldsymbol{\mathcal{L}}_{\Gamma}^\dagger \mathbf{B}(\mathbf{S})_2^T- \mathbf{B}(\mathbf{S})_2\boldsymbol{\mathcal{L}}_{\Gamma}^\dagger \mathbf{B}(\mathbf{S})_1^T,
\end{equation}
thus $\mathbf{S} \in \mathcal{C}$ iff $\mathbf{B}(\mathbf{S})_1\boldsymbol{\mathcal{L}}_{\Gamma}^\dagger \mathbf{B}(\mathbf{S})_2^T$ is symmetric. At this point, we note that
\begin{equation}
    \begin{matrix}
        \mathbf{B}(\mathbf{S})_1 & = & [ & \mathbf{X}_1 & \mathbf{0}_{d \times n_2} & \mathbf{X}_3 & -(\mathbf{X}_1\mathbf{1}_{n_1} + \mathbf{X}_3\mathbf{1}_{n_3}) & \mathbf{0}_{d} & ]\\
        \mathbf{B}(\mathbf{S})_2 & = & [ & \mathbf{0}_{d \times n_1} & \mathbf{Y}_2 & \mathbf{Y}_3 & \mathbf{0}_{d} & -(\mathbf{Y}_2\mathbf{1}_{n_2} + \mathbf{Y}_3\mathbf{1}_{n_3}) & ]
    \end{matrix}
\end{equation}
where (see Remark~\ref{rmk:L0DB}) $\mathbf{X}_1 \in \mathbb{R}^{d \times n_1}$ and $\mathbf{X}_3 \in \mathbb{R}^{d \times n_3}$ correspond to the local coordinates, due to the first view, of the $n_1$ points that lie exclusively in the first view and the $n_3$ points that lie on the overlap of both views, respectively. Similarly, $\mathbf{Y}_2 \in \mathbb{R}^{d \times n_2}$ and $\mathbf{Y}_3 \in \mathbb{R}^{d \times n_3}$ correspond to the local coordinates, due to the second view, of the $n_2$ points that lie exclusively in the second view and the $n_3$ points which lie on the overlap of both views, respectively. In particular, $\mathbf{X}_3 = \mathbf{B}(\mathbf{S})_{1,2}$ and $\mathbf{Y}_3 = \mathbf{B}(\mathbf{S})_{2,1}$ (perhaps after permuting the points). Moreover,
\begin{equation}
    \boldsymbol{\mathcal{L}}_{\Gamma} = \begin{bmatrix}
    \mathbf{I}_{n_1} & & & -\mathbf{1}_{n_1} & \mathbf{0}_{n_1}\\
     & \mathbf{I}_{n_2} & & \mathbf{0}_{n_2} & -\mathbf{1}_{n_2}\\
    &  & 2\mathbf{I}_{n_3} & -\mathbf{1}_{n_3} & -\mathbf{1}_{n_3}\\
    -\mathbf{1}_{n_1}^T & \mathbf{0}_{n_2}^T & -\mathbf{1}_{n_3}^T & n_1+n_3 & \\
    \mathbf{0}_{n_1}^T & -\mathbf{1}_{n_2}^T & -\mathbf{1}_{n_3}^T &  & n_2+n_3\end{bmatrix}.
\end{equation}
Through simple calculations, we obtain
\begin{equation}
    \boldsymbol{\mathcal{L}}_{\Gamma}^\dagger = \frac{1}{2n_3}\begin{bmatrix}
    2n_3\mathbf{I}_{n_1} + \mathbf{1}_{n_1}\mathbf{1}_{n_1}^T & -\mathbf{1}_{n_1}\mathbf{1}_{n_2}^T & & \mathbf{1}_{n_1} & -\mathbf{1}_{n_1}\\
    -\mathbf{1}_{n_2}\mathbf{1}_{n_1}^T & 2n_3\mathbf{I}_{n_2} + \mathbf{1}_{n_2}\mathbf{1}_{n_2}^T & & -\mathbf{1}_{n_2} & \mathbf{1}_{n_2}\\
    &  & n_3\mathbf{I}_{n_3} &  & \\
    \mathbf{1}_{n_1}^T & -\mathbf{1}_{n_2}^T &  & 1 & -1\\
    -\mathbf{1}_{n_1}^T & \mathbf{1}_{n_2}^T & &  -1 & 1\end{bmatrix}.
\end{equation}
Thus,
$$\mathbf{B}(\mathbf{S})_1\boldsymbol{\mathcal{L}}_{\Gamma}^\dagger = \begin{bmatrix}\mathbf{X}_1 - \frac{\mathbf{X}_3\mathbf{1}_{n_3}\mathbf{1}_{n_1}^T}{2n_3}, & \frac{\mathbf{X}_3\mathbf{1}_{n_3}\mathbf{1}_{n_2}^T}{2n_3}, & \frac{1}{2}\mathbf{X}_3, & -\frac{\mathbf{X}_3\mathbf{1}_{n_3}}{2n_3}, & \frac{\mathbf{X}_3\mathbf{1}_{n_3}}{2n_3} \end{bmatrix}.$$
Then, using Definition~\ref{def:BSicapj}
\begin{align}
    \mathbf{B}(\mathbf{S})_1\boldsymbol{\mathcal{L}}_{\Gamma}^\dagger \mathbf{B}(\mathbf{S})_2^T &= \frac{1}{2}\mathbf{X}_3\left(\mathbf{I}_{n_3}-\frac{1}{n_3}\mathbf{1}_{n_3}\mathbf{1}_{n_3}^T\right)\mathbf{Y}_3^T\\
    &=\frac{1}{2}\mathbf{B}(\mathbf{S})_{1,2}\left(\mathbf{I}_{n_3}-\frac{1}{n_3}\mathbf{1}_{n_3}\mathbf{1}_{n_3}^T\right)\left(\mathbf{I}_{n_3}-\frac{1}{n_3}\mathbf{1}_{n_3}\mathbf{1}_{n_3}^T\right)^T\mathbf{B}(\mathbf{S})_{2,1}^T\\
    &= \frac{1}{2}\overline{\mathbf{B}(\mathbf{S})}_{1,2}\overline{\mathbf{B}(\mathbf{S})}_{2,1}^T
\end{align}
Since $\mathbf{S} \in \mathcal{C}$ iff $\mathbf{B}(\mathbf{S})_1\boldsymbol{\mathcal{L}}_{\Gamma}^\dagger \mathbf{B}(\mathbf{S})_2^T$ is symmetric, implies $\mathbf{S} \in \mathcal{C}$ iff $\overline{\mathbf{B}(\mathbf{S})}_{1,2}\overline{\mathbf{B}(\mathbf{S})}_{2,1}^T$ is symmetric. 

\noindent \underline{\textbf{Part 2}}. For $\mathbf{S} \in \mathcal{C}$, from the Remark~\ref{rmk:C_hat_L_structure} and Part 1, we have,
\begin{equation}
    \mathbf{L}(\mathbf{S}) = \begin{bmatrix}
    \mathbf{B}(\mathbf{S})_1\boldsymbol{\mathcal{L}}_{\Gamma}^\dagger \mathbf{B}(\mathbf{S})_2^T & -\mathbf{B}(\mathbf{S})_1\boldsymbol{\mathcal{L}}_{\Gamma}^\dagger \mathbf{B}(\mathbf{S})_2^T\\
    -\mathbf{B}(\mathbf{S})_1\boldsymbol{\mathcal{L}}_{\Gamma}^\dagger \mathbf{B}(\mathbf{S})_2^T & \mathbf{B}(\mathbf{S})_1\boldsymbol{\mathcal{L}}_{\Gamma}^\dagger \mathbf{B}(\mathbf{S})_2^T
    \end{bmatrix}. \label{supp:eq:LS_two_views}
\end{equation}
Let $\boldsymbol{\Omega}_1,\boldsymbol{\Omega}_2 \in \Skew(d)$ such that $\boldsymbol{\Omega}_1 + \boldsymbol{\Omega}_2 = 0$. Then, using the above equations, $$\Tr\left(\begin{bmatrix}
        \boldsymbol{\Omega}_1^T & \boldsymbol{\Omega}_2^T
    \end{bmatrix} \mathbf{L}(\mathbf{S}) \begin{bmatrix}
        \boldsymbol{\Omega}_1\\\boldsymbol{\Omega}_2
    \end{bmatrix}\right) = \Tr\left(\begin{bmatrix}
        -\boldsymbol{\Omega}_1 & \boldsymbol{\Omega}_1
    \end{bmatrix} \mathbf{L}(\mathbf{S}) \begin{bmatrix}
        \boldsymbol{\Omega}_1\\-\boldsymbol{\Omega}_1
    \end{bmatrix}\right) = 2 \Tr(\boldsymbol{\Omega}_1^T\overline{\mathbf{B}(\mathbf{S})}_{1,2}\overline{\mathbf{B}(\mathbf{S})}_{2,1}^T\boldsymbol{\Omega}_1).$$
Combining the above and Part 1 with Proposition~\ref{prop:HessFSZZ}, we conclude that $\pi(\mathbf{S})$ is a local minimum of $\widetilde{F}$ iff the first two conditions of the statement are met.

\noindent \underline{\textbf{Part 3}}. Here we deal with the non-degeneracy of $\widetilde{\mathbf{S}} = \pi(\mathbf{S})$. For $d=1$, $\widetilde{\mathbf{S}}$ is trivially non-degenerate. So we assume that $d \geq 2$. From Part 2, we note that for $\boldsymbol{\Omega}_1,\boldsymbol{\Omega}_2 \in \Skew(d)$ such that $\boldsymbol{\Omega}_1 + \boldsymbol{\Omega}_2 = 0$, $\mathbf{L}(\mathbf{S})[\boldsymbol{\Omega}_i]_1^2 = 0$ iff $\overline{\mathbf{B}(\mathbf{S})}_{1,2}\overline{\mathbf{B}(\mathbf{S})}_{2,1}^T\boldsymbol{\Omega} = 0$. Thus $\widetilde{\mathbf{S}}$ is non-degenerate iff $\overline{\mathbf{B}(\mathbf{S})}_{1,2}\overline{\mathbf{B}(\mathbf{S})}_{2,1}^T\boldsymbol{\Omega} = 0 \iff \boldsymbol{\Omega} = 0$. It suffices to show that $\overline{\mathbf{B}(\mathbf{S})}_{1,2}\overline{\mathbf{B}(\mathbf{S})}_{2,1}^T\boldsymbol{\Omega} = 0$  iff $\rank (\overline{\mathbf{B}(\mathbf{S})}_{1,2}\overline{\mathbf{B}(\mathbf{S})}_{2,1}^T) \geq d-1$.

($\impliedby$) Suppose $\rank (\overline{\mathbf{B}(\mathbf{S})}_{1,2}\overline{\mathbf{B}(\mathbf{S})}_{2,1}^T) \geq d-1$ then null space of $\overline{\mathbf{B}(\mathbf{S})}_{1,2}\overline{\mathbf{B}(\mathbf{S})}_{2,1}^T$ is at most one-dimensional. Moreover, rank of a non-zero skew symmetric matrix of size $d \geq 2$, is at least two. Thus $\overline{\mathbf{B}(\mathbf{S})}_{1,2}\overline{\mathbf{B}(\mathbf{S})}_{2,1}^T\boldsymbol{\Omega} = 0 \iff \boldsymbol{\Omega} = 0$. We conclude that
$\widetilde{\mathbf{S}}$ is non-degenerate.

$(\implies)$ Suppose $\rank (\overline{\mathbf{B}(\mathbf{S})}_{1,2}\overline{\mathbf{B}(\mathbf{S})}_{2,1}^T) \leq d-2$, then there exist non-zero vectors $\mathbf{u},\mathbf{v} \in \mathbb{R}^d$ in the kernel of $\overline{\mathbf{B}(\mathbf{S})}_{1,2}\overline{\mathbf{B}(\mathbf{S})}_{2,1}^T$ such that $\mathbf{u}^T\mathbf{v} = 0$. Let $\boldsymbol{\Omega} = \mathbf{u}\mathbf{v}^T - \mathbf{v}\mathbf{u}^T$ then clearly $\boldsymbol{\Omega} \in \Skew(d)$, $\boldsymbol{\Omega} \neq 0$ and $\overline{\mathbf{B}(\mathbf{S})}_{1,2}\overline{\mathbf{B}(\mathbf{S})}_{2,1}^T\boldsymbol{\Omega} = 0$. 

\proofof{Proposition~\ref{prop:noiseless_setting1}}
Since $\mathbf{S}$ is a perfect alignment $F(\mathbf{S}) = \Tr(\mathbf{C}\mathbf{S}\mathbf{S}^T) = 0$. Since $\mathbf{C} \succeq 0$, the columns of $\mathbf{S}$ lie in the kernel of $\mathbf{C}$. In particular $\mathbf{C}\mathbf{S} = \mathbf{0}$. It follows that $\widehat{\mathbf{C}}(\mathbf{S}) = \mathbf{0}$ (see Eq.~(\ref{eq:C_hat})). We conclude that $\mathbf{L}(\mathbf{S}) = \mathbf{C}(\mathbf{S})$.

\proofof{Proposition~\ref{prop:non_deg_views}} \revadd{WLOG assume that each views is centered at the origin i.e. $\mathbf{B}_{i,i}\mathbf{1}_{n_i} = 0$ where $n_i$ is the number of points in the $i$th view. Due to Assumption~\ref{assump:non_deg_views}, the matrix $\mathbf{B}_{i,i}\mathbf{B}_{i,i}^T$ has a rank of $d$ and consequently $\sigma_{\min}(\mathbf{B}_{i,i}\mathbf{B}_{i,i}^T) > 0$. Let $\Theta(\mathbf{S})_i$ and $\Theta(\mathbf{O})_i$ be the realizations of the points in the $i$th views i.e. of $\mathbf{B}_{i,i}$, due to the perfect alignments $\mathbf{S}$ and $\mathbf{O}$, respectively. Also, the optimal translation of the $i$th view due to $\mathbf{S}$ is given by $\mathbf{S}^T\mathbf{B}\boldsymbol{\mathcal{L}}_{\Gamma}^\dagger\mathbf{e}^{n+m}_{n+i}$ and the translation of all the views so that $\Theta(\mathbf{S})$ is centered at the origin is given by $-\mathbf{S}^T\mathbf{B}\boldsymbol{\mathcal{L}}_{\Gamma}^\dagger\mathbf{1}^{n+m}_{n}$. For brevity, define the vector $\mathbf{v}_i\coloneqq \mathbf{B}\boldsymbol{\mathcal{L}}_{\Gamma}^\dagger(\mathbf{e}^{n+m}_{n+i}-\mathbf{1}^{n+m}_{n})$. The net translation for the $i$th view due to $\mathbf{S}$ is the sum of the two translations $\mathbf{S}^T\mathbf{v}_i$ (similarly for $\mathbf{O}$). Then the result follows from,}
\begin{align}
    &\left\|\Theta(\mathbf{S}) - \Theta(\mathbf{O})\right\|^2_F \geq \left\|\Theta(\mathbf{S})_i - \Theta(\mathbf{O})_i\right\|^2_F\\
    &= \left\|(\mathbf{S}_i^T\mathbf{B}_{i,i} + \mathbf{S}^T\mathbf{v}_i\mathbf{1}_{n_i}^T) -(\mathbf{O}_i^T\mathbf{B}_{i,i} + \mathbf{O}^T\mathbf{v}_i\mathbf{1}_{n_i}^T)  \right\|_F^2\\
    &= \left\|(\mathbf{S}_i-\mathbf{O}_i)^T\mathbf{B}_{i,i}\right\|_F^2 + \left\|(\mathbf{S}-\mathbf{O})^T\mathbf{v}_i\mathbf{1}_{n_i}^T\right\|_F^2 - 2\Tr((\mathbf{S}_i-\mathbf{O}_i)^T\mathbf{B}_{i,i}\mathbf{1}_{n_i}\mathbf{v}_i^T(\mathbf{S}-\mathbf{O}))\\
    &\geq \left\|(\mathbf{S}_i-\mathbf{O}_i)^T\mathbf{B}_{i,i}\right\|_F^2 \geq \left\|\mathbf{S}_i-\mathbf{O}_i\right\|_F^2\sigma_{\min}(\mathbf{B}_{i,i}\mathbf{B}_{i,i}^T).
\end{align}

\proofof{Theorem~\ref{thm:inf_rigid}}
\revadd{Suppose $\mathbf{S}$ is degenerate then, due to Theorem~\ref{prop:noiseless_setting1} and Definition~\ref{def:LScertificate}, there exist $\boldsymbol{\Omega}$ such that $\mathbf{C}(\mathbf{S}) \boldsymbol{\Omega} = 0$. Since $\mathbf{C}(\mathbf{S}) \succeq 0$ therefore $\Tr(\mathbf{C}(\mathbf{S}) \boldsymbol{\Omega}\boldsymbol{\Omega}^T) = 0$. Following Eq.~(\ref{eq:opt_Z}), we set the perturbations to be $\mathbf{p}_k \coloneqq \boldsymbol{\Omega}^T\blockdiag(\mathbf{S}^T)\mathbf{B}\boldsymbol{\mathcal{L}}_{\Gamma}^\dagger \mathbf{e}^{n+m}_{k}$ for $k \in [1,n+m]$. Then, following Eq.~(\ref{eq:A_1_}) and Eq.~(\ref{eq:A_1}) in that order, and the fact that $\Tr(\mathbf{C}(\mathbf{S}) \boldsymbol{\Omega}\boldsymbol{\Omega}^T) = 0$,  we obtain $\mathbf{p}_{k} = \boldsymbol{\Omega}_i^T\mathbf{S}_i^T \mathbf{x}_{k,i} + \mathbf{p}_{n+i}$. Consequently, for $(k_1,i), (k_2,i)\in E(\Gamma)$, $\mathbf{p}_{k_1} - \mathbf{p}_{k_2} = \boldsymbol{\Omega}_i^T\mathbf{S}_i^T (\mathbf{x}_{k_1,i}-\mathbf{x}_{k_2,i})$. Similarly, since $\mathbf{S}$ is a perfect alignment, $\mathbf{x}_{k_1}(\mathbf{S}) - \mathbf{x}_{k_2}(\mathbf{S}) = \mathbf{S}_i^T (\mathbf{x}_{k_1,i}-\mathbf{x}_{k_2,i})$. Finally, $(\mathbf{p}_{k_1} - \mathbf{p}_{k_2})^T(\mathbf{x}_{k_1}(\mathbf{S}) - \mathbf{x}_{k_2}(\mathbf{S})) = (\mathbf{x}_{k_1,i}-\mathbf{x}_{k_2,i})^T\mathbf{S}_i\boldsymbol{\Omega}_i\mathbf{S}_i^T (\mathbf{x}_{k_1,i}-\mathbf{x}_{k_2,i}) = 0$ since $\Tr(\boldsymbol{\Omega}_i \mathbf{A}) = 0$ for any symmetric matrix $\mathbf{A}$.}

\revadd{Now suppose $\Theta(\mathbf{S}) = (\mathbf{x}_k(\mathbf{S}))_1^n$ is not infinitesimally rigid. From Definition~\ref{def:inf_rigid}, there exist a non-trivial perturbation $(\mathbf{p}_k)_1^n$ such that $(\mathbf{p}_{k_1}-\mathbf{p}_{k_2})^T(\mathbf{x}_{k_1}(\mathbf{S})-\mathbf{x}_{k_2}(\mathbf{S})) = 0$ for all $(k_1,i), (k_2,i) \in E(\Gamma)$. Combining this with the Assumption~\ref{assump:non_deg_views} that each view is affinely non-degenerate, it follows from  \citea{schulze2010symmetric, asimow1978rigidity} that for each $i \in [1,m]$ there exist $\boldsymbol{\Omega}_i \in \Skew(d)$ and and $\mathbf{t}_i \in \mathbf{R}^d$ such that for each $(k,i) \in E(\Gamma)$, $\mathbf{p}_k = \boldsymbol{\Omega}_i^T\mathbf{x}_k(\mathbf{S}) + \mathbf{t}_i$. Therefore, following Eq.(\ref{eq:A_1}, \ref{eq:A_1_}, \ref{eq:GPOP}), we conclude that $\mathbf{C}(\mathbf{S})\boldsymbol{\Omega} = 0$. From Proposition~\ref{prop:noiseless_setting1} and Proposition~\ref{prop:non_deg_triv_cert}, it suffices to show that $\boldsymbol{\Omega}$ is non-trivial i.e. not all $\boldsymbol{\Omega}_i$'s are equal. In fact, since $\Gamma$ is connected (Assumption~\ref{assump:connected_gamma}) and for $(k,i), (k,j) \in E(\Gamma)$, $\mathbf{p}_k = \boldsymbol{\Omega}_i^T\mathbf{x}_k(\mathbf{S}) + \mathbf{t}_i = \boldsymbol{\Omega}_j^T\mathbf{x}_k(\mathbf{S}) + \mathbf{t}_j$, therefore if $\boldsymbol{\Omega}_i = \boldsymbol{\Omega}_j$ for all $i,j \in [1,m]$ then $\mathbf{t}_i = \mathbf{t}_j$ too. As a result, the perturbation $(\mathbf{p}_k)_1^n$ ends up being trivial, a contradiction.}

\proofof{Theorem~\ref{thm:loc_rigid}}
First, under Assumption~\ref{assump:non_deg_views}, the equation in Definition~\ref{def:loc_rigid}, $\Theta(\mathbf{O}) = \Theta(\mathbf{S}\mathbf{Q})$, is equivalent to $\pi(\mathbf{S}) = \pi(\mathbf{O})$.

($\impliedby$) Suppose $\mathbf{S}$ is a perfect alignment but $\pi(\mathbf{S})$ is not a strict global minimum of $\widetilde{F}$. 
Define 
$$\eta \coloneqq \left\|\mathbf{B}\boldsymbol{\mathcal{L}}_{\Gamma}^\dagger(:,1:n)\left(\mathbf{I}_n - n^{-1}\mathbf{1}_{n}\mathbf{1}_{n}^T\right)\right\|_F$$
and let $\epsilon > 0$ be arbitrary.
Then there exists another perfect alignment $\mathbf{O} \in \mathbb{O}(d)^m$ such that $\left\|\mathbf{S}-\mathbf{O}\right\|_F < \epsilon / \eta$ and $\pi(\mathbf{S}) \neq \pi(\mathbf{O})$. Due to Assumption~\ref{assump:non_deg_views}, we have $\Theta(\mathbf{O}) \neq \Theta(\mathbf{S}\mathbf{Q})$ for any $\mathbf{Q} \in \mathbb{O}(d)$, however, from Definition~\ref{def:realization},
$$\left\|\Theta(\mathbf{O})-\Theta(\mathbf{S})\right\|_F \leq \left\|\mathbf{S}-\mathbf{O}\right\|_F\left\|\mathbf{B}\boldsymbol{\mathcal{L}}_{\Gamma}^\dagger(:,1:n)\left(\mathbf{I}_n - n^{-1}\mathbf{1}_{n}\mathbf{1}_{n}^T\right)\right\|_F = \eta\left\|\mathbf{S}-\mathbf{O}\right\|_F < \epsilon.$$
Since $\epsilon$ is arbitrary, we conclude that $\Theta(\mathbf{S})$ is not locally rigid.

($\implies$) Let $\epsilon > 0$ be arbitrary. Suppose $\Theta(\mathbf{S})$ is not locally rigid, then there exist another perfect alignment $\mathbf{O}_\epsilon \in \mathbb{O}(d)^m$ such that $\left\|\Theta(\mathbf{O}_\epsilon)-\Theta(\mathbf{S})\right\|_F < \epsilon$ but $\Theta(\mathbf{O}_\epsilon) \neq \Theta(\mathbf{S}\mathbf{Q})$. \revadd{Due to Proposition~\ref{prop:non_deg_views}, $\left\|\mathbf{O}_\epsilon-\mathbf{S}\right\|_F < \varrho\epsilon$ where the constant $\varrho > 0$, but $\pi(\mathbf{O}_\epsilon) \neq \pi(\mathbf{S})$}. Since this true for all $\epsilon > 0$, we conclude that $\pi(\mathbf{S})$ is not a strict global minimum of $\widetilde{F}$.

\revadd{Finally, we note that a non-degenerate extremum is a strict extremum and the result follows.}

\proofof{Theorem~\ref{thm:nec_cond_loc_rigid_of_views}}
Consider a partition of $[1,m]$ into two non-empty subsets $A$ and $B$. Suppose $\rank(\overline{\mathbf{B}(\mathbf{S})}_{A, B})$ is at most $d-2$. Let $\boldsymbol{\Omega}_0 \in \Skew(d)$ be such that $\boldsymbol{\Omega}_0 \neq 0$ and $\overline{\mathbf{B}(\mathbf{S})}_{A,B}^T\boldsymbol{\Omega}_0 = 0$ (its existence follows from the third part of the proof of Theorem~\ref{thm:non_deg_two_views_gen_setting}). WLOG assume that $\mathbf{B}(\mathbf{S})_{A,B}\mathbf{1}_{n'} = 0$ (here $n'$ is as in Definition~\ref{def:BSAcapB}) (perhaps by translating all aligned views by $-\mathbf{B}(\mathbf{S})_{A,B}\mathbf{1}_{n'}$). Then $\mathbf{B}(\mathbf{S})_{A,B}^T\boldsymbol{\Omega}_0 = 0$.

Let $\boldsymbol{\Omega} = [\boldsymbol{\Omega}_i]_1^m$ be such that $\boldsymbol{\Omega}_i = \boldsymbol{\Omega}_0$ for $i \in A$ and $\boldsymbol{\Omega}_i = -\boldsymbol{\Omega}_0$ for $i \in B$. Clearly, $\boldsymbol{\Omega} \in \Skew(d)^m$ such that not all $\boldsymbol{\Omega}_i$ are equal. It suffices to show that $\boldsymbol{\Omega}$ is a nontrivial certificate of $\mathbf{L}(\mathbf{S})$, equivalently $\Tr(\boldsymbol{\Omega}^T\mathbf{L}(\mathbf{S})\boldsymbol{\Omega}) = 0$. First we observe that for $i \in A$, 
$$[\mathbf{L}(\mathbf{S}) \boldsymbol{\Omega}]_i = (\mathbf{B}(\mathbf{S})_i\boldsymbol{\mathcal{L}}_{\Gamma}^\dagger \mathbf{B}(\mathbf{S})^T\mathbf{I}^m_d - \textstyle\sum_{1}^{m}(-1)^{\mathbf{1}_{B}(j)}\mathbf{B}(\mathbf{S})_i\boldsymbol{\mathcal{L}}_{\Gamma}^\dagger \mathbf{B}(\mathbf{S})_j^T)\boldsymbol{\Omega}_0 = 2 \textstyle\sum_{j \in B}\mathbf{B}(\mathbf{S})_i\boldsymbol{\mathcal{L}}_{\Gamma}^\dagger \mathbf{B}(\mathbf{S})_j^T \boldsymbol{\Omega}_0$$
where $\mathbf{1}_{B}(j) = 1$ iff $j \in B$.
Similarly, for $i \in B$, 
$$[\mathbf{L}(\mathbf{S}) \boldsymbol{\Omega}]_i = -2 \textstyle\sum_{j \in A}\mathbf{B}(\mathbf{S})_i\boldsymbol{\mathcal{L}}_{\Gamma}^\dagger \mathbf{B}(\mathbf{S})_j^T \boldsymbol{\Omega}_0.$$ Denote by $\mathbf{B}_A$ and $\mathbf{B}_B$, the matrices $\textstyle\sum_{i \in A}\mathbf{B}(\mathbf{S})_i$ and $\textstyle\sum_{j \in B}\mathbf{B}(\mathbf{S})_j$, respectively. Thus,
\begin{equation}
    \Tr(\boldsymbol{\Omega}^T\mathbf{L}(\mathbf{S})\boldsymbol{\Omega}) = 4\Tr(\boldsymbol{\Omega}_0^T \mathbf{B}_A\boldsymbol{\mathcal{L}}_{\Gamma}^\dagger \mathbf{B}_B^T \boldsymbol{\Omega}_0). \label{supp:eq:Omega0TLSOmega0}
\end{equation}

We are going to show that the above evaluates to zero. WLOG assume that the first $n_1$ points lie in the views with indices in $A \setminus B$, next $n_2$ points  lie in the views with indices in $B \setminus A$ and the remaining $n_3$ points lie in the views with indices in $A \cap B$. Note that $n_1 + n_2 + n_3 = n$ and $|A| + |B| = m$. Then the matrices $\mathbf{B}_A$ and $\mathbf{B}_B$ (perhaps after permuting the views) have the following structure. 
\begin{align}
    \begin{matrix}
        \mathbf{B}_A & = & [ & \mathbf{X}_{1} & \mathbf{0}_{d \times n_2} & \mathbf{X}_{3} & \mathbf{U}_{1} + \mathbf{U}_{3} & \mathbf{0}_{d \times |B|} & ]\\
        \mathbf{B}_B & = & [ & \mathbf{0}_{d \times n_1} & \mathbf{Y}_{2} & \mathbf{Y}_{3} & \mathbf{0}_{d \times |A|} & \mathbf{V}_{2} + \mathbf{V}_{3} & ]
    \end{matrix}
\end{align}
where $\mathbf{X}_{1} \in \mathbb{R}^{d \times n_1}$ and $\mathbf{Y}_{2} \in \mathbb{R}^{d \times n_2}$ contain the sum of the local coordinates of the $n_1$ and $n_2$ points, respectively. Also, $\mathbf{X}_{3} \in \mathbb{R}^{d \times n_3}$ and $\mathbf{Y}_{3} \in \mathbb{R}^{d \times n_3}$ contain the sum of the local coordinates of the remaining $n_3$ points due to the views with indices in $A$ and $B$ respectively. The matrices $\mathbf{U}_{1} \in \mathbb{R}^{d \times |A|}$, $\mathbf{V}_{2} \in \mathbb{R}^{d \times |B|}$, $\mathbf{U}_{3} \in \mathbb{R}^{d \times |A|}$ and $\mathbf{V}_{3} \in \mathbb{R}^{d \times |B|}$ follow from Remark~\ref{rmk:L0DB}. Further define
\begin{equation}
    \begin{matrix}
         \mathbf{B}_{A \setminus B} & \coloneqq & [ & \mathbf{X}_{1} & 0 & 0 & \mathbf{U}_{1} & 0 & ]\\
         \mathbf{B}_{A,B} & \coloneqq & [ & 0 & 0 & \mathbf{X}_{3} & \mathbf{U}_{3} & 0 & ]\\
         \mathbf{B}_{B \setminus A} & \coloneqq & [ & 0 & \mathbf{Y}_{2} & 0 & 0 & \mathbf{V}_{2} & ]\\
         \mathbf{B}_{B,A} & \coloneqq & [ & 0 & 0& \mathbf{Y}_{3} & 0 & \mathbf{V}_{3} & ]
    \end{matrix} \label{supp:eq:eq8}
\end{equation}
then $\mathbf{B}_A = \mathbf{B}_{A \setminus B} + \mathbf{B}_{A,B}$ and $\mathbf{B}_B = \mathbf{B}_{B \setminus A} + \mathbf{B}_{B,A}$. Note that $\mathbf{1}_{n+m}$ lies in the kernel of the four matrices defined above (see Remark~\ref{rmk:L0DB}) and, $\mathbf{B}_{A \setminus B}\mathbf{B}_{B\setminus A}^T = 0$, $\mathbf{B}_{A \setminus B}\mathbf{B}_{B,A}^T = 0$ and $\mathbf{B}_{B \setminus A}\mathbf{B}_{A,B}^T = 0$. Now, the structure of $\boldsymbol{\mathcal{L}}_{\Gamma}$ is as follows,
\begin{equation}
    \boldsymbol{\mathcal{L}}_{\Gamma} = \begin{bmatrix}
        \boldsymbol{\mathcal{D}}_{1} & & & -\boldsymbol{\mathcal{K}}_{1} & \\
        & \boldsymbol{\mathcal{D}}_{2} & &  & -\boldsymbol{\mathcal{K}}_{2}\\
        & &\boldsymbol{\mathcal{D}}_{3} + \boldsymbol{\mathcal{D}}_{3} & -\boldsymbol{\mathcal{K}}_{A_3} & -\boldsymbol{\mathcal{K}}_{B_3}\\
        -\boldsymbol{\mathcal{K}}_{1}^T & & -\boldsymbol{\mathcal{K}}_{A_3}^T & \boldsymbol{\mathcal{D}}_{A} & \\
        & -\boldsymbol{\mathcal{K}}_{2}^T & -\boldsymbol{\mathcal{K}}_{B_3}^T & & \boldsymbol{\mathcal{D}}_{B}
    \end{bmatrix}
\end{equation}
where $\boldsymbol{\mathcal{K}}_{1} \in \mathbb{R}^{n_1 \times |A|}$ is the adjacency between the first $n_1$ points and the views with indices in $A$, $\boldsymbol{\mathcal{K}}_{2} \in \mathbb{R}^{n_2 \times |B|}$ is the adjacency between the next $n_2$ points and the views with indices in $B$, and $\boldsymbol{\mathcal{K}}_{A_3} \in \mathbb{R}^{n_3 \times |A|}$ and $\boldsymbol{\mathcal{K}}_{B_3} \in \mathbb{R}^{n_3 \times |B|}$ are the adjacencies between the remaining $n_3$ points and the views with indices in $A$ and $B$ respectively. As for the remaining matrices, $\boldsymbol{\mathcal{D}}_{i} = \diag(\boldsymbol{\mathcal{K}}_{i}\mathbf{1}_{|A|})$ represents the degrees of the points in the bipartite adjacency $\boldsymbol{\mathcal{K}}_{i}$, $\boldsymbol{\mathcal{D}}_{A} = \diag(\boldsymbol{\mathcal{K}}_{1}^T\mathbf{1}_{n_1} + \boldsymbol{\mathcal{K}}_{A_3}^T\mathbf{1}_{n_3})$ represents the degree of the views i.e. the number of points contained in the views with indices in $A$, Similarly, $\boldsymbol{\mathcal{D}}_{j} = \diag(\boldsymbol{\mathcal{K}}_{j}\mathbf{1}_{|B|})$. and $\boldsymbol{\mathcal{D}}_{B} = \diag(\boldsymbol{\mathcal{K}}_{2}^T\mathbf{1}_{n_2} + \boldsymbol{\mathcal{K}}_{B_3}^T\mathbf{1}_{n_3})$. 

Since $\mathbf{S}$ is a perfect alignment, the local coordinates of a point due to the views are the same. Using the fact that $\mathbf{B}(\mathbf{S})_{A,B}$ represent the local coordinates of the $n_3$ points contained in views with indices in $A \cap B$ (see Definition~\ref{def:BSAcapB}), we obtain
\begin{align}
    \begin{matrix}
        \mathbf{X}_{3} &=& \mathbf{B}(\mathbf{S})_{A,B}\boldsymbol{\mathcal{D}}_{3}, & & \mathbf{Y}_{3} &=& \mathbf{B}(\mathbf{S})_{A,B}\boldsymbol{\mathcal{D}}_{3},\\
        \mathbf{U}_{3} &=& -\mathbf{B}(\mathbf{S})_{A,B}\boldsymbol{\mathcal{K}}_{A_3},& &
        \mathbf{V}_{3} &=& -\mathbf{B}(\mathbf{S})_{A,B}\boldsymbol{\mathcal{K}}_{B_3}.
    \end{matrix}\label{supp:eq:eq9}
\end{align} 
Similarly, it follows that $\mathbf{U}_1 = \mathbf{X}_1\boldsymbol{\mathcal{D}}_{1}^{-1}\boldsymbol{\mathcal{K}}_{1}$ and $\mathbf{V}_2 = \mathbf{Y}_2\boldsymbol{\mathcal{D}}_{2}^{-1}\boldsymbol{\mathcal{K}}_{2}$. Thus,
\begin{align}
    \begin{matrix}
         \mathbf{B}_{A \setminus B} & = & [ & \mathbf{X}_{1}\boldsymbol{\mathcal{D}}_{1}^{-1} & \mathbf{0}_{d \times n_2} & \mathbf{0}_{d \times n_3} & \mathbf{0}_{d \times |A|} & \mathbf{0}_{d \times |B|} & ]\boldsymbol{\mathcal{L}}_{\Gamma}\\
         \mathbf{B}_{B \setminus A} & = & [ & \mathbf{0}_{d \times n_1} & \mathbf{Y}_{2}\boldsymbol{\mathcal{D}}_{2}^{-1} & \mathbf{0}_{d \times n_3} & \mathbf{0}_{d \times |A|} & \mathbf{0}_{d \times |B|} & ]\boldsymbol{\mathcal{L}}_{\Gamma}
    \end{matrix}\label{supp:eq:eq14_} \text{ and }
\end{align}
\begin{align}
    \begin{matrix}
         \mathbf{B}_{A \setminus B}\boldsymbol{\mathcal{L}}_{\Gamma}^\dagger & = & [ & \mathbf{X}_{1}\boldsymbol{\mathcal{D}}_{1}^{-1} & \mathbf{0}_{d \times n_2} & \mathbf{0}_{d \times n_3} & \mathbf{0}_{d \times |A|} & \mathbf{0}_{d \times |B|} & ] + \mathbf{t}_{A \setminus B}\mathbf{1}_{n+m}^T\\
         \mathbf{B}_{B \setminus A}\boldsymbol{\mathcal{L}}_{\Gamma}^\dagger & = & [ & \mathbf{0}_{d \times n_1} & \mathbf{Y}_{2}\boldsymbol{\mathcal{D}}_{2}^{-1} & \mathbf{0}_{d \times n_3} & \mathbf{0}_{d \times |A|} & \mathbf{0}_{d \times |B|} & ] +  \mathbf{t}_{B \setminus A}\mathbf{1}_{n+m}^T
    \end{matrix}\label{supp:eq:eq14}
\end{align}
for some translation vectors $\mathbf{t}_{A \setminus B}, \mathbf{t}_{B \setminus A} \in \mathbb{R}^d$. Since $\mathbf{1}_{n+m}$ lies in $\ker(\mathbf{B}_{A \setminus B})$ and $\ker(\mathbf{B}_{A \setminus B})$ (see the paragraph after Eq.~(\ref{supp:eq:eq8})), thus
\begin{equation}
    \mathbf{B}_{A \setminus B}\boldsymbol{\mathcal{L}}_{\Gamma}^\dagger \mathbf{B}_{B\setminus A}^T = \mathbf{B}_{A \setminus B}\boldsymbol{\mathcal{L}}_{\Gamma}^\dagger \mathbf{B}_{B,A}^T = \mathbf{B}_{B \setminus A}\boldsymbol{\mathcal{L}}_{\Gamma}^\dagger \mathbf{B}_{A,B}^T = 0. \label{supp:eq:eq11}
\end{equation}

Since $\mathbf{B}(\mathbf{S})_{A,B}^T\boldsymbol{\Omega}_0 = 0$ (by assumption), combining with Eq.~(\ref{supp:eq:eq9}, \ref{supp:eq:eq8}) yields $\mathbf{B}_{A,B}^T\boldsymbol{\Omega}_0 = 0$ and $ \mathbf{B}_{B,A}^T\boldsymbol{\Omega}_0 = 0$. Substituting the above and Eq.~(\ref{supp:eq:eq11}) into Eq.~(\ref{supp:eq:Omega0TLSOmega0}), we obtain 
$$\Tr(\boldsymbol{\Omega}_0^T \mathbf{B}_A\boldsymbol{\mathcal{L}}_{\Gamma}^\dagger \mathbf{B}_B^T \boldsymbol{\Omega}_0) = \Tr(\boldsymbol{\Omega}_0^T (\mathbf{B}_{A \setminus B}+\mathbf{B}_{A,B})\boldsymbol{\mathcal{L}}_{\Gamma}^\dagger (\mathbf{B}_{B\setminus A} + \mathbf{B}_{B,A})^T \boldsymbol{\Omega}_0)
    = \Tr(\boldsymbol{\Omega}_0^T \mathbf{B}_{A,B}\boldsymbol{\mathcal{L}}_{\Gamma}^\dagger\mathbf{B}_{B,A}^T \boldsymbol{\Omega}_0) = 0.$$
We conclude that $\boldsymbol{\Omega}$ is a non-trivial certificate of $\mathbf{L}(\mathbf{S})$ and thus $\pi(\mathbf{S})$ is degenerate.

\proofof{Lemma~\ref{lem:subproblem_cert}} 
WLOG, let the $m$th vertex be removed. Let $\Gamma$ be the bipartite graph representing the correspondence between $m$ views and $n$ vertices, as described in Section~\ref{sec:setup}. Let $\Gamma_{-}$ be the bipartite graph obtained after the removal of the vertices representing the $m$th view and the points which lie exclusively in it. Let $\mathbf{D} \in \mathbb{R}^{md \times md}$, $\mathbf{B} \in \mathbb{R}^{md \times (n+m)}$, $\mathbf{D}_{-} \in \mathbb{R}^{(m-1)d\times (m-1)d}$ and $\mathbf{B}_{-} \in \mathbb{R}^{(m-1)d \times (n_1+n_2+m-1)}$ be the matrices defined in Remark~\ref{rmk:L0DB} for graphs $\Gamma$ and $\Gamma_{-}$. Also, let
\begin{itemize}[leftmargin=*]
    \item $\boldsymbol{\mathcal{K}}_1 \in \mathbb{R}^{n_1 \times (m-1)}$ is the bipartite adjacency matrix between the first $m-1$ views and the $n_1$ points which lie exclusively in them. Note that the adjacency between such points and the $m$th view is $\mathbf{0}_{n_1}$.
    \item $\boldsymbol{\mathcal{K}}_2 \in \mathbb{R}^{n_2 \times (m-1)}$ is the bipartite adjacency between the first $m-1$ views and the $n_2$ points which lie on the overlap of the $m$th view and the union of the first $m-1$ views. Note that the adjacency between such points and the $m$th view is $\mathbf{1}_{n_2}$. Also note that since $\Gamma$ is connected by Assumption~\ref{assump:connected_gamma}, $n_2 > 0$.
    \item the fifth and the third column in $\boldsymbol{\mathcal{L}}_{\Gamma}$ correspond to the  $m$th view and the $n_3$ points that lie exclusively in it, respectively. The adjacency between such points and the first $m-1$ views is $\mathbf{0}_{n_3 \times (m-1)}$, and that with the $m$th view is $\mathbf{1}_{n_3}$.
    \item $\boldsymbol{\mathcal{D}}_1 = \diag (\boldsymbol{\mathcal{K}}_1\mathbf{1}_{m-1})$, $\boldsymbol{\mathcal{D}}_2 = \diag (\boldsymbol{\mathcal{K}}_2\mathbf{1}_{m-1})$ and $\overline{\boldsymbol{\mathcal{D}}} = \diag (\boldsymbol{\mathcal{K}}_1^T\mathbf{1}_{n_1} + \boldsymbol{\mathcal{K}}_2^T\mathbf{1}_{n_2})$.
\end{itemize}
Then the structure of the combinatorial Laplacian of $\Gamma$ and $\Gamma_{-}$ are
\begin{equation}
    \boldsymbol{\mathcal{L}}_{\Gamma} = \begin{bmatrix}
        \boldsymbol{\mathcal{D}}_1 &  &  & -\boldsymbol{\mathcal{K}}_1 & \mathbf{0}_{n_1}\\
         & \boldsymbol{\mathcal{D}}_2 + \mathbf{I}_{n_2} &  & -\boldsymbol{\mathcal{K}}_2 & -\mathbf{1}_{n_2}\\
         &  & \mathbf{I}_{n_3} & \mathbf{0}_{n_3 \times (m-1)} & -\mathbf{1}_{n_3}\\
        -\boldsymbol{\mathcal{K}}_1^T & -\boldsymbol{\mathcal{K}}_2^T & \mathbf{0}_{n_3}^T & \overline{\boldsymbol{\mathcal{D}}} & \mathbf{0}_{m-1}\\
        \mathbf{0}_{n_1}^T & -\mathbf{1}_{n_2}^T & -\mathbf{1}_{n_3}^T & \mathbf{0}_{m-1}^T & n_2+n_3
    \end{bmatrix}
\end{equation}
and
$\boldsymbol{\mathcal{L}}_{\Gamma_{-}} = \begin{bsmallmatrix}
    \boldsymbol{\mathcal{D}}_1 &  & -\boldsymbol{\mathcal{K}}_1\\
     & \boldsymbol{\mathcal{D}}_2 & -\boldsymbol{\mathcal{K}}_2\\
    -\boldsymbol{\mathcal{K}}_1^T & -\boldsymbol{\mathcal{K}}_2^T & \overline{\boldsymbol{\mathcal{D}}}
    \end{bsmallmatrix}$.
Using a permutation matrix
$\boldsymbol{\mathcal{P}}_{0} = \begin{bsmallmatrix}
        \mathbf{I}_{n_1+n_2} &  &  & \\
        & & \mathbf{I}_{n_3}& \\
        & \mathbf{I}_{m-1}& & \\
        & & & 1
    \end{bsmallmatrix}$
and a diagonal matrix $\boldsymbol{\mathcal{D}}_0 = \diag ((\mathbf{0}_{n_1},\mathbf{1}_{n_2},\mathbf{0}_{m-1}))$, we obtain
\begin{equation}
    \boldsymbol{\mathcal{P}}_{0}\boldsymbol{\mathcal{L}}_{\Gamma}\boldsymbol{\mathcal{P}}_{0}^T = \begin{bmatrix}\mathbf{A}_{11}&\mathbf{A}_{12}\\\mathbf{A}_{21}&\mathbf{A}_{22}\end{bmatrix} = \begin{bmatrix}
        &&& \mathbf{0}_{n_1} & \mathbf{0}_{n_1}\\
        & \boldsymbol{\mathcal{L}}_{\Gamma_{-}} + \boldsymbol{\mathcal{D}}_0 & & \mathbf{0}_{n_2} & -\mathbf{1}_{n_2}\\
        &&& \mathbf{0}_{m-1} & \mathbf{0}_{m-1}\\
        \mathbf{0}_{n_1}^T & \mathbf{0}_{n_2}^T & \mathbf{0}_{m-1}^T & \mathbf{I}_{n_3} & -\mathbf{1}_{n_3}\\
        \mathbf{0}_{n_1}^T & -\mathbf{1}_{n_2}^T & \mathbf{0}_{m-1}^T & -\mathbf{1}_{n_3}^T & n_2+n_3
    \end{bmatrix} \label{supp:eq:eq15}
\end{equation}

The rest is divided into three parts. First, we derive the pseudoinverse of the above block matrix using \citeb[Section 3.6.2]{gentle2007matrix}. Then we show that $\mathbf{S}_{-} \coloneqq \mathbf{S}_{-m}$ is a perfect alignment of the $m-1$ views and finally we show that $[\boldsymbol{\Omega}_i]_{1}^{m-1}$ is a certificate of $\mathbf{L}_{-}(\mathbf{S}_{-}) \coloneqq \mathbf{L}_{-m}(\mathbf{S}_{-m})$ when $[\boldsymbol{\Omega}_i]_{1}^{m}$ is a certificate of $\mathbf{L}(\mathbf{S})$.

\noindent \underline{\textbf{Part 1}}. Here we derive the pseudoinverse of the matrix in Eq.~(\ref{supp:eq:eq15}). First, we note
\begin{prop}
\label{supp:prop:LplusD0}
$\boldsymbol{\mathcal{L}}_{\Gamma_{-}} + \boldsymbol{\mathcal{D}}_0 \succ 0$.
\end{prop}
\textit{Proof}. Since $\boldsymbol{\mathcal{L}}_{\Gamma_{-}} \succeq 0$ and $\boldsymbol{\mathcal{D}}_0 \succeq 0$, it suffices to show that $\ker (\boldsymbol{\mathcal{L}}_{\Gamma_{-}}) \cap \ker (\boldsymbol{\mathcal{D}}_0) = \{0\}$. Recall that the $m$th view contains $n_2 + n_3$ points where $n_2>0$ points lie on the overlap of $m$th view and the union of first $m-1$ views, and $n_3$ points lie exclusively in the $m$th view. Removal of the $m$th view and the $n_3$ points that lie exclusively in it may disconnect $\Gamma$ to produce $\Gamma_{-}$ with at most $n_2$ connected components. The vectors $\mathbf{u}_i$ with ones at the indices of the vertices in the $i$th component and zeros elsewhere, form an orthogonal basis of $\ker(\boldsymbol{\mathcal{L}}_{\Gamma_{-}})$. Since there exists at least one $k \in [n_1+1, n_1+n_2]$ with $\mathbf{u}_i(k) = 1$, thus $\mathbf{u}_i^T\boldsymbol{\mathcal{D}}_0\mathbf{u}_i > 0$. Also, for $i \neq j$, $\mathbf{u}_i^T\boldsymbol{\mathcal{D}}_0\mathbf{u}_j = 0$. The result follows.~$\blacksquare$

Since $\boldsymbol{\mathcal{L}}_{\Gamma_{-}} + \boldsymbol{\mathcal{D}}_0 \succ 0$, thus $(\boldsymbol{\mathcal{L}}_{\Gamma_{-}}+\boldsymbol{\mathcal{D}}_0)^\dagger = (\boldsymbol{\mathcal{L}}_{\Gamma_{-}}+\boldsymbol{\mathcal{D}}_0)^{-1}$ and
\begin{equation}
    (\boldsymbol{\mathcal{L}}_{\Gamma_{-}}+\boldsymbol{\mathcal{D}}_0)  \begin{bmatrix}\mathbf{1}_{n_1}\\ \mathbf{1}_{n_2} \\ \mathbf{1}_{m-1}\end{bmatrix} =  \begin{bmatrix}\mathbf{0}_{n_1}\\ \mathbf{1}_{n_2} \\ \mathbf{0}_{m-1}\end{bmatrix} \implies (\boldsymbol{\mathcal{L}}_{\Gamma_{-}}+\boldsymbol{\mathcal{D}}_0)^\dagger \begin{bmatrix}\mathbf{0}_{n_1}\\ \mathbf{1}_{n_2} \\ \mathbf{0}_{m-1}\end{bmatrix} = \begin{bmatrix}\mathbf{1}_{n_1}\\ \mathbf{1}_{n_2} \\ \mathbf{1}_{m-1}\end{bmatrix}. \label{supp:eq:L_Gamma_minus__plus_D_0_1s}
\end{equation}
Using the above equation, the matrix $\mathbf{Z} \coloneqq [(\boldsymbol{\mathcal{P}}_{0}\boldsymbol{\mathcal{L}}_{\Gamma}\boldsymbol{\mathcal{P}}_{0}^T)^\dagger]_{22} = \mathbf{A}_{22} - \mathbf{A}_{21}\mathbf{A}_{11}^\dagger \mathbf{A}_{12}$ and its pseudoinverse are, $\mathbf{Z} = \begin{bmatrix}
    \mathbf{I}_{n_3} & -\mathbf{1}_{n_3}\\
    -\mathbf{1}_{n_3}^T & n_3
\end{bmatrix}$ and $\mathbf{Z}^\dagger = \begin{bmatrix}
    \mathbf{I}_{n_3}  & \mathbf{0}_{n_3}\\
    \mathbf{0}_{n_3}^T & 0
\end{bmatrix}$. Next, we have 
$$[(\boldsymbol{\mathcal{P}}_{0}\boldsymbol{\mathcal{L}}_{\Gamma}\boldsymbol{\mathcal{P}}_{0}^T)^\dagger]_{11} = (\boldsymbol{\mathcal{L}}_{\Gamma_{-}}+\boldsymbol{\mathcal{D}}_0)^\dagger + ((\boldsymbol{\mathcal{L}}_{\Gamma_{-}}+\boldsymbol{\mathcal{D}}_0)^\dagger \mathbf{A}_{12})\mathbf{Z}^{\dagger}(\mathbf{A}_{21}(\boldsymbol{\mathcal{L}}_{\Gamma_{-}}+\boldsymbol{\mathcal{D}}_0)^\dagger).$$
Using Eq.~(\ref{supp:eq:eq15}, \ref{supp:eq:L_Gamma_minus__plus_D_0_1s}),  we obtain,
\begin{align}
[(\boldsymbol{\mathcal{P}}_{0}\boldsymbol{\mathcal{L}}_{\Gamma}\boldsymbol{\mathcal{P}}_{0}^T)^\dagger]_{11} &= (\boldsymbol{\mathcal{L}}_{\Gamma_{-}}+\boldsymbol{\mathcal{D}}_0)^\dagger. \label{supp:eq:pinvA11}\\
[(\boldsymbol{\mathcal{P}}_{0}\boldsymbol{\mathcal{L}}_{\Gamma}\boldsymbol{\mathcal{P}}_{0}^T)^\dagger]_{12} &= -(\mathbf{A}_{11}^\dagger \mathbf{A}_{12}) \mathbf{Z}^\dagger = 0\\
[(\boldsymbol{\mathcal{P}}_{0}\boldsymbol{\mathcal{L}}_{\Gamma}\boldsymbol{\mathcal{P}}_{0}^T)^\dagger]_{21} &= 0.
\end{align}
Thus,
$$(\boldsymbol{\mathcal{P}}_{0}\boldsymbol{\mathcal{L}}_{\Gamma}\boldsymbol{\mathcal{P}}_{0}^T)^\dagger = \blockdiag((\boldsymbol{\mathcal{L}}_{\Gamma_{-}}+\boldsymbol{\mathcal{D}}_0)^\dagger, \mathbf{I}_{n_3}, 0).$$

\noindent \underline{\textbf{Part 2}}. Now, let $\mathbf{S} = [\mathbf{S}_i]_1^{m}$ and $\mathbf{S}_{-} = [\mathbf{S}_i]_1^{m-1}$. Intuitively, it should be clear that $\mathbf{S}_{-}$ is a perfect alignment for the $m-1$ views. Since $\mathbf{S}$ is a perfect alignment, the alignment error (see Eq.~(\ref{eq:GPOP}))
$$\Tr(\mathbf{S}^T(\mathbf{D}-\mathbf{B}\boldsymbol{\mathcal{L}}_{\Gamma}^\dagger\mathbf{B}^T)\mathbf{S}) = 0.$$ We show that the error after the removal of the $m$th view is still zero i.e. 
$$\Tr(\mathbf{S}_{-}^T(\mathbf{D}_{-}-\mathbf{B}_{-}\boldsymbol{\mathcal{L}}_{\Gamma_{-}}\mathbf{B}_{-}^T)\mathbf{S}_{-}) = 0.$$
Let $\mathbf{B}^*_{1} \in \mathbb{R}^{d \times n_1}$, $\mathbf{B}^*_{2} \in \mathbb{R}^{d \times n_2}$ and $\mathbf{B}^*_{3} \in \mathbb{R}^{d \times n_3}$ contain the coordinates (after alignment with $\mathbf{S}$) of the $n_1$ points that lie exclusively in the first $m-1$ views, of the $n_2$ points that lie on the overlap of the $m$th view with the remaining views, and of the $n_3$ points that lie exclusively in the $m$th view, respectively. Then it suffices to show
\begin{prop}
\label{supp:prop:subproblem_perf_alignment}
(i) $\mathbf{S}^T\mathbf{B}\boldsymbol{\mathcal{L}}_{\Gamma}^\dagger\mathbf{B}^T\mathbf{S} =  \mathbf{S}_{-}^T\mathbf{B}_{-}\boldsymbol{\mathcal{L}}_{\Gamma_{-}}^\dagger \mathbf{B}_{-}^T\mathbf{S}_{-} + \mathbf{B}^*_{2}\mathbf{B}^{*^T}_{2} + \mathbf{B}^*_{3}\mathbf{B}^{*^T}_{3}$ and (ii) $\mathbf{S}^T\mathbf{D}\mathbf{S} = \mathbf{S}_{-}^T\mathbf{D}_{-}\mathbf{S}_{-} + \mathbf{B}^*_{2}\mathbf{B}^{*^T}_{2} + \mathbf{B}^*_{3}\mathbf{B}^{*^T}_{3}$. By taking the trace of the difference of these equations, the main result follows.
\end{prop}
\textit{Proof}. The second equation follows from Eq.~(\ref{eq:D}), Remark~\ref{rmk:L0DB} and the fact that, since the $m$ views are perfectly aligned, the local coordinates of the points are the same as those in the matrices $\mathbf{B}^*_{1}$,  $\mathbf{B}^*_{2}$ and  $\mathbf{B}^*_{3}$. We proceed to prove the first equation.

Since $\mathbf{\mathcal{P}}_0$ is a permutation matrix, $\mathbf{S}^T\mathbf{B}\boldsymbol{\mathcal{L}}_{\Gamma}^\dagger\mathbf{B}^T\mathbf{S} = (\mathbf{S}^T\mathbf{B}\boldsymbol{\mathcal{L}}_{\Gamma}^\dagger\boldsymbol{\mathcal{P}}_{0}^T)(\boldsymbol{\mathcal{P}}_{0}\mathbf{B}^T\mathbf{S})$. Then, using the same idea as in Eq.~(\ref{supp:eq:eq9}, \ref{supp:eq:eq14_}, \ref{supp:eq:eq14}) in the proof of Theorem~\ref{thm:nec_cond_loc_rigid_of_views}, we obtain
\begin{align}
    \mathbf{S}^T\mathbf{B} &= \begin{bmatrix}
        \mathbf{B}^*_{1}\boldsymbol{\mathcal{D}}_1 & \mathbf{B}^*_{2}(\boldsymbol{\mathcal{D}}_2 + \mathbf{I}_{n_2}) & \mathbf{B}^*_{3} & -(\mathbf{B}^*_{1}\boldsymbol{\mathcal{K}}_1+\mathbf{B}^*_{2}\boldsymbol{\mathcal{K}}_2) & -(\mathbf{B}^*_{2}\mathbf{1}_{n_2}+\mathbf{B}^*_{3}\mathbf{1}_{n_3})
    \end{bmatrix}\\
    &= \begin{bmatrix}
        \mathbf{B}^*_{1} & \mathbf{B}^*_{2} & \mathbf{B}^*_{3} & \mathbf{0}_{d \times (m-1)} & \mathbf{0}_{d}
    \end{bmatrix}\boldsymbol{\mathcal{L}}_{\Gamma}\label{supp:eq:STB}
\end{align}
and thus,
\begin{equation}
    \mathbf{S}^T\mathbf{B}\boldsymbol{\mathcal{L}}_{\Gamma}^\dagger = \begin{bmatrix}
        \mathbf{B}^*_{1} & \mathbf{B}^*_{2} & \mathbf{B}^*_{3} & \mathbf{0}_{d \times (m-1)} & \mathbf{0}_{d}
    \end{bmatrix} + \mathbf{t}\mathbf{1}_{n+m}^T \label{supp:eq:STBL_GammaBTS}
\end{equation}
for some translation vector $\mathbf{t} \in \mathbb{R}^d$. Similarly,
\begin{equation}
    \mathbf{S}_{-}^T\mathbf{B}_{-} = \begin{bmatrix}
        \mathbf{B}^*_{1}\boldsymbol{\mathcal{D}}_1 & \mathbf{B}^*_{2} \boldsymbol{\mathcal{D}}_2 & -(\mathbf{B}^*_{1}\boldsymbol{\mathcal{K}}_1+\mathbf{B}^*_{2}\boldsymbol{\mathcal{K}}_2)
    \end{bmatrix} = \begin{bmatrix}
        \mathbf{B}^*_{1} & \mathbf{B}^*_{2} & \mathbf{0}_{d \times (m-1)}
    \end{bmatrix}\boldsymbol{\mathcal{L}}_{\Gamma_{-}} \label{supp:eq:SmTBm1}
\end{equation}
and thus 
$$\mathbf{S}_{-}^T\mathbf{B}_{-}\boldsymbol{\mathcal{L}}_{\Gamma_{-}}^\dagger  = \begin{bmatrix}
        \mathbf{B}^*_{1} & \mathbf{B}^*_{2} & \mathbf{0}_{d \times (m-1)}
    \end{bmatrix} + \mathbf{t}_{-}\mathbf{v}_{-}^T$$    
for some translation vector $\mathbf{t}_{-} \in \mathbb{R}^d$ and $\mathbf{v}_{-} \in \ker(\boldsymbol{\mathcal{L}}_{\Gamma_{-}})$. From Proposition~\ref{prop:kerB}, $\ker(\boldsymbol{\mathcal{L}}_{\Gamma_{-}}) \subseteq \ker(\mathbf{B}_{-})$, therefore,
\begin{equation}
    \mathbf{S}_{-}^T\mathbf{B}_{-}\boldsymbol{\mathcal{L}}_{\Gamma_{-}}^\dagger \mathbf{B}_{-}^T\mathbf{S}_{-}  = \begin{bmatrix}
        \mathbf{B}^*_{1} & \mathbf{B}^*_{2} & \mathbf{0}_{d \times (m-1)}
    \end{bmatrix}\mathbf{B}_{-}^T\mathbf{S}_{-}. \label{supp:eq:eq20}
\end{equation}
Now, combining Eq.~(\ref{supp:eq:STB}) and Eq.~(\ref{supp:eq:SmTBm1}) we can write 
\begin{align}
    \mathbf{S}^T\mathbf{B}\boldsymbol{\mathcal{P}}_{0}^T &= \begin{bmatrix}
    \mathbf{B}^*_{1}\boldsymbol{\mathcal{D}}_1 & \mathbf{B}^*_{2}(\boldsymbol{\mathcal{D}}_2 + \mathbf{I}_{n_2}) & -(\mathbf{B}^*_{1}\boldsymbol{\mathcal{K}}_1+\mathbf{B}^*_{2}\boldsymbol{\mathcal{K}}_2)  & \mathbf{B}^*_{3} & -(\mathbf{B}^*_{2}\mathbf{1}_{n_2}+\mathbf{B}^*_{3}\mathbf{1}_{n_3})
\end{bmatrix}\\
    &= \begin{bmatrix}
    \mathbf{S}_{-}^T\mathbf{B}_{-} & \mathbf{0}_{d \times n_3} & \mathbf{0}_{d}
\end{bmatrix} + \begin{bmatrix}
    \mathbf{0}_{d \times n_1} & \mathbf{B}^*_{2} &  \mathbf{0}_{d \times (m-1)} &  \mathbf{B}^*_{3} & -(\mathbf{B}^*_{2}\mathbf{1}_{n_2}+\mathbf{B}^*_{3}\mathbf{1}_{n_3})
\end{bmatrix}.
\end{align}
Finally, due to Eq.~(\ref{supp:eq:STBL_GammaBTS}) and $\mathbf{1}_{n+m} \in \ker(\mathbf{S}^T\mathbf{B}\boldsymbol{\mathcal{P}}_{0}^T)$, the above equation and Eq.~(\ref{supp:eq:eq20}),
\begin{align}
    \mathbf{S}^T\mathbf{B}\boldsymbol{\mathcal{L}}_{\Gamma}^\dagger\mathbf{B}^T\mathbf{S} &= (\mathbf{S}^T\mathbf{B}\boldsymbol{\mathcal{L}}_{\Gamma}^\dagger\boldsymbol{\mathcal{P}}_{0}^T)(\boldsymbol{\mathcal{P}}_{0}\mathbf{B}^T\mathbf{S})\\
    &= \left(\begin{bmatrix}
    \mathbf{B}^*_{1} & \mathbf{B}^*_{2} &  \mathbf{0}_{d \times (m-1)} & \mathbf{B}^*_{3} & \mathbf{0}_{d}
\end{bmatrix} + \mathbf{t}\mathbf{1}_{n+m}^T\right)(\boldsymbol{\mathcal{P}}_{0}\mathbf{B}^T\mathbf{S})\\
&=\begin{bmatrix}
        \mathbf{B}^*_{1} & \mathbf{B}^*_{2} &  \mathbf{0}_{d \times (m-1)} & \mathbf{B}^*_{3} & \mathbf{0}_{d}
    \end{bmatrix}\begin{bmatrix}
        \mathbf{S}_{-}^T\mathbf{B}_{-} & \mathbf{0}_{d \times n_3} & \mathbf{0}_{d}
    \end{bmatrix}^T + \mathbf{B}^*_{2}\mathbf{B}^{*^T}_{2} + \mathbf{B}^*_{3}\mathbf{B}^{*^T}_{3}\\
    &=\mathbf{S}_{-}^T\mathbf{B}_{-}\boldsymbol{\mathcal{L}}_{\Gamma_{-}}^\dagger \mathbf{B}_{-}^T\mathbf{S}_{-} + \mathbf{B}^*_{2}\mathbf{B}^{*^T}_{2} + \mathbf{B}^*_{3}\mathbf{B}^{*^T}_{3},
\end{align}
\hfill $\blacksquare$

\noindent \underline{\textbf{Part 3}}. Now let $\mathbf{L}(\mathbf{S})$ and $\mathbf{L}_{-}(\mathbf{S}_{-})$ be the matrices, as described in Eq.~(\ref{eq:L_of_S}) for the two graphs $\Gamma$ and $\Gamma_{-}$ and the corresponding views. Let $\boldsymbol{\Omega} = [\boldsymbol{\Omega}_i]_1^m$ be a certificate of $\mathbf{L}(\mathbf{S})$. By Remark~\ref{rmk:C_hat_L_structure}, $\boldsymbol{\Omega}^{'} = \boldsymbol{\Omega} - [\boldsymbol{\Omega}_m]_1^m$ is also a certificate of $\mathbf{L}(\mathbf{S})$ and in particular $\boldsymbol{\Omega}^{'}_m = 0$. Define $\boldsymbol{\Omega}_{-} = [\boldsymbol{\Omega}^{'}_i]_1^{m-1}$. We are going to show that $\Tr(\boldsymbol{\Omega}_{-}^T\mathbf{L}_{-}(\mathbf{S}_{-})\boldsymbol{\Omega}_{-}) = 0$ i.e. $\boldsymbol{\Omega}_{-}$ is a certificate of $\mathbf{L}_{-}(\mathbf{S}_{-})$. Then using Remark~\ref{rmk:C_hat_L_structure}, it follows that $[\boldsymbol{\Omega}_i]_1^{m-1}$ (which equals $\boldsymbol{\Omega}_{-} + [\boldsymbol{\Omega}_m]_1^{m-1}$) is a certificate of $\mathbf{L}_{-}(\mathbf{S}_{-})$.
First, we note that
\begin{align}
    \mathbf{B}(\mathbf{S})\boldsymbol{\mathcal{P}}_{0}^T &= \begin{bmatrix}
        \mathbf{B}_{-}(\mathbf{S}_{-}) &  \mathbf{0}_{(m-1)d \times n_3} & \mathbf{0}_{(m-1)d}\\
        \begin{bmatrix}\mathbf{0}_{d \times n_1} & \mathbf{B}^*_{2} & \mathbf{0}_{d \times (m-1)} \end{bmatrix} &\mathbf{B}^*_{3} & -(\mathbf{B}^*_{2}\mathbf{1}_{n_2}+\mathbf{B}^*_{3}\mathbf{1}_{n_3})
    \end{bmatrix}\\
    \mathbf{D}(\mathbf{S}) &= \blockdiag(\mathbf{D}_{-}(\mathbf{S}_{-}), \mathbf{B}^*_{2}\mathbf{B}^{*^T}_{2} + \mathbf{B}^*_{3}\mathbf{B}^{*^T}_{3}).
\end{align}
Since $\boldsymbol{\Omega}^{'}$ is a certificate of $\mathbf{L}(\mathbf{S})$, $\Tr(\boldsymbol{\Omega}^{'^T}\mathbf{L}(\mathbf{S})\boldsymbol{\Omega}^{'}) = 0$. Then, using the definition of $\mathbf{L}(\mathbf{S})$, the above equations, the fact that $\boldsymbol{\Omega}^{'}_m = 0$, and Eq.~(\ref{supp:eq:pinvA11}), we obtain
\begin{align}
    0 &= \Tr(\boldsymbol{\Omega}^{'^T}(\mathbf{D}(\mathbf{S}) - \mathbf{B}(\mathbf{S})\boldsymbol{\mathcal{P}}_{0}^T\boldsymbol{\mathcal{P}}_{0}\boldsymbol{\mathcal{L}}_{\Gamma}^\dagger\boldsymbol{\mathcal{P}}_{0}^T\boldsymbol{\mathcal{P}}_{0} \mathbf{B}(\mathbf{S})^T)\boldsymbol{\Omega}^{'})\\
    &= \Tr(\boldsymbol{\Omega}_{-}^T(\mathbf{D}_{-}(\mathbf{S}) - \mathbf{B}_{-}(\mathbf{S}_{-}) (\boldsymbol{\mathcal{P}}_{0}\boldsymbol{\mathcal{L}}_{\Gamma}^\dagger\boldsymbol{\mathcal{P}}_{0}^T)_{11} \mathbf{B}_{-}(\mathbf{S}_{-})^T)\boldsymbol{\Omega}_{-})\\
    &= \Tr(\boldsymbol{\Omega}_{-}^T(\mathbf{D}_{-}(\mathbf{S}) - \mathbf{B}_{-}(\mathbf{S}_{-})(\boldsymbol{\mathcal{L}}_{\Gamma_{-}}+\boldsymbol{\mathcal{D}}_0)^\dagger \mathbf{B}_{-}(\mathbf{S}_{-})^T)\boldsymbol{\Omega}_{-}).
\end{align}

From Proposition~\ref{prop:noiseless_setting1} and Eq.~(\ref{eq:C_of_S}), 
$$\mathbf{L}_{-}(\mathbf{S}_{-}) = \mathbf{D}_{-}(\mathbf{S}_{-}) - \mathbf{B}_{-}(\mathbf{S}_{-})\boldsymbol{\mathcal{L}}_{\Gamma_{-}}^\dagger \mathbf{B}_{-}(\mathbf{S}_{-})^T,$$ thus 
$$\Tr(\boldsymbol{\Omega}_{-}^T\mathbf{L}_{-}(\mathbf{S}_{-})\boldsymbol{\Omega}_{-}) - \Tr(\boldsymbol{\Omega}_{-}^T(\mathbf{B}_{-}(\mathbf{S}_{-}) ((\boldsymbol{\mathcal{L}}_{\Gamma_{-}}+\boldsymbol{\mathcal{D}}_0)^\dagger - \boldsymbol{\mathcal{L}}_{\Gamma_{-}}^\dagger) \mathbf{B}_{-}(\mathbf{S}_{-})^T )\boldsymbol{\Omega}_{-}) = 0.$$
Since $\mathbf{L}_{-}(\mathbf{S}_{-}) \succeq 0$, to show that $\Tr(\boldsymbol{\Omega}_{-}^T\mathbf{L}_{-}(\mathbf{S}_{-})\boldsymbol{\Omega}_{-}) = 0$, it suffices to show that
\begin{prop} $\mathbf{B}_{-}(\mathbf{S}_{-}) ((\boldsymbol{\mathcal{L}}_{\Gamma_{-}}+\boldsymbol{\mathcal{D}}_0)^\dagger - \boldsymbol{\mathcal{L}}_{\Gamma_{-}}^\dagger) \mathbf{B}_{-}(\mathbf{S}_{-})^T \preceq 0$.
\end{prop}
\textit{Proof}. Since $\boldsymbol{\mathcal{L}}_{\Gamma_{-}} \succeq 0$, consider $\boldsymbol{\mathcal{L}}_{\Gamma_{-}} = \mathbf{U}\boldsymbol{\Lambda}\mathbf{U}^T = \begin{bsmallmatrix}
        \mathbf{U}_1 & \mathbf{U}_2
    \end{bsmallmatrix}\begin{bsmallmatrix}
        \mathbf{\Lambda}_1 & \\
        & \mathbf{0}
    \end{bsmallmatrix}\begin{bsmallmatrix}
        \mathbf{U}_1^T\\
        \mathbf{U}_2^T
    \end{bsmallmatrix}$
where the $1 \leq n' \leq n_2$ columns of $\mathbf{U}_2 \in \mathbb{R}^{(n_1+n_2+m-1) \times n'}$ form an orthogonal basis of the $\ker (\boldsymbol{\mathcal{L}}_{\Gamma_{-}})$ (see proof of Proposition~\ref{supp:prop:LplusD0}). Also, $\mathbf{\Lambda}_1 \succ 0$. From Proposition~\ref{prop:kerB} and Eq.~(\ref{eq:BofS}), $\mathbf{B}_{-}(\mathbf{S}_{-})\mathbf{U}_2 = 0$. Thus, $\mathbf{B}_{-}(\mathbf{S}_{-})\mathbf{U} = \begin{bmatrix}
    \mathbf{B}_{-}(\mathbf{S}_{-})\mathbf{U}_1 & 0
\end{bmatrix}$. Then note that
\begin{equation}
    ((\boldsymbol{\mathcal{L}}_{\Gamma_{-}}+\boldsymbol{\mathcal{D}}_0)^\dagger - \boldsymbol{\mathcal{L}}_{\Gamma_{-}}^\dagger) = \mathbf{U} \left\{\left(\begin{bmatrix}
        \mathbf{\Lambda}_1 & \\
        & \mathbf{0}
    \end{bmatrix} + \mathbf{U}^T\boldsymbol{\mathcal{D}}_0\mathbf{U}\right)^\dagger-\begin{bmatrix}
        \mathbf{\Lambda}_1^{-1} & \\
        & \mathbf{0}
    \end{bmatrix}\right\}\mathbf{U}^T. \label{supp:eq:eq30}
\end{equation}
Using \citeb[Eq.~(10, 11, 17, 19)]{kovanic1979pseudoinverse} and simple calculations, we obtain
$$\left(\blockdiag(\mathbf{\Lambda}_1, \mathbf{0}) + \mathbf{U}^T\boldsymbol{\mathcal{D}}_0\mathbf{U}\right)^\dagger = (\blockdiag(\mathbf{\Lambda}_1^{-1}, \mathbf{0}) + \blockdiag(\mathbf{W}_1, \mathbf{W}_2))$$
where 
$$\mathbf{W}_1 = -\mathbf{\Lambda}_1^{-1}\mathbf{U}_1^T\boldsymbol{\mathcal{D}}_0(\mathbf{I} + \boldsymbol{\mathcal{D}}_0\mathbf{U}_1\mathbf{\Lambda}_1^{-1}\mathbf{U}_1^T\boldsymbol{\mathcal{D}}_0)^{-1}\boldsymbol{\mathcal{D}}_0\mathbf{U}_1\mathbf{\Lambda}_1^{-1}$$
and $\mathbf{W}_2 = (\mathbf{U}_2^T\boldsymbol{\mathcal{D}}_0\mathbf{U}_2)^\dagger$.
In particular $\mathbf{W}_1 \preceq 0$ and $\mathbf{W}_2 \succeq 0$.
Combining above with Eq.~(\ref{supp:eq:eq30}), we deduce that 
$$\mathbf{B}_{-}(\mathbf{S}_{-}) ((\boldsymbol{\mathcal{L}}_{\Gamma_{-}}+\boldsymbol{\mathcal{D}}_0)^\dagger - \boldsymbol{\mathcal{L}}_{\Gamma_{-}}^\dagger) \mathbf{B}_{-}(\mathbf{S}_{-})^T = \mathbf{B}_{-}(\mathbf{S}_{-})\mathbf{U}_1\mathbf{W}_1\mathbf{U}_1^T\mathbf{B}_{-}(\mathbf{S}_{-})^T \preceq 0.$$
\hfill $\blacksquare$

\noindent We conclude that $\Tr(\boldsymbol{\Omega}_{-}^T\mathbf{L}_{-}(\mathbf{S}_{-})\boldsymbol{\Omega}_{-}) = 0$ and thus $\boldsymbol{\Omega}_{-}$ is a certificate of $\mathbf{L}_{-}(\mathbf{S}_{-})$.

\proofof{Proposition~\ref{prop:same_conn_comp_non_deg}}
It suffices to show that if the $i$th and $j$th vertices are adjacent in $\mathbb{G}$ then $\boldsymbol{\Omega}_i = \boldsymbol{\Omega}_j$. For $m=2$, the result is a direct consequence of Theorem~\ref{thm:nec_suff_cond_loc_rigid_two_views}. Suppose the result holds for $m-1$ views for some $m > 2$. If there are no edges in $\mathbb{G}$, then the result holds trivially for $m$ views. Suppose $i$th and $j$th vertices are adjacent in $\mathbb{G}$. Let $r \in [1,m] \setminus \{i,j\}$. We remove the $r$th view and the points which lie exclusively in it. Then by Lemma~\ref{lem:subproblem_cert}, $[\boldsymbol{\Omega}_k]_{k \in [1,m] \setminus r}$ is a certificate of $\mathbf{L}_{-r}(\mathbf{S}_{-r})$. Now, construct $\mathbb{G}_{-r}$ (in the same way as $\mathbb{G}$) and note that the $i$th and $j$th vertices are still adjacent. Thus, by the induction hypothesis, $\boldsymbol{\Omega}_i = \boldsymbol{\Omega}_j$.

\proofof{Theorem~\ref{thm:G_star_1}}
The result holds for two views (see Theorem~\ref{thm:nec_suff_cond_loc_rigid_two_views}). Suppose it holds for $m-1$ views for some $m > 2$. Let $\boldsymbol{\Omega}$ be a certificate of $\mathbf{L}(\mathbf{S})$. We need to show that $\boldsymbol{\Omega}$ is trivial. Since $|\mathbb{G}^*(\mathbf{S})|=1$, $\mathbb{G}$ must have a connected component with at least two views. Pick one such component and note that there exist a view in it such that removing it will not disconnect the component. Let it be the $i$th view. Consider removing the $i$th view and the points which lie exclusively in it. For the new set of views we still have $|\mathbb{G}^*(\mathbf{S}_{-i})|=1$ (where $\mathbb{G}^*_{-i}(\mathbf{S}_{-i})$ is constructed in the same manner as $\mathbb{G}^*(\mathbf{S})$). By Lemma~\ref{lem:subproblem_cert} and the induction hypothesis, we conclude that $[\boldsymbol{\Omega}_j]_{j \in [1,m]\setminus \{i\}}$ must be trivial. By Proposition~\ref{prop:same_conn_comp_non_deg} we conclude that $\boldsymbol{\Omega}$ is trivial.

\proofof{Theorem~\ref{thm:nec_cond_glob_rigid_views}}
The calculations up to Eq.~(\ref{supp:eq:eq11}) in the proof of Theorem~\ref{thm:nec_cond_loc_rigid_of_views} are reused here.
Consider a partition of $[1,m]$ into two non-empty subsets $A$ and $B$. Suppose the rank of $\overline{\mathbf{B}(\mathbf{S})}_{A, B}$ is at most $d-1$. As in the proof of Theorem~\ref{thm:nec_cond_loc_rigid_of_views}, WLOG assume that $\mathbf{B}(\mathbf{S})_{A, B}\mathbf{1}_{n'} = 0$ (where $n'$ is as in Definition~\ref{def:BSAcapB}). Then the rank of $\mathbf{B}(\mathbf{S})_{A, B}$ is at most $d-1$. We are going to construct another perfect alignment $\mathbf{S}'$ such that $\pi(\mathbf{S}) \neq \pi(\mathbf{S}')$, thus concluding that $\mathbf{S}$ is not unique.

Let $\mathbf{V}_1, \mathbf{V}_2 \in \mathbb{O}(d)$ and $\mathbf{\Sigma}$ be the diagonal matrix containing the singular values of $\mathbf{B}(\mathbf{S})_{A, B}$ such that $\mathbf{B}(\mathbf{S})_{A, B} = \mathbf{V}_1\mathbf{\Sigma}\mathbf{V}_2^T$. Since rank of $\mathbf{B}(\mathbf{S})_{A, B} \leq d-1$, there exist $\mathbf{U} \in \mathbb{O}(d)$ such that $\mathbf{U} \neq \mathbf{I}_d$ and $\mathbf{\Sigma} = \mathbf{U}\mathbf{\Sigma}$. Define $\mathbf{Q} = \mathbf{V}_1\mathbf{U}^T\mathbf{V}_1^T$. Then $\mathbf{Q} \neq \mathbf{I}_d$ and 
\begin{equation}
    \mathbf{Q}^T\mathbf{B}(\mathbf{S})_{A, B} = \mathbf{B}(\mathbf{S})_{A, B}. \label{supp:eq:QTBSAB}
\end{equation}
Define $\mathbf{S}' \in \mathbb{O}(d)^m$ such that $\mathbf{S}'_i = \mathbf{S}_i\mathbf{Q}$ for all $i \in A$ and $\mathbf{S}'_j = \mathbf{S}_j$ for all $j \in B$. Clearly, $\mathbf{S}' \neq \mathbf{S}$. We will show that $\mathbf{S}'$ is another perfect alignment. It is easy to see that
$$\Tr(\mathbf{S}'^T\mathbf{C}\mathbf{S}') = \textstyle\sum_{\substack{i \in A, j \in A\\i \in B, j \in B}}\Tr(\mathbf{S}_i^T\mathbf{C}_{ij}\mathbf{S}_j) + 2 \textstyle\sum_{i \in A, j \in B}\Tr(\mathbf{Q}^T\mathbf{S}_i^T\mathbf{C}_{ij}\mathbf{S}_j).$$
Since, for $i \in A$ and $j \in B$, $\mathbf{C}_{ij} = \mathbf{B}_i \boldsymbol{\mathcal{L}}_{\Gamma}^\dagger \mathbf{B}_j^T$, it suffices to show that
$$\Tr\left(\left(\textstyle\sum_{i \in A}\mathbf{B}(\mathbf{S})_i\right)\boldsymbol{\mathcal{L}}_{\Gamma}^\dagger \left(\textstyle\sum_{j \in B}\mathbf{B}(\mathbf{S})_j^T\right)\right) = \Tr\left(\mathbf{Q}^T\left(\textstyle\sum_{i \in A}\mathbf{B}(\mathbf{S})_i\right)\boldsymbol{\mathcal{L}}_{\Gamma}^\dagger \left(\textstyle\sum_{j \in B}\mathbf{B}(\mathbf{S})_j^T\right)\right).$$
Define $\mathbf{B}_A$, $\mathbf{B}_B$,  $\mathbf{B}_{A\setminus B}$, $\mathbf{B}_{A,B}$, $\mathbf{B}_{B\setminus A}$ and $\mathbf{B}_{B,A}$ as in the proof of Theorem~\ref{thm:nec_cond_loc_rigid_of_views}, then it suffices to show that 
$\Tr(\mathbf{Q}^T \mathbf{B}_A\boldsymbol{\mathcal{L}}_{\Gamma}^\dagger \mathbf{B}_B^T) = \Tr( \mathbf{B}_A\boldsymbol{\mathcal{L}}_{\Gamma}^\dagger \mathbf{B}_B^T).$
Using Eq.~(\ref{supp:eq:QTBSAB}, \ref{supp:eq:eq8}, \ref{supp:eq:eq9}) we obtain $\mathbf{Q}^T \mathbf{B}_{A,B} = \mathbf{B}_{A,B}$ and this combined with Eq.~(\ref{supp:eq:eq11}), yields
\begin{align}
    \Tr(\mathbf{Q}^T \mathbf{B}_A\boldsymbol{\mathcal{L}}_{\Gamma}^\dagger \mathbf{B}_B^T) &= \Tr(\mathbf{Q}^T (\mathbf{B}_{A \setminus B}+\mathbf{B}_{A,B})\boldsymbol{\mathcal{L}}_{\Gamma}^\dagger (\mathbf{B}_{B\setminus A} + \mathbf{B}_{B,A})^T ) \\
    &= \Tr(\mathbf{Q}^T \mathbf{B}_{A,B}\boldsymbol{\mathcal{L}}_{\Gamma}^\dagger\mathbf{B}_{B,A}^T) = \Tr( \mathbf{B}_{A,B}\boldsymbol{\mathcal{L}}_{\Gamma}^\dagger\mathbf{B}_{B,A}^T)\\
    &=\Tr((\mathbf{B}_{A \setminus B}+\mathbf{B}_{A,B})\boldsymbol{\mathcal{L}}_{\Gamma}^\dagger (\mathbf{B}_{B\setminus A} + \mathbf{B}_{B,A})^T )\\
    &= \Tr(\mathbf{B}_A\boldsymbol{\mathcal{L}}_{\Gamma}^\dagger \mathbf{B}_B^T)
\end{align}

\proofof{Proposition~\ref{prop:same_conn_comp_uniq}} By replacing $\mathbb{G}$ with $\overline{\mathbb{G}}$ and Theorem~\ref{thm:nec_suff_cond_loc_rigid_two_views} with \ref{thm:nec_suff_cond_glob_rigid_two_views}, the inductive proof is the same as of Proposition~\ref{prop:same_conn_comp_non_deg}.

\proofof{Theorem~\ref{thm:overline_G_star_1}} By replacing $\mathbb{G}$ with $\overline{\mathbb{G}}$, Theorem~\ref{thm:nec_suff_cond_loc_rigid_two_views} with \ref{thm:nec_suff_cond_glob_rigid_two_views} and Proposition~\ref{prop:same_conn_comp_non_deg} with \ref{prop:same_conn_comp_uniq}, the inductive proof is the same as of Theorem~\ref{thm:G_star_1}.


\proofof{Lemma~\ref{lem:retraction}}
For $\mathbf{S} \in \pi^{-1}(\widetilde{\mathbf{S}})$ and $\mathbf{Z} = [\mathbf{S}_i\boldsymbol{\Omega}_i]_1^m \in T_{\mathbf{S}}\mathbb{O}(d)^m$, the horizontal lift of $\widetilde{\mathbf{Z}} = [\widetilde{\mathbf{S}}_i\widetilde{\boldsymbol{\Omega}}_i] \in T_{\widetilde{\mathbf{S}}}\mathbb{O}(d)^m/_{\sim}$,
$$\pi(R_{\EXP}(\mathbf{S}, \mathbf{Z})) = [\mathbf{S}_i\exp(\boldsymbol{\Omega}_i)(\mathbf{S}_1\exp(\boldsymbol{\Omega}_1))^T]_1^m.$$
It suffices to show that $\mathbf{S}_{i+1}\exp(\boldsymbol{\Omega}_{i+1})(\mathbf{S}_1\exp(\boldsymbol{\Omega}_1))^T$ depends only on $\widetilde{\mathbf{S}}$ and $\widetilde{\boldsymbol{\Omega}}$ for all $i \in [1,m-1]$. Using Proposition~\ref{prop:hlift_char} and expanding the expression, we obtain 
$$\mathbf{S}_{i+1}\exp(\boldsymbol{\Omega}_{i+1})(\mathbf{S}_1\exp(\boldsymbol{\Omega}_1))^T = \mathbf{S}_{i+1}\exp(\mathbf{S}_1^T \widetilde{\boldsymbol{\Omega}}_i\mathbf{S}_1 + \boldsymbol{\Omega}_1) (\mathbf{S}_1\exp(\boldsymbol{\Omega}_1))^T.$$
Since $\mathbf{S}_1 \in \mathbb{O}(d)$, we have $\exp(\mathbf{S}_1^T \widetilde{\boldsymbol{\Omega}}_i\mathbf{S}_1) = \mathbf{S}_1^T\exp(\widetilde{\boldsymbol{\Omega}}_i)\mathbf{S}_1$. Also, the following identites hold: $\exp(\boldsymbol{\Omega}_1)^T = \exp(\boldsymbol{\Omega}_1^T) = \exp(-\boldsymbol{\Omega}_1)$ and $\exp(\mathbf{A}_1+\mathbf{A}_1) = \exp(\mathbf{A}_1)\exp(\mathbf{A}_2)$. By substituting back into the expression, we obtain 
$$\mathbf{S}_{i+1}\exp(\mathbf{S}_1^T \widetilde{\boldsymbol{\Omega}}_i\mathbf{S}_1 + \boldsymbol{\Omega}_1) (\mathbf{S}_1\exp(\boldsymbol{\Omega}_1))^T = \mathbf{S}_{i+1}\mathbf{S}_1^T\exp(\widetilde{\boldsymbol{\Omega}}_i) = \widetilde{\mathbf{S}}_i\exp(\widetilde{\boldsymbol{\Omega}}_i).$$

\proofof{Proposition~\ref{prop:liu_pf}}
Since $\mathbf{Z}_i = \mathbf{S}_i\boldsymbol{\Omega}_i$ where $\boldsymbol{\Omega}_i \in \Skew(d)$,
$$
    \left\|\mathbf{S}_i\exp(\mathbf{S}_i^T\mathbf{Z}_i) - (\mathbf{S}_i + \mathbf{Z}_i)\right\|_F = \left\|\exp(\boldsymbol{\Omega}_i) - (\mathbf{I}_d + \boldsymbol{\Omega}_i)\right\|_F \leq \left(\sum_{2}^{\infty}\frac{1}{k !}\right)\left\|\boldsymbol{\Omega}_i\right\|_F^2 = (e-1)\left\|\boldsymbol{\Omega}_i\right\|_F^2
$$
where the inequality follows from the triangle inequality and the fact that $\left\|\boldsymbol{\Omega}_i\right\|_F = \left\|\mathbf{Z}_i\right\|_F \leq 1$.

\proofof{Proposition~\ref{prop:second_order_boundedness_of_Rtilde}}
The proof of part (a) follows from Proposition~\ref{prop:liu_pf}, Proposition~\ref{prop:hlift_frob_ineq} and 
$$\left\|R_\EXP(\mathbf{S}, \boldsymbol{\xi}) - (\mathbf{S} + \boldsymbol{\xi})\right\|_F \leq \sum_1^m \left\|\mathbf{S}_i\exp(\mathbf{S}_i^T\boldsymbol{\xi}_i) - (\mathbf{S}_i + \boldsymbol{\xi}_i)\right\|_F.$$
For part (b), let $\mathbf{S} \in \pi^{-1}(\widetilde{\mathbf{S}})$ and $\mathbf{Z} = [\mathbf{S}_i\boldsymbol{\Omega}_i]_1^m \in T_{\mathbf{S}}\mathbb{O}(d)^m$ be the horizontal lift of $\widetilde{\mathbf{Z}} = [\widetilde{\mathbf{S}}_i\widetilde{\boldsymbol{\Omega}}_i]_1^m$ at $\mathbf{S}$. Note that $\widetilde{\mathbf{S}}_i = \mathbf{S}_{i+1}\mathbf{S}_1^T$ (Eq.~(\ref{eq:pi_inv_wtS})) and $\boldsymbol{\Omega}_{i+1} -  \boldsymbol{\Omega}_{1}= \mathbf{S}_1^T\widetilde{\boldsymbol{\Omega}}_{i}\mathbf{S}_1$ (Eq.~(\ref{eq:hlifti})). Then, using the identity $\exp(\boldsymbol{\Omega}_{i+1})\exp(\boldsymbol{\Omega}_{1})^T = \exp(\boldsymbol{\Omega}_{i+1}-\boldsymbol{\Omega}_{1})$, the inequality $\left\|\boldsymbol{\Omega}_{i+1}-\boldsymbol{\Omega}_{1}\right\|_F \leq \left\|\boldsymbol{\Omega}_{i+1}\right\|_F + \left\|\boldsymbol{\Omega}_{1}\right\|_F \leq 1$ and Proposition~\ref{prop:liu_pf}, we obtain
\begin{align}
    \left\|\widetilde{R}_\EXP(\widetilde{\mathbf{S}}, \widetilde{\mathbf{Z}})_i - (\widetilde{\mathbf{S}}_i + \widetilde{\mathbf{Z}}_i)\right\|_F &= \left\|\mathbf{S}_{i+1}(\exp(\boldsymbol{\Omega}_{i+1})\exp(\boldsymbol{\Omega}_{1})^T - (\mathbf{I}_d +\boldsymbol{\Omega}_{i+1} -  \boldsymbol{\Omega}_{1}))\mathbf{S}_{1}^T\right\|_F\\
    &\leq (e-1)\left\|\boldsymbol{\Omega}_{i+1} -  \boldsymbol{\Omega}_{1}\right\|_F^2 = (e-1)\left\|\widetilde{\boldsymbol{\Omega}}_{i}\right\|_F^2 = (e-1)\left\|\widetilde{\mathbf{Z}}_i\right\|_F^2.
\end{align}
Then the result follows from $\left\|\widetilde{R}_\EXP(\widetilde{\mathbf{S}}, \widetilde{\mathbf{Z}}) - (\widetilde{\mathbf{S}} + \widetilde{\mathbf{Z}})\right\|_F \leq \sum_1^m \left\|\widetilde{R}_\EXP(\widetilde{\mathbf{S}}, \widetilde{\mathbf{Z}})_i - (\widetilde{\mathbf{S}}_i + \widetilde{\mathbf{Z}}_i)\right\|_F$.

\proofof{Proposition~\ref{prop:alpha_grad}}
\revadd{The proof is similar to the one in \citea{liu2019quadratic}. Due to \textbf{(A2)} in Section~\ref{subsec:loc_lin_conv}, WLOG we assume that $\grad F(\mathbf{S}^k) \neq 0$ for all $k \geq 0$. Then, from the proof of Proposition~\ref{prop:gradFS}, since 
$$\grad F(\mathbf{S})_i = 0.5\left(\nabla F(\mathbf{S})_i - \mathbf{S}_i\nabla F(\mathbf{S})_i^T\mathbf{S}_i\right),$$
we obtain 
\begin{align}
g(\nabla F(\mathbf{S}^k)_i, \grad F(\mathbf{S}^k)_i) &= 0.5\left(\left\|\nabla F(\mathbf{S}^k)_i\right\|_F^2 -g(\nabla F(\mathbf{S}^k)_i,\mathbf{S}^k_i\nabla F(\mathbf{S}^k)_i^T\mathbf{S}^k_i)\right)\\
&= 0.25 \left\|\nabla F(\mathbf{S}^k)_i - \mathbf{S}^k_i\nabla F(\mathbf{S}^k)_i^T\mathbf{S}^k_i\right\|_F^2\\
&= \left\|\grad F(\mathbf{S}^k)_i\right\|_F^2.
\end{align}
Thus, 
$$g(\nabla F(\mathbf{S}^k), \grad F(\mathbf{S}^k)) = \textstyle\sum_1^m g(\nabla F(\mathbf{S}^k)_i, \grad F(\mathbf{S}^k)_i) = \left\|\grad F(\mathbf{S}^k)\right\|_F^2.$$ 
This together with Algorithm~\ref{algo:rgd}, and Eq.~(\ref{eq:armijo_step}), implies
\begin{equation}
    F(\mathbf{S}^{k+1})-F(\mathbf{S}^{k}) \leq -\gamma \alpha_k \left\|\grad F(\mathbf{S}^k)\right\|_F^2 \label{supp:eq:eq1}
\end{equation}
Since $\mathbb{O}(d)^m$ is compact and $F$ is analytic, thus $F$ is bounded. Since $\alpha_k \in (0,1]$ for all $k \geq 0$, 
$$\textstyle\sum_0^\infty \alpha_k^2 \left\|\grad F(\mathbf{S}^k)\right\|_F^2 \leq \textstyle\sum_0^\infty \alpha_k \left\|\grad F(\mathbf{S}^k)\right\|_F^2 \leq \frac{1}{\gamma}(F(\mathbf{S}^0) - \lim F(\mathbf{S}^k)) < \infty.$$
Thus, $\lim \alpha_k \left\|\grad F(\mathbf{S}^k)\right\|_F = 0$ and from Proposition~\ref{prop:hlift_frob_ineq}, we obtain $\lim \alpha_k \left\|\grad \widetilde{F}(\widetilde{\mathbf{S}}^k)\right\|_F = 0$.}

\proofof{(A1) in Section~\ref{subsec:loc_lin_conv}}
From Proposition~\ref{prop:alpha_grad}, there exists $k_1 \geq 0$ such that $\alpha_k\left\|\grad \widetilde{F}(\widetilde{\mathbf{S}}^k)\right\|_F \leq 1/2$ for all $k \geq k_1$. Then note that 
$$\left\|\widetilde{\mathbf{S}}^{k+1} - \widetilde{\mathbf{S}}^k\right\|_F = \left\|\widetilde{R}_\EXP(\widetilde{\mathbf{S}}^k, -\alpha_k\grad \widetilde{F}(\widetilde{\mathbf{S}}^k)) - \widetilde{\mathbf{S}}^k\right\|_F.$$
From Proposition~\ref{prop:second_order_boundedness_of_Rtilde}, for all $k \geq k_1$,
\begin{align}
    \left\|\widetilde{R}_\EXP(\widetilde{\mathbf{S}}^k, -\alpha_k\grad \widetilde{F}(\widetilde{\mathbf{S}}^k)) - \widetilde{\mathbf{S}}^k\right\|_F &\leq (e-1) \alpha_k^2 \left\|\grad \widetilde{F}(\widetilde{\mathbf{S}}^k)\right\|_F^2 + \alpha_k\left\|\grad \widetilde{F}(\widetilde{\mathbf{S}}^k)\right\|_F\\
    &\leq \frac{1}{2}(e+1)\alpha_k \left\|\grad \widetilde{F}(\widetilde{\mathbf{S}}^k)\right\|_F
\end{align}
Finally, using Eq~(\ref{supp:eq:eq1}) and Proposition~\ref{prop:hlift_frob_ineq}, for all $k \geq k_1$, 
$$\widetilde{F}(\widetilde{\mathbf{S}}^{k+1})-\widetilde{F}(\widetilde{\mathbf{S}}^{k}) \leq -2\gamma(e+1)^{-1}(m+1)^{-1/2} \left\|\grad \widetilde{F}(\widetilde{\mathbf{S}}^k)\right\|_F \cdot \left\|\widetilde{\mathbf{S}}^{k+1}-\widetilde{\mathbf{S}}^k\right\|_F.$$
Thus, (\textbf{A1}) holds for $\kappa_0 = 2\gamma(e+1)^{-1}(m+1)^{-1/2}$. The result follows.

\proofof{(A3) in Section~\ref{subsec:loc_lin_conv}}
Since $\nabla F(\mathbf{S}) = 2\mathbf{C}\mathbf{S}$ is Lipschitz with parameter $L_F \leq 2\left\|\mathbf{C}\right\|_F$, the proof of \textbf{(A3)} is same as in \citeb[Pg. 235]{liu2019quadratic} (alternatively \citeb[Theorem 2.10]{schneider2015convergence}). For the sake of completeness, we present an adaptation of their proof to our setting.
WLOG we assume that $\alpha_k \left\|\grad \widetilde{F}(\widetilde{\mathbf{S}}^{k}))\right\|_F \neq 0$ for all $k \geq 0$. Then, $\widetilde{\mathbf{S}}^{k+1} - \widetilde{\mathbf{S}}^{k} = \widetilde{R}_\EXP(\widetilde{\mathbf{S}}^{k}, -\alpha_k \grad \widetilde{F}(\widetilde{\mathbf{S}}^{k})) - \widetilde{\mathbf{S}}^{k}$ and we have
\begin{align}
    \left\|\widetilde{\mathbf{S}}^{k+1}-\widetilde{\mathbf{S}}^{k}\right\|_F &\geq \alpha_k \left\|\grad \widetilde{F}(\widetilde{\mathbf{S}}^{k}))\right\|_F - \left\|\widetilde{R}_\EXP(\widetilde{\mathbf{S}}^{k}, -\alpha_k \grad \widetilde{F}(\widetilde{\mathbf{S}}^{k})) - (\widetilde{\mathbf{S}}^{k} - \grad \widetilde{F}(\widetilde{\mathbf{S}}^{k})))\right\|_F\\
    \left\|\widetilde{\mathbf{S}}^{k+1}-\widetilde{\mathbf{S}}^{k}\right\|_F &\leq \alpha_k \left\|\grad \widetilde{F}(\widetilde{\mathbf{S}}^{k}))\right\|_F + \left\|\widetilde{R}_\EXP(\widetilde{\mathbf{S}}^{k}, -\alpha_k \grad \widetilde{F}(\widetilde{\mathbf{S}}^{k})) - (\widetilde{\mathbf{S}}^{k} - \grad \widetilde{F}(\widetilde{\mathbf{S}}^{k})))\right\|_F.
\end{align}
Using Proposition~\ref{prop:second_order_boundedness_of_Rtilde}, we obtain 
$$\lim \frac{\left\|\widetilde{\mathbf{S}}^{k+1}-\widetilde{\mathbf{S}}^{k}\right\|_F}{\alpha_k\left\|\grad \widetilde{F}(\widetilde{\mathbf{S}}^{k}))\right\|_F} = 1.$$
It suffices to show that $\liminf \alpha_k > 0$. Let $\overline{\alpha}_k = \overline{\alpha}(\mathbf{S}^k) > 0$ where 
\begin{equation}
    \overline{\alpha}(\mathbf{S}) = \inf\{\alpha > 0| F(R_\EXP(\mathbf{S}, -\alpha \grad F(\mathbf{S}))) - F(\mathbf{S}) = -\gamma \alpha g(\nabla F(\mathbf{S}), \grad F(\mathbf{S}))\}.
\end{equation}
where $\overline{\alpha}(\mathbf{S})$ is well defined due to \citeb[Proposition 2.8]{schneider2015convergence} since $F$ extends to a continuously differentiable non-negative function on $\mathbb{R}^{md \times d}$ containing $\mathbb{O}(d)^m$. By Eq.~(\ref{eq:armijo_step}) and above equation, we have $\alpha_k = 1$ if $\overline{\alpha}_k \geq 1$ and $\alpha_k \geq \beta \overline{\alpha}_k$ if $\overline{\alpha}_k < 1$.
It suffices to assume that $\overline{\alpha}_k < 1$ for all $k \geq 0$ and show that $\liminf \overline{\alpha}_k > 0$. From the above equation, it follows that $\overline{\alpha}_k \left\|\grad F(\mathbf{S}^{k}))\right\|_F \leq (\alpha_k/\beta)\left\|\grad F(\mathbf{S}^{k}))\right\|_F$ which combined with Proposition~\ref{prop:alpha_grad} implies that $\lim\overline{\alpha}_k \left\|\grad F(\mathbf{S}^k)\right\|_F = 0$.
By mean value theorem and the definition of $\overline{\alpha}_k$, there exist $\zeta_k \in (0,1)$ such that 
$$\mathbf{U}^k = \zeta_k(R_\EXP(\mathbf{S}^k, -\overline{\alpha}_k \grad F(\mathbf{S}^k)) - \mathbf{S}^k)$$
satisfies
\begin{align}
    (R_\EXP(\mathbf{S}, -\overline{\alpha}_k \grad F(\mathbf{S}))- \mathbf{S}^k)^T\nabla F(\mathbf{S}^k + \mathbf{U}_k) &= F(R_\EXP(\mathbf{S}, -\overline{\alpha}_k \grad F(\mathbf{S}))) - F(\mathbf{S}^k)\\
    &= -\gamma \overline{\alpha}_k g(\nabla F(\mathbf{S}^k), \grad F(\mathbf{S}^k)). \label{eq:Uk}
\end{align}
Moreover, for sufficiently large $k \geq 0$, $\overline{\alpha}_k\left\|\grad F(\mathbf{S}^k)\right\|_F < 1$, therefore using Proposition~\ref{prop:liu_pf} and the triangle inequality,
\begin{equation}
    \left\|\mathbf{U}^k\right\|_F \leq \left\|R_\EXP(\mathbf{S}^k, -\overline{\alpha}_k \grad F(\mathbf{S}^k)) - \mathbf{S}^k\right\|_F \leq e\overline{\alpha}_k\left\|\grad F(\mathbf{S}^k)\right\|_F.
\end{equation}
Then we obtain the following set of inequalities using the above inequality, the fact that $\nabla F$ is Lipschitz continuous with parameter $L_F \leq 2 \left\|\mathbf{C}\right\|_F$, using Cauchy-Schwarz inequality, Eq~(\ref{eq:Uk}) and the triangle inequality,

\begin{align}
    \overline{\alpha}_k^2&\left\|\grad F(\mathbf{S}^k)\right\|_F^2 \geq e^{-1}\left\|\mathbf{U}_k\right\|_F \left\|R_\EXP(\mathbf{S}^k, - \alpha_k \grad F(\mathbf{S}^k)) - \mathbf{S}^k\right\|_F\\
    &\geq (eL_F)^{-1} \left\|\nabla F(\mathbf{S}^k) - \nabla F(\mathbf{S}^k + \mathbf{U}_k)\right\|_F \left\|R_\EXP(\mathbf{S}^k, - \overline{\alpha}_k\grad F(\mathbf{S}^k)) - \mathbf{S}^k\right\|_F\\
    &\geq (eL_F)^{-1} |g(\nabla F(\mathbf{S}^k) - \nabla F(\mathbf{S}^k + \mathbf{U}_k), R_\EXP(\mathbf{S}^k, - \overline{\alpha}_k\grad F(\mathbf{S}^k)) - \mathbf{S}^k)|\\
    &= (eL_F)^{-1}|g(\nabla F(\mathbf{S}^k), R_\EXP(\mathbf{S}^k, - \overline{\alpha}_k\grad F(\mathbf{S}^k)) - \mathbf{S}^k) + \gamma \overline{\alpha}_k g(\nabla F(\mathbf{S}^k), \grad F(\mathbf{S}^k))|\\
    &\geq (1-\gamma)\overline{\alpha}_k(eL_F)^{-1} |g(\nabla F(\mathbf{S}^k),\grad F(\mathbf{S}^k))| -\\
    &\qquad (eL_F)^{-1}|g(\nabla F(\mathbf{S}^k), R_\EXP(\mathbf{S}^k, - \overline{\alpha}_k\grad F(\mathbf{S}^k)) - (\mathbf{S}^k - \overline{\alpha}_k\grad F(\mathbf{S}^k))|.
\end{align}
Combining with $g(\nabla F(\mathbf{S}^k),\grad F(\mathbf{S}^k)) = \left\|\grad F(\mathbf{S}^k)\right\|_F^2$ (proof of Proposition~\ref{prop:alpha_grad}), using Cauchy-Schwarz inequality and dividing by $\overline{\alpha}_k\left\|\grad F(\mathbf{S}^k)\right\|_F^2$, we obtain
\begin{equation}
    \overline{\alpha}_k \geq (1-\gamma)(eL_{F})^{-1} - (eL_{F})^{-1}  \frac{\left\|\nabla F(\mathbf{S})\right\|_F\left\|R_\EXP(\mathbf{S}^k, - \overline{\alpha}_k\grad F(\mathbf{S}^k)) - (\mathbf{S}^k - \overline{\alpha}_k\grad F(\mathbf{S}^k))\right\|_F}{\overline{\alpha}_k\left\|\grad F(\mathbf{S}^k)\right\|_F^2}.
\end{equation}
Finally, since $\left\|\nabla F(\mathbf{S})\right\|_F \leq 2\sqrt{md} \left\|\mathbf{C}\right\|_F$ and using Proposition~\ref{prop:second_order_boundedness_of_Rtilde},
we obtain 
\begin{equation}
    \liminf \overline{\alpha}_k \geq \frac{1-\gamma}{eL_F + 2\sqrt{md}\left\|\mathbf{C}\right\|_F} > 0.
\end{equation}


\proofof{Theorem~\ref{thm:rgd_conv2}}
Below, we provide an adaptation of the proof of [33, Chapter 7, Theorem 4.2] tailored to our setting. First, we need the following lemma which will be utilized in the application of Taylor's theorem.
\begin{lem}
\label{rev2:lem1}
    Let $\mathbf{O}_1, \mathbf{O}_2 \in \mathbb{O}(d)^m$ and $\widetilde{\mathbf{O}}_1 = \pi(\mathbf{O}_1)$ and $\widetilde{\mathbf{O}}_2 = \pi(\mathbf{O}_2)$. Let $\widetilde{\mathbf{O}} = \pi(\mathbf{O})$ be on the minimal geodesic between $\widetilde{\mathbf{O}}_1$ and $\widetilde{\mathbf{O}}_2$. Suppose $\min_{\mathbf{Q} \in \mathbb{O}(d)}\left\|\mathbf{O}_j - \mathbf{S}^*\mathbf{Q}\right\|_F < \min\{2,\frac{2}{\pi}\zeta\delta(\mathbf{S}^*)\}$ for both $j =1, 2$. Then, $\min_{\mathbf{Q} \in \mathbb{O}(d)}\left\|\mathbf{O} - \mathbf{S}^*\mathbf{Q}\right\|_F < \min\{\pi, \zeta\delta(\mathbf{S}^*)\}$.
\end{lem}
\begin{proof}
Define $d(\cdot,\cdot)$ to be the geodesic distance on $\mathbb{O}(d)$ equipped with the Euclidean metric induced by the ambient space i.e. $\mathbb{R}^{d \times d}$. Then, from the properties of a quotient and a product manifold, for a general $\widetilde{\mathbf{S}} = \pi(\mathbf{S})$ it trivially follows that 
\begin{align}
 d_{\widetilde{g}}(\widetilde{\mathbf{S}},\widetilde{\mathbf{S}}^{*}) &= \min_{\mathbf{Q}\in \mathbb{O}(d)} d_{g}(\mathbf{S},\mathbf{S}^{*}\mathbf{Q}) \geq \min_{\mathbf{Q}\in \mathbb{O}(d)} \left\|\mathbf{S} -\mathbf{S}^{*}\mathbf{Q}\right\|_F. \label{rev2:eq7}
\end{align}
Moreover, if $\min_{\mathbf{Q} \in \mathbb{O}(d)}\left\|\mathbf{S} - \mathbf{S}^*\mathbf{Q}\right\|_F < 2$ (meaning $\mathbf{S}_i$ and $\mathbf{S}^*_i\mathbf{Q}^*$ where $\mathbf{Q}^*$ is the minimizer, lie in the same component of $\mathbb{O}(d)$ for all $i \in [1,m]$), then
\begin{align}
 d_{\widetilde{g}}(\widetilde{\mathbf{S}},\widetilde{\mathbf{S}}^{*})^2 &= \min_{\mathbf{Q}\in \mathbb{O}(d)} d_{g}(\mathbf{S},\mathbf{S}^{*}\mathbf{Q})^2 \leq d_{g}(\mathbf{S},\mathbf{S}^{*}\mathbf{Q}^*)^2 = \sum_{i=1}^{m} d(\mathbf{S}_i,\mathbf{S}^{*}_i\mathbf{Q}^*)^2\\
 &\leq  \frac{\pi^2}{4} \sum_{i=1}^{m}\left\|\mathbf{S}_i -\mathbf{S}^{*}_i\mathbf{Q}^*\right\|_F^2 = \frac{\pi^2}{4} \left\|\mathbf{S} -\mathbf{S}^{*}\mathbf{Q}^*\right\|_F^2.\\
 &= \min_{\mathbf{Q}\in \mathbb{O}(d)}\frac{\pi^2}{4} \left\|\mathbf{S} -\mathbf{S}^{*}\mathbf{Q}\right\|_F^2.
\end{align}
Here, the second inequality follows from the following fact: if $\mathbf{Q}_1, \mathbf{Q}_2 \in \mathbb{O}(d)$ and are in the same component, then $d(\mathbf{Q}_1,\mathbf{Q}_2)^2 \leq \frac{\pi^2}{4} \left\|\mathbf{Q}_1-\mathbf{Q}_2\right\|_F^2$. This is because, since $\mathbf{Q}_1^T\mathbf{Q}_2 \in \mathbb{SO}(d)$ and therefore we can write $\mathbf{Q}_1^T\mathbf{Q}_2 = \mathbf{V}\mathbf{R}\mathbf{V}^T$ where $\mathbf{V} \in \mathbb{O}(d)$,
$\mathbf{R} = \blockdiag(R(\theta_1),\ldots,R(\theta_k),\mathbf{I}_r)$
and $R(\theta_j) = \begin{bmatrix}\cos(\theta_j) & -\sin(\theta_j)\\ \sin(\theta_j) & \cos(\theta_j)\end{bmatrix}$. Here $2k + r = d$ and $\theta_j \in (-\pi,\pi]$. Also, $\bm{\Omega} = \text{Log}(\mathbf{Q}_1^T\mathbf{Q}_2)$ has nonzero eigenvalues $\{\pm \iota \theta_j\}_{j=1}^{k}$. Putting all together,
\begin{align}
 d(\mathbf{Q}_1,\mathbf{Q}_2)^2 = \left\|\bm{\Omega}\right\|_F^2 = 2\sum_{j=1}^{k}\theta_j^2 \leq \frac{\pi^2}{4} \sum_{j=1}^{k} 4(1-\cos(\theta_j)) = \frac{\pi^2}{4} \left\|\mathbf{I}_d - \mathbf{Q}_1^T\mathbf{Q}_2\right\|_F^2 = \frac{\pi^2}{4} \left\|\mathbf{Q}_1-\mathbf{Q}_2\right\|_F^2.
\end{align}
It follows that $d_{\widetilde{g}}(\widetilde{\mathbf{O}}_j,\widetilde{\mathbf{S}}^{*}) < \min\{\pi, \zeta \delta(\mathbf{S}^*)\}$ for both $j=1,2$. Then, since a geodesic ball of radius less then $\pi$ is geodesically convex in $\mathbb{O}(d)$, the same result holds for geodesic balls on $\mathbb{O}(d)^m$. Consequently, $d_{\widetilde{g}}(\widetilde{\mathbf{O}},\widetilde{\mathbf{S}}^{*}) < \min\{\pi, \zeta \delta(\mathbf{S}^*)\}$. From Eq.~(\ref{rev2:eq7}), it follows that $\min_{\mathbf{Q}\in \mathbb{O}(d)} \left\|\mathbf{O} -\mathbf{S}^{*}\mathbf{Q}\right\|_F < \min\{\pi,\zeta \delta(\mathbf{S}^*)\}$.
\end{proof}

Finally, the proof [33, Chapter 7, Theorem 4.2] adapted to our setting follows.
For convenience, denote $a = (1-\zeta)\lambda_{-}(\mathbf{S}^*)$ and $b = (\lambda_{+}(\mathbf{S}^*)+\zeta \lambda_{-}(\mathbf{S}^*))$.
Due to the assumption on $\mathbf{S}^*$ and Proposition~\ref{prop:HessVicinity}, 
\begin{equation}
    a\widetilde{g}(\widetilde{\mathbf{Z}}, \widetilde{\mathbf{Z}}) \leq \widetilde{g}(\Hess \widetilde{F}(\widetilde{\mathbf{O}})[\widetilde{\mathbf{Z}}],\widetilde{\mathbf{Z}}) \leq  b\widetilde{g}(\widetilde{\mathbf{Z}}, \widetilde{\mathbf{Z}})
\end{equation}
for all $\mathbf{O} \in \mathbb{O}(d)^m$, $\widetilde{\mathbf{O}} = \pi(\mathbf{O})$ and $\widetilde{\mathbf{Z}} \in T_{\widetilde{\mathbf{O}}}\mathbb{O}(d)^m/_{\sim}$ such that $\min_{\mathbf{Q}\in\mathbb{O}(d)}\left\|\mathbf{O}-\mathbf{S}^*\mathbf{Q}\right\|_F < \min\{2, \frac{2}{\pi} \zeta\delta(\mathbf{S}^*)\} < \zeta\delta(\mathbf{S}^*)$. Let $\widetilde{\mathbf{Z}} \in T_{\widetilde{\mathbf{O}}}\mathbb{O}(d)^m/_{\sim}$ be such that $\widetilde{\mathbf{S}}^* = \widetilde{R}_{\EXP }(\widetilde{\mathbf{O}}, \widetilde{\mathbf{Z}})$. Then, using the Taylor's theorem centered at $\widetilde{\mathbf{O}}$,
\begin{align}
    \widetilde{F}(\widetilde{\mathbf{S}}^*) &= \widetilde{F}(\widetilde{\mathbf{O}}) + \widetilde{g}(\grad \widetilde{F}(\widetilde{\mathbf{O}}), \widetilde{\mathbf{Z}}) +\frac{1}{2} \widetilde{g}(\Hess \widetilde{F}(\widetilde{\mathbf{O}}')[\widetilde{\mathbf{Z}}'],\widetilde{\mathbf{Z}}')
\end{align}
where $\widetilde{\mathbf{O}}'$ is a point on the minimial geodesic between $\widetilde{\mathbf{O}}$ and $\widetilde{\mathbf{S}}^*$ and $\widetilde{\mathbf{Z}}'$ is the parallel transport of $\widetilde{\mathbf{Z}}$ from $\widetilde{\mathbf{O}}$ to $\widetilde{\mathbf{O}}'$. Consequently, using Lemma~\ref{rev2:lem1} (here $\left\|\cdot\right\|$ is the norm induced by $\widetilde{g}$),
\begin{align}
    \widetilde{F}(\widetilde{\mathbf{O}}) - \widetilde{F}(\widetilde{\mathbf{S}}^*) &= -\widetilde{g}(\grad \widetilde{F}(\widetilde{\mathbf{O}}), \widetilde{\mathbf{Z}}) - \frac{1}{2}\widetilde{g}(\Hess \widetilde{F}(\widetilde{\mathbf{O}}')[\widetilde{\mathbf{Z}}'],\widetilde{\mathbf{Z}}')\\
    &\leq \left\|\grad \widetilde{F}(\widetilde{\mathbf{O}})\right\|\left\|\widetilde{\mathbf{Z}}\right\| - \frac{1}{2}a \left\|\widetilde{\mathbf{Z}}'\right\|^2\\
    &= \left\|\grad \widetilde{F}(\widetilde{\mathbf{O}})\right\|\left\|\widetilde{\mathbf{Z}}\right\| - \frac{1}{2}a \left\|\widetilde{\mathbf{Z}}\right\|^2\\
    &= \left\|\grad \widetilde{F}(\widetilde{\mathbf{O}})\right\|d_{\widetilde{g}}(\widetilde{\mathbf{O}},\widetilde{\mathbf{S}}^*) - \frac{1}{2}a d_{\widetilde{g}}(\widetilde{\mathbf{O}},\widetilde{\mathbf{S}}^*)^2. \label{rev2:eq1}
\end{align}
Similarly, since $\grad \widetilde{F}(\widetilde{\mathbf{S}}^*) = 0$,
\begin{align}
    \frac{1}{2}a d_{\widetilde{g}}(\widetilde{\mathbf{O}},\widetilde{\mathbf{S}}^*)^2 \leq \widetilde{F}(\widetilde{\mathbf{O}}) - \widetilde{F}(\widetilde{\mathbf{S}}^*) \leq \frac{1}{2}b d_{\widetilde{g}}(\widetilde{\mathbf{O}},\widetilde{\mathbf{S}}^*)^2. \label{rev2:eq2}
\end{align}
From the above equation and Eq.~(\ref{rev2:eq1}), we obtain
\begin{equation}
    \frac{1}{2}a d_{\widetilde{g}}(\widetilde{\mathbf{O}},\widetilde{\mathbf{S}}^*)^2 \leq \left\|\grad \widetilde{F}(\widetilde{\mathbf{O}})\right\|d_{\widetilde{g}}(\widetilde{\mathbf{O}},\widetilde{\mathbf{S}}^*) - \frac{1}{2}a d_{\widetilde{g}}(\widetilde{\mathbf{O}},\widetilde{\mathbf{S}}^*)^2. 
\end{equation}
This implies
\begin{equation}
    d_{\widetilde{g}}(\widetilde{\mathbf{O}},\widetilde{\mathbf{S}}^*) \leq \frac{\left\|\grad \widetilde{F}(\widetilde{\mathbf{O}})\right\|}{a}.
\end{equation}
From Eq.~(\ref{rev2:eq2}) we also obtain
\begin{equation}
    d_{\widetilde{g}}(\widetilde{\mathbf{O}},\widetilde{\mathbf{S}}^*)^2 \geq \frac{2}{b}(\widetilde{F}(\widetilde{\mathbf{O}}) - \widetilde{F}(\widetilde{\mathbf{S}}^*)).
\end{equation}
The above two equations, combined with Eq.~(\ref{rev2:eq1}) results in
\begin{equation}
    \widetilde{F}(\widetilde{\mathbf{O}}) - \widetilde{F}(\widetilde{\mathbf{S}}^*) \leq \frac{\left\|\grad \widetilde{F}(\widetilde{\mathbf{O}})\right\|^2}{a} - \frac{a}{b}\left(\widetilde{F}(\widetilde{\mathbf{O}}) - \widetilde{F}(\widetilde{\mathbf{S}}^*)\right). \label{rev2:eq}
\end{equation}
This simplifies to
\begin{equation}
    \left\|\grad \widetilde{F}(\widetilde{\mathbf{O}})\right\|^2 \geq a \left(1+\frac{a}{b}\right)\left(\widetilde{F}(\widetilde{\mathbf{O}}) - \widetilde{F}(\widetilde{\mathbf{S}}^*)\right).
\end{equation}
Moreover, following the proof of Proposition~\ref{prop:alpha_grad} and using the fact that $\left\|\grad \widetilde{F}(\widetilde{\mathbf{S}}^k)\right\| = \left\|\grad F (\mathbf{S}^k)\right\|_F$ (see Proposition~\ref{prop:g_tilde}), it is easy to deduce that
\begin{equation}
    \widetilde{F}(\widetilde{\mathbf{S}}^{k+1})-\widetilde{F}(\widetilde{\mathbf{S}}^{k}) \leq -\gamma \alpha_k \left\|\grad \widetilde{F}(\widetilde{\mathbf{S}}^{k})\right\|^2. \label{rev2:eq3}
\end{equation}

To this end we note that $\min_{\mathbf{Q} \in \mathbb{O}(d)}\left\|\mathbf{S}^{k} -\mathbf{S}^{*}\mathbf{Q}\right\|_F < \min\left\{2, \frac{2}{\pi}\zeta \delta(\mathbf{S}^*)\right\}$ follows from the assumption of the theorem. As a consequence, Eq.~(\ref{rev2:eq}) holds at $\widetilde{\mathbf{S}}^{k}$ for all $k$. Combining Eq.~(\ref{rev2:eq}) and Eq.~(\ref{rev2:eq3}),
\begin{equation}
    \widetilde{F}(\widetilde{\mathbf{S}}^{k+1})-\widetilde{F}(\widetilde{\mathbf{S}}^{k}) \leq -\gamma \alpha_k a \left(1+\frac{a}{b}\right)\left(\widetilde{F}(\widetilde{\mathbf{S}}^k) - \widetilde{F}(\widetilde{\mathbf{S}}^*)\right). \label{rev2:eq4}
\end{equation}
Adding and subtracting $\widetilde{F}(\widetilde{\mathbf{S}}^*)$ on the l.h.s of Eq.~(\ref{rev2:eq4}), we obtain
\begin{align}
     \widetilde{F}(\widetilde{\mathbf{S}}^{k+1})-\widetilde{F}(\widetilde{\mathbf{S}}^*) &\leq \left(1-\gamma \alpha_k a \left(1+\frac{a}{b}\right)\right)\left(\widetilde{F}(\widetilde{\mathbf{S}}^k) - \widetilde{F}(\widetilde{\mathbf{S}}^*)\right)\\
     &\leq q_k(\mathbf{S}^*, \gamma) \left(\widetilde{F}(\widetilde{\mathbf{S}}^k) - \widetilde{F}(\widetilde{\mathbf{S}}^*)\right)\\
     &\leq \prod_{j=0}^{k}q_j\left(\widetilde{F}(\widetilde{\mathbf{S}}^0) - \widetilde{F}(\widetilde{\mathbf{S}}^*)\right)
\end{align}
where $q_j < 1$ for all $k \in \{0,\ldots,k\}$. Combining the above inequality with Eq.~(\ref{rev2:eq2}),
\begin{align}
    d_{\widetilde{g}}(\widetilde{\mathbf{S}}^{k+1}, \widetilde{\mathbf{S}}^*) &\leq \left(\frac{2}{a}\right)^{1/2} \left(\widetilde{F}(\widetilde{\mathbf{S}}^{k+1})-\widetilde{F}(\widetilde{\mathbf{S}}^*)\right)^{1/2}\\
    &\leq \left(\frac{2}{a}\right)^{1/2}\left(\prod_{j=0}^{k}q_j^{1/2}\right) \left(\widetilde{F}(\widetilde{\mathbf{S}}^0)-\widetilde{F}(\widetilde{\mathbf{S}}^*)\right)^{1/2}
\end{align}
Finally,  we derive a common upper bound for $q_j$. By Taylor's theorem, there exist $\widetilde{\mathbf{O}}$ on the minimal geodesic between $\widetilde{\mathbf{S}}^{k+1}$ and $\widetilde{\mathbf{S}}^{k}$ such that $\mathbf{\widetilde{Z}}$ is the parallel transport of $-\grad \widetilde{F}(\widetilde{\mathbf{S}}^k)$ to $\widetilde{\mathbf{O}}$. Then, it follows from Lemma~\ref{rev2:lem1},
\begin{align}
    \widetilde{F}(\widetilde{\mathbf{S}}^{k+1}) - \widetilde{F}(\widetilde{\mathbf{S}}^k) &=  -\alpha_k\left\|\grad \widetilde{F}(\widetilde{\mathbf{S}}^k)\right\|^2 +\frac{\alpha_k^2}{2} \widetilde{g}(\Hess \widetilde{F}(\widetilde{\mathbf{O}})[\widetilde{\mathbf{Z}}], \widetilde{\mathbf{Z}})\\
    &\leq -\alpha_k \left\|\grad \widetilde{F}(\widetilde{\mathbf{S}}^k)\right\|^2 + \frac{\alpha_k^2b}{2}\left\|\widetilde{\mathbf{Z}}\right\|^2\\
    &= -\alpha_k \left\|\grad \widetilde{F}(\widetilde{\mathbf{S}}^k)\right\|^2 + \frac{\alpha_k^2b}{2}\left\|\grad \widetilde{F}(\widetilde{\mathbf{S}}^k)\right\|^2\\
    &= -\alpha_k\left(1-\frac{\alpha_kb}{2}\right)\left\|\grad \widetilde{F}(\widetilde{\mathbf{S}}^k)\right\|^2
\end{align}
Therefore, Eq.~(\ref{rev2:eq3}) is satisfied if
\begin{equation}
    1 - \frac{\alpha_k b}{2} \geq \gamma \implies \alpha_k \leq \alpha^* = \frac{2(1-\gamma)}{b}.
\end{equation}
Consequently,
\begin{align}
    q_j \leq q &= 1 - \gamma \alpha^* a \left(1 + \frac{a}{b}\right)\\
    &= 1 - 2\frac{\gamma(1-\gamma)a}{b}\left(1 + \frac{a}{b}\right)\\
    &= 1-2\gamma(1-\gamma)r(1+r).
\end{align}

\proofof{Lemma~\ref{lem:quadgrowth}}
Since $\mathbf{S}_0$ is a unique perfect alignment, all other perfect alignments are of the form $\mathbf{S}_0\mathbf{Q}$ where $\mathbf{Q}\in \mathbb{O}(d)$. Moreover, the null space of $\mathbf{C}_0$ is exactly the span of columns of $\mathbf{S}_0$. Therefore, we obtain the following decomposition of $\mathbf{C}_0 = \mathbf{U}_0\boldsymbol{\Lambda}_0 \mathbf{U}_0^T$ where $\mathbf{U}_0^T\mathbf{U}_0 = \mathbf{I}_{(m-1)d}$, $\mathbf{S}_0^T \mathbf{U}_0 = 0$ and $\boldsymbol{\Lambda}_0$ is a diagonal matrix containing the strictly positive eigenvalues of $\mathbf{C}_0$. Using the above decomposition, we have 
$$\Tr(\mathbf{C}_0\mathbf{S}\mathbf{S}^T) = \Tr(\mathbf{U}_0\boldsymbol{\Lambda}_0 \mathbf{U}_0^T\mathbf{S}\mathbf{S}^T) \geq \lambda_{d+1}(\mathbf{C}_0) \left\|\mathbf{U}_0^T\mathbf{S}\right\|_F^2.$$

\textbf{Claim:} $\left\|\mathbf{U}_0^T\mathbf{S}\right\|_F^2 \geq \frac{1}{2}\min_{\mathbf{Q} \in \mathbb{O}(d)}\left\|\mathbf{S}- \mathbf{S}_0\mathbf{Q}\right\|_F^2$. Since the union of the columns of $\mathbf{S}_0$ and $\mathbf{U}_0$ span $\mathbb{R}^{md}$ therefore there exist $\mathbf{R}_1 \in \mathbb{R}^{d \times d}$ and $\mathbf{R}_2 \in \mathbb{R}^{(m-1)d \times d}$ such that $\mathbf{S} = \mathbf{S}_0 \mathbf{R}_1 + \mathbf{U}_0 \mathbf{R}_2$. Let $\mathbf{Q}^* \in \mathbb{O}(d)$ be such that $\left\|\mathbf{S}- \mathbf{S}_0\mathbf{Q}^*\right\|_F^2 = \min_{\mathbf{Q} \in \mathbb{O}(d)^m} \left\|\mathbf{S}- \mathbf{S}_0\mathbf{Q}\right\|_F^2$. Then $\mathbf{Q}^* = \mathbf{U}_1\mathbf{V}_1^T$ where $\mathbf{R}_1 = \mathbf{U}_1\boldsymbol{\Sigma}_1\mathbf{V}_1^T$ is a singular vector decomposition of $\mathbf{R}_1$. Using $\mathbf{S}^T\mathbf{S} = \mathbf{S}_0^T\mathbf{S}_0 = m\mathbf{I}_d$ and $\mathbf{R}_1 = \mathbf{U}_1\boldsymbol{\Sigma}_1\mathbf{V}_1^T$,

\begin{align}
    \left\|\mathbf{U}_0^T\mathbf{S}\right\|_F^2 = \left\|\mathbf{R}_2\right\|_F^2 = \left\|\mathbf{S}-\mathbf{S}_0\mathbf{R}_1\right\|_F^2 &= md + m \left\|\mathbf{R}_1\right\|_F^2 - 2\Tr(\mathbf{S}\mathbf{R}_1^T\mathbf{S}_0^T)\\
    &= md + m\left\|\boldsymbol{\Sigma}\right\|_F^2 - 2\Tr(\mathbf{S}\mathbf{V}_1\boldsymbol{\Sigma}_1\mathbf{U}_1^T\mathbf{S}_0^T)
\end{align}

Moreover, using the fact that $\mathbf{S}_0^T\mathbf{S}_0 = m\mathbf{I}_d$ and the definition of $\mathbf{Q}^*$,
{\allowdisplaybreaks
\begin{align}
    \left\|\mathbf{S} - \mathbf{S}_0\mathbf{Q}^*\right\|_F^2 &= \left\|\mathbf{U}_0\mathbf{R}_2\right\|_F^2 + \left\|\mathbf{S}_0(\mathbf{R}_1 - \mathbf{Q}^*)\right\|_F^2 = \left\|\mathbf{S} - \mathbf{S}_0\mathbf{R}_1\right\|_F^2 + m\left\|\mathbf{R}_1 - \mathbf{Q}^*\right\|_F^2\\
    &= \left\|\mathbf{S} - \mathbf{S}_0\mathbf{R}_1\right\|_F^2 + m\left\|\mathbf{I}_d-\boldsymbol{\Sigma}_1\right\|_F^2\\
    &= \left\|\mathbf{S} - \mathbf{S}_0\mathbf{R}_1\right\|_F^2 + m(d + \left\|\boldsymbol{\Sigma}_1\right\|_F^2 - 2\Tr(\boldsymbol{\Sigma}_1))\\
    &= 2(md + m\left\|\boldsymbol{\Sigma}_1\right\|_F^2) - 2(\Tr(\mathbf{S}\mathbf{V}_1\boldsymbol{\Sigma}_1\mathbf{U}_1^T\mathbf{S}_0^T) + m\Tr(\boldsymbol{\Sigma}_1))
\end{align}
}
Overall, 
$$\left\|\mathbf{U}_0^T\mathbf{S}\right\|_F^2 - \frac{1}{2}\left\|\mathbf{S} - \mathbf{S}_0\mathbf{Q}^*\right\|_F^2 = m\Tr(\boldsymbol{\Sigma}_1) - \Tr(\mathbf{S}\mathbf{V}_1\boldsymbol{\Sigma}_1\mathbf{U}_1^T\mathbf{S}_0^T) = \sum_{1}^{m}(\Tr(\boldsymbol{\Sigma}_1) - \Tr(\mathbf{U}_1^T\mathbf{S}_{0_i}^T\mathbf{S}_i\mathbf{V}_1\boldsymbol{\Sigma}_1)).$$
Since $\max_{\mathbf{Q} \in \mathbb{O}(d)}\Tr(\mathbf{Q} \boldsymbol{\Sigma}_1) = \Tr(\boldsymbol{\Sigma}_1)$, the result follows.

\proofof{Lemma~\ref{lem:distS_0Sstar}}
The proof is motivated from \citeb[Proposition~4.32]{bonnans2013perturbation}. Define 
$$H(\mathbf{S}) = \Tr(\mathbf{C}\mathbf{S}\mathbf{S}^T) - \Tr(\mathbf{C}_0\mathbf{S}\mathbf{S}^T)$$
and note that $H(\mathbf{S})$ is Lipschitz with a Lipschitz constant bounded by $\left\|\nabla H(\mathbf{S})\right\|_F \leq 2m \left\|\mathbf{C}-\mathbf{C}_0\right\|_F$. 

Let $\mathbf{Q}^* \in \mathbb{O}(d)$ be such that 
$$\left\|\mathbf{S}^* - \mathbf{S}_0\mathbf{Q}^*\right\|_F = \min_{\mathbf{Q}\in\mathbb{O}(d)}\left\|\mathbf{S}^* - \mathbf{S}_0\mathbf{Q}\right\|_F.$$
Then, using the mean value theorem and the fact the $\mathbf{S}^*$ is an optimal alignment in the noisy setting, meaning $\Tr(\mathbf{C}\mathbf{S}^*\mathbf{S}^{*^T}) \leq \Tr(\mathbf{C}\mathbf{S}_0\mathbf{Q}^*(\mathbf{S}_0\mathbf{Q}^*)^T)$,
\begin{align}
    \Tr(\mathbf{C}_0\mathbf{S}^*\mathbf{S}^{*^T}) - \Tr(\mathbf{C}_0\mathbf{S}_0\mathbf{Q}^*(\mathbf{S}_0\mathbf{Q}^*)^T) &= H(\mathbf{S}_0\mathbf{Q}^*) - H(\mathbf{S}^*) + (\Tr(\mathbf{C}\mathbf{S}^*\mathbf{S}^{*^T}) - \Tr(\mathbf{C}\mathbf{S}_0\mathbf{Q}^*(\mathbf{S}_0\mathbf{Q}^*)^T))\\
    &\leq 2m \left\|\mathbf{C}-\mathbf{C}_0\right\|_F \left\|\mathbf{S}^*-\mathbf{S}_0\mathbf{Q}^*\right\|_F
\end{align}
Combining with Lemma~\ref{lem:quadgrowth} and $\mathbf{C}_0\mathbf{S}_0 = 0$, we obtain 
$$(\lambda_{d+1}(\mathbf{C}_0)/2) \left\|\mathbf{S}^*- \mathbf{S}_0\mathbf{Q}^*\right\|_F^2 \leq 2m \left\|\mathbf{C}-\mathbf{C}_0\right\|_F \left\|\mathbf{S}^*-\mathbf{S}_0\mathbf{Q}^*\right\|_F.$$
The result follows.

\proofof{Theorem~\ref{thm:rgd_noise_stability}} \revadd{The following bound holds from \citeb[Eq.~(5.8, 5.11, 5.12)]{chaudhury2015global}
\begin{equation}
    \min_{\mathbf{Q} \in \mathbb{O}(d)} \left\|\mathbf{S}_{spec}(\mathbf{C}) - \mathbf{S}_0\mathbf{Q}\right\|_F \leq \frac{4\pi\sqrt{md(d+1)}}{\lambda_{d+1}(\mathbf{C})}(K_1 \varepsilon + K_2\varepsilon^2). \label{eq:spec_bound}
\end{equation}
Let $\mathbf{Q}^*$ be such that $\left\|\mathbf{S}^* - \mathbf{S}_0\mathbf{Q}^*\right\|_F = \min_{\mathbf{Q}\in\mathbb{O}(d)}\left\|\mathbf{S}^* - \mathbf{S}_0\mathbf{Q}\right\|_F$. Then,
\begin{align}
    \left\|\mathbf{S}_{spec}(\mathbf{C}) - \mathbf{S}^*\mathbf{Q}\right\|_F &\leq \left\|\mathbf{S}_{spec}(\mathbf{C}) - \mathbf{S}_0\mathbf{Q}^*\mathbf{Q}\right\|_F + \left\|\mathbf{S}^*\mathbf{Q} - \mathbf{S}_0\mathbf{Q}^*\mathbf{Q}\right\|_F\\
    &\leq \left\|\mathbf{S}_{spec}(\mathbf{C}) - \mathbf{S}_0\mathbf{Q}^*\mathbf{Q}\right\|_F + \left\|\mathbf{S}^* - \mathbf{S}_0\mathbf{Q}^*\right\|_F
\end{align}
Minimizing the above over $\mathbf{Q}$, using the fact that 
$$\min_{\mathbf{Q}\in\mathbb{O}(d)}\left\|\mathbf{S}_{spec}(\mathbf{C}) - \mathbf{S}_0\mathbf{Q}^*\mathbf{Q}\right\|_F = \min_{\mathbf{Q}\in\mathbb{O}(d)}\left\|\mathbf{S}_{spec}(\mathbf{C}) - \mathbf{S}_0\mathbf{Q}\right\|_F,$$
followed by Eq.~(\ref{eq:spec_bound}), Lemma~\ref{lem:distS_0Sstar}, Theorem~\ref{thm:rgd_conv2} and \citeb[Eq.~(5.12)]{chaudhury2015global}, the result follows.}
\bibliographystyle{siamplain}
\bibliography{loc_rigid}

\end{document}